\documentclass{amsart}
\pdfoutput=1
\usepackage{fullpage}
\usepackage{natbib}
\usepackage{amssymb,amsmath}
\usepackage{color}
\usepackage{marginnote} % for \marinate
\usepackage{lineno}
\usepackage{graphicx}
\newcommand{\bone}{\mathbf{1}}

\newenvironment{biology}[1]{\vskip 12pt\textbf{Biological interpretation#1.} \it}{\vskip 12pt}
\newtheorem{theorem}{Theorem}[section]
\newtheorem{lemma}[theorem]{Lemma}
\newtheorem{corollary}[theorem]{Corollary}
\newtheorem{proposition}[theorem]{Proposition}

 \newtheoremstyle{example}{5pt}{5pt}%
     {}%         Body font
     {}%         Indent amount (empty = no indent, \parindent = para indent)
     {\bfseries}% Thm head font
     {.}%        Punctuation after thm head
     {5pt}%     Space after thm head (\newline = linebreak)
     {\thmname{#1}\thmnumber{ #2}\thmnote{ #3}}%         Thm head spec

   \theoremstyle{example}
   \newtheorem{example}{Example}[section]

\newcommand{\Tr}{\mathop{\mbox{Tr}}}

\def\diag{\mathrm{diag}}
\def\id{\mathbf{id}}
\def\trchar{\kappa_{\mbox{\rm tr}}} % trivial character
\def\gtimes{} % group multiplication
\newcommand{\bV}{\mathbf{V}}
\newcommand{\bE}{\mathbf{E}}
\newcommand{\bX}{\mathbf{X}}
\newcommand{\bY}{\mathbf{Y}}
\newcommand{\bZ}{\mathbf{Z}}

\newcommand{\bB}{\mathbf{B}}
\newcommand{\bH}{\mathbf{H}}
\newcommand{\bU}{ \mathbf{U} }

\newcommand{\bmu}{\mathbf{\mu} }
\newcommand{\bx}{\mathbf{x}}

\newcommand{\bM}{\mathbf{M}}

\newcommand{\Var}{\mathrm{Var}}
\newcommand{\E}{\mathbb{E}}
\newcommand{\R}{\mathbb{R}}

\def\Z{\mathbb{Z}}
\def\P{\mathbb{P}}

\begin{document}

\title{Stochastic population growth in spatially heterogeneous environments}
% \author{Steven N. Evans \and Peter L. Ralph \and Sebastian J. Schreiber \and Arnab Sen}
% \institute{S.N. Evans \at Department of Statistics \#3860, 367 Evans Hall, University of California, Berkeley, CA  94720-3860 USA \email{evans@stat.berkeley.edu}
% \and 
% P.L. Ralph \and S.J. Schreiber \at 
% Department of Evolution and Ecology, University of California,  Davis, CA 956116  USA 
% \email{plralph@ucdavis.edu, sschreiber@ucdavis.edu}
% \and
% A. Sen \at  Statistical Laboratory, Centre for Mathematical Sciences,
% Wilberforce Road, Cambridge CB3 0WB, United Kingdom \email{A.Sen@statslab.cam.ac.uk}
% }

\author{Steven N. Evans}
\address{Department of Statistics \#3860 \\
367 Evans Hall \\
University of California \\
Berkeley, CA  94720-3860 \\
USA}
\email{evans@stat.berkeley.edu}
\thanks{SNE was supported in part by NSF grant DMS-0907630}

\author{Peter L. Ralph}
\address{Department of Evolution and Ecology \\
University of California \\
Davis, CA 956116 \\
USA}
\email{plralph@ucdavis.edu}
\thanks{PLR was supported by funds to SJS from the Dean's Office of the College of Biological Sciences, University of California, Davis, and by NIH fellowship F32-GM096686}

\author{Sebastian J. Schreiber}
\address{Department of Evolution and Ecology \\
University of California \\
Davis, CA 956116 \\
USA}
\email{sschreiber@ucdavis.edu}
\thanks{SJS was supported in part by NSF grants EF-0928987 and DMS-1022639}

\author{Arnab Sen}
\address{Statistical Laboratory \\
Centre for Mathematical Sciences \\
Wilberforce Road \\
Cambridge CB3 0WB \\
United Kingdom}
\email{A.Sen@statslab.cam.ac.uk}
\thanks{AS was supported by EPSRC grant EP/G055068/1}

\bibliographystyle{plainnat}

\maketitle

\begin{abstract}
Classical ecological theory predicts that environmental stochasticity increases extinction risk by reducing the average per-capita growth rate of populations. For sedentary populations in a spatially homogeneous yet temporally variable environment, a simple model of population growth is  
a stochastic differential equation $dZ_t=\mu Z_t dt + \sigma Z_t dW_t$, $t \ge 0$,
where the conditional law of $Z_{t + \Delta t}- Z_t$
given $Z_t = z$ has mean and variance approximately $z \mu \Delta t$ and $z^2 \sigma^2 \Delta t$ when the time increment
$\Delta t$ is small.
The long-term stochastic  growth rate  
$\lim_{t \to \infty} t^{-1} \log Z_t$ for such a population equals $\mu -\frac{\sigma^2}{2}$. Most populations, however, experience spatial as well as temporal variability. To understand the interactive effects of environmental stochasticity, spatial heterogeneity, and dispersal on population growth, we study an analogous model $\bX_t = (X_t^1, \ldots, X_t^n)$, $t \ge 0$, for the 
population abundances in $n$ patches:  the conditional law of
$\bX_{t+\Delta t}$ given $\bX_t = x$ is such that 
 the conditional mean of
$X_{t+\Delta t}^i - X_t^i$ is approximately 
$[x^i \mu_i + \sum_j (x^j D_{ji} - x^i D_{ij})] \Delta t$
where $\mu_i$ is the per capita growth rate in the $i$-th patch
and $D_{ij}$ is the dispersal rate from the $i$-th patch to the
$j$-th patch, and the conditional covariance of
$X_{t+\Delta t}^i - X_t^i$ and $X_{t+\Delta t}^j - X_t^j$
is approximately $x^i x^j \sigma_{ij} \Delta t$ for
some covariance matrix $\Sigma = (\sigma_{ij})$.
We show  for such a 
spatially extended population that if $S_t = X_t^1 + \cdots + X_t^n$ denotes
 the total population abundance, then $\bY_t = \bX_t/S_t$, the
vector of patch proportions, converges in law to
a random vector $\bY_\infty$ as $t \to \infty$, and the stochastic growth rate
$\lim_{t \to \infty} t^{-1} \log S_t$ equals the space-time average per-capita growth rate $\sum_i \mu_i \E[Y_\infty^i]$ experienced by the population minus half of the space-time average temporal variation $\E[\sum_{i,j} \sigma_{ij} Y_\infty^i Y_\infty^j]$ experienced by the population.  Using this characterization of the stochastic growth rate, we derive an explicit expression for the stochastic growth rate for populations living in two patches, determine which choices of the dispersal matrix $D$ produce the maximal stochastic growth rate for a freely dispersing population, derive an analytic approximation of the stochastic growth rate for dispersal limited populations, and use group theoretic techniques to approximate the stochastic growth rate for populations living in multi-scale landscapes (e.g. insects on plants in meadows on islands). Our results provide fundamental insights into ``ideal free'' movement in the face of uncertainty,  the persistence of coupled sink populations, the evolution of dispersal rates, and the single large or several small (SLOSS) debate in conservation biology.  For example, our analysis implies that even in the absence of density-dependent feedbacks, ideal-free dispersers occupy multiple patches in spatially heterogeneous environments provided environmental fluctuations are sufficiently strong and sufficiently weakly correlated across space. In contrast, for diffusively dispersing populations living in similar environments, intermediate dispersal rates maximize their stochastic growth rate.  
\end{abstract}


\section{Introduction}

\linenumbers

Environmental conditions (e.g. light, precipitation, nutrient availability) vary in space and time. Since these conditions influence survivorship and fecundity of an organism, all organisms whether they be plants, animals, or viruses are faced with a fundamental quandary of ``Should I stay or should I go?'' On the one hand, if individuals disperse in a spatially heterogeneous environment, then they may  arrive in locations with poorer environmental conditions. On the other hand, if individuals do not disperse, then they may fare poorly due to temporal fluctuations in local environmental conditions. The consequences of this
 interaction between dispersal and environmental heterogeneity for population growth has been studied extensively from theoretical, experimental, and applied perspectives~\citep{hastings-83,petchy-etal-97,lundberg-etal-00,gonzalez-holt-02,schmidt-04,roy-etal-05,boyce-etal-06,matthews-gonzalez-07,prsb-10,durrett-remenik-11}. Here, we provide a mathematically rigorous  perspective on these interactive effects using spatially explicit models of stochastic population growth. 

Population growth is inherently stochastic due to numerous unpredictable causes. For a single, unstructured population with overlapping generations, the simplest model accounting for these fluctuations is a linear stochastic differential equation of the form
\begin{equation}\label{eq:intro}
d Z_t= \mu Z_t dt+\sigma Z_t dB_t,
\end{equation}
where $Z_t$ is the population abundance at time $t$, $\mu$ is the mean per-capita growth rate (that is, $\E[Z_{t + \Delta t} - Z_t \, | \, Z_t = z]
\approx z \mu \Delta t$), $\sigma^2$ is the ``infinitesimal'' variance of fluctuations in the per-capita growth rate (that is, $\E[(Z_{t + \Delta t} - Z_t - z \mu \Delta t)^2\, | \, Z_t = z]
\approx z^2 \sigma^2  \Delta t$), and $B_t$ is a standard Brownian motion. 
Equivalently, the log population abundance $\log Z_t$ is normally distributed with mean $\log Z_0 + (\mu-\sigma^2/2)t$ and variance $\sigma^2 t$.   Hence, even if
the mean per-capita growth rate $\mu$ is positive
these populations decline exponentially towards
extinction when $\sigma^2/2 > \mu$ due to the predominance of the
stochastic fluctuations.  Despite its simplicity, the model \eqref{eq:intro} is used extensively for projecting future population sizes and estimating extinction risk~\citep{dennis-etal-91,foley-94,lande-etal-03}. For example, \citet{dennis-etal-91} estimated $\mu$ and $\sigma$ for six endangered species. These estimates provided a favorable outlook for  the continued recovery of the Whooping Crane (i.e. $\mu \gg\sigma^2/2$), but unfavorable prospects for the Yellowstone Grizzly Bear. 

Individuals cannot avoid being subject to temporal heterogeneity, but it is only when they disperse that they are affected by spatial variation in the environment. The effect of spatial heterogeneity on population growth depends, intuitively, on how individuals respond to environmental cues~\citep{hastings-83,cantrell-cosner-91,dockery-etal-98,chesson-00,cantrell-etal-06,siap-06,amnat-09b}. When movement is towards regions with superior habitat quality, the presence of spatial heterogeneity increases the rate of population growth~\citep{chesson-00,amnat-09b}. The most extreme form of this phenomenon occurs when individuals are able to disperse freely and ideally; that is, they can move instantly to the locations that maximize their per-capita growth rate~\citep{fretwell-lucas-70,cantrell-etal-07}.  Anthropogenically altered habitats, however, can cause a disassociation between cues  used by organisms to assess habitat quality and the actual habitat quality. This disassociation can result in negative associations between movement patterns and habitat quality and a corresponding reduction in the rate of  population growth~\citep{remes-00,delibes-etal-01,amnat-09b}. For ``random diffusive movement'' (that is, no association between movement patterns and habitat quality), spatial heterogeneity increases population growth rates due to the influence of patches of higher quality.  
However, this boost in growth rate is most potent for sedentary populations~\citep{hastings-83,dockery-etal-98,siap-06,amnat-09a}. This dilutionary effect of dispersal on population growth was observed in the invasion of a woody weed, \emph{Mimosa pigra}, into the wetlands of tropical Australia~\citep{lonsdale-93}. A relatively fast disperser, this weed had a population doubling time of 1.2 years on favorable patches, but 
it exhibited much slower growth at the regional scale (doubling time of 6.7 years)  due to the separation of  suitable wetland habitats by unsuitable eucalyptus savannas.

Despite these substantial analytic advances in understanding separately the effects of spatial and temporal heterogeneity on population growth, there are few analytic studies that consider the combined effects. For well-mixed populations with non-overlapping generations living in patchy environments, \citet{metz-etal-83} showed that population growth is determined by the geometric mean in time of the spatially (arithmetically) averaged per-capita growth rates. A surprising consequence of this expression is that populations coupled by dispersal can persist even though they are extinction prone in every patch~\citep{jansen-yoshimura-98}. This ``rescue effect'', however, only occurs when spatial correlations are sufficiently weak~\citep{harrison-quinn-89}. \citet{prsb-10} extended these results by deriving an analytic approximation for stochastic growth rates for partially mixing populations. This approximation reveals that positive temporal correlations can inflate population growth rates at intermediate dispersal rates, a conclusion consistent with simulation and empirical studies~\citep{roy-etal-05,matthews-gonzalez-07}. For example, \citet{matthews-gonzalez-07} manipulated metapopulations of \emph{Paramecium aurelia} by varying spatial-temporal patterns of temperature.  In spatially uncorrelated environments, the populations coupled by dispersal always persisted for the duration of the experiment, while some of the uncoupled populations went extinct. Moreover, metapopulations experiencing positive temporal correlations exhibited higher growth rates than metapopulations living in temporally uncorrelated environments.

Here, we introduce and analyze stochastic models of populations that continuously experience uncertainty in time and space. For these models, our analysis answers some fundamental questions in population biology such as: 
\begin{itemize}
\item How is the long-term spatial distribution of a population related to its rate of growth? 
\item
When are population growth rates maximized at low, high, or intermediate dispersal rates for populations exhibiting diffusive movement? 
\item
What is ideal free movement for individuals constantly facing uncertainty about local environmental conditions? 
\item To what extent do spatial correlations in temporal fluctuations hamper population persistence?  
\item How do multiple spatial scales of environmental heterogeneity influence population persistence? 
\end{itemize}

In Section~\ref{S:model} we introduce our model for population growth in a patchy environment.  It describes temporal fluctuations in the qualities of the various patches using  multivariate Brownian motions with
correlated components. 

In Section~\ref{S:distribution_growth}, we first consider the vector-valued stochastic process given
by the proportions of the population in each  patch.  These
proportions converge in distribution to a (random) equilibrium at large times.
The probability that this equilibrium spatial distribution
is in some given subset of the set of possible patch
proportions is just the long-term average amount of time that the
process spends in that subset.  We derive
a simple expression for the stochastic growth of the population in terms of
the first and second moments of this equilibrium spatial distribution.
We also show that this equilibrium spatial distribution is characterized by a solution of a PDE that we solve in the case of two patches and use to examine how the equilibrium spatial distribution depends on the dispersal mechanism.  We then present some numerical simulations 
to give a first indication of the interesting range of phenomena
that can occur when there is spatial heterogeneity in per-capita growth rates and biased movement between patches. 

We use the results from Section~\ref{S:distribution_growth} 
in Section~\ref{S:ideal_free} to
investigate ideal free dispersal in stochastic environments.
That is, we determine which forms of dispersal  maximize the
stochastic growth rate for given mean per-capita
growth rates in each of the patches
and given infinitesimal covariances for their temporal fluctuations.

We consider the effect of constraints on  dispersal
 in Section~\ref{S:constraints}.
  We suppose that the dispersal rates are fixed up to
a scalar multiple $\delta$ and establish an analytic approximation 
for the stochastic growth rate of the form $a + b/\delta$ for large $\delta$. We use this approximation to give criteria for whether low, intermediate, or high dispersal rates maximize the stochastic growth rate. In particular, 
we combine this analysis with tools from group representation theory 
to obtain results on the stochastic growth rate
for environments with multiple spatial scales.

We discuss how our results relate to existing literature
in Section~\ref{S:discussion}.  We end with a collection
of Appendices where, for the sake of streamlining
the presentation of our results in the remainder of the paper, 
we collect most of the proofs.

\section{The Model}
\label{S:model}

We consider a population with overlapping generations living in a spatially heterogeneous environment consisting of
 $n$ distinct patches
and suppose that the per-capita growth rates within each patch are determined by a mixture of deterministic and stochastic environmental inputs. Let $X_t^i$ denote the abundance of the population in the $i$-th patch at time $t$
and write $\bX_t = (X_t^1, \ldots, X_t^n)^T$ for the resulting column vector
(we will use the superscript $T$ throughout to denote the transpose
of a vector or a matrix).
If there was no dispersal between patches, 
it is appropriate to model $\bX$ as a Markov process
with the following specifications for $\Delta t$ small:
\[
\E[X_{t+\Delta t}^i - X^i_t \, | \, \bX_t=x] 
\approx 
\mu_i x^i \Delta t,
\]
where $\mu_i$ is the mean per-capita  growth rate in patch $i$,
and 
\[
\mathrm{Cov}[X_{t+\Delta t}^i - X^i_t, \, X_{t+\Delta t}^j - X^j_t \, | \, \bX_t=x] \approx   
\sigma_{ij} x^i x^j \Delta t, 
\] 
where $\Sigma=(\sigma_{ij})$ is a covariance matrix that
captures the spatial dependence between the temporal fluctuations
in patch quality. More formally, we consider 
the system of stochastic differential equations of the form 
 \[
dX^i_t = X_t^i\left( \mu_i  dt + dE_t^i\right),
\]
where $\bE_t=\Gamma^T \bB_t$, $\Gamma$ is an $n \times n$ matrix such
that $\Gamma^T \Gamma= \Sigma$, and 
$\bB_t = (B_t^1, \ldots, B_t^n)^T$, $t \ge 0$,
is a vector of independent standard Brownian motions. 

In order to incorporate dispersal that couples 
the dynamics between patches, let $D_{ij} \ge 0$ for $j\neq i$ be the per-capita rate at which the population in patch $i$ disperses to patch $j$.  Define $-D_{ii} := \sum_{j\neq i} D_{ij}$ to be  the
total per-capita immigration rate out of patch $i$.  
The resulting matrix $D$ has zero row sums and non-negative off-diagonal entries.  We call such matrices \emph{dispersal matrices}. It is worth noting that any dispersal matrix $D$ can be viewed as a generator of a continuous time Markov chain; that is, if we write 
$P_t := \exp(t D)$ for $t \ge 0$, so that $P_t$, $t \ge 0$, solves
the matrix-valued ODE 
\[
\frac{d}{dt} P_t = P_t D,
\]
then the matrix $P_t$ has nonnegative entries, its rows sum to one,  
and 
the Chapman-Kolmogorov relations $P_s P_t = P_{s+t}$ hold for all $s,t \ge 0$.
The $(i,j)$-th entry of $P_t$ gives the proportion of the population
that was originally in patch $i$ at time $0$ but has dispersed to patch $j$
at time $t$.

Adding dispersal to the regional dynamics leads to the 
system of stochastic differential equations 
\begin{equation}\label{eq:main1}
dX^i_t = X_t^i(\mu_i  dt + dE_t^i)+\sum_{j=1}^n D_{ji} X_t^j dt.
\end{equation}
We can write this system more compactly as the vector-valued
stochastic differential equation
\begin{equation}\label{eq:main2}
\begin{split}
d\bX_t & = \diag(\bX_t)\left(\bmu dt +d \bE_t\right) +D^T \bX_t\,dt \\
& = \diag(\bX_t)\left(\bmu dt +\Gamma^T d \bB_t\right) +D^T \bX_t\,dt, \\
\end{split}
\end{equation}
where $\bmu:=(\mu_1,\dots,\mu_n)^T$, 
and, given a vector $u$, we write $\diag(u)$ for the diagonal
matrix that has the entries of $u$ along the diagonal.

We implicitly assume in the above set-up that all dispersing individuals arrive in some patch on the landscape. To account for dispersal induced mortality, we can add fictitious patches in which dispersing individuals enter and experience a mortality rate before dispersing to their final destination. 

Also, our model does not include density-dependent effects on population growth.   However, one can view it as a linearization of a density-dependent model  about the extinction equilibrium  $(0,\dots,0)^T$ and, therefore, \eqref{eq:main2} determines how the population grows when abundances are low.  Moreover, for discrete-time analogues of our model, positive population growth for this linearization  implies persistence in the sense that there exists a unique positive stationary distribution for corresponding
models with compensating density-dependence~\citep{tpb-09}.  We conjecture that the same conclusion  holds for our continuous time model.

 \textbf{From now on we  assume} that the dispersal matrix
$D$ is {\em irreducible} (that is, that it can not be put into
block upper-triangular form by a re-labeling of the patches).
This is equivalent to assuming that the entries of the matrix
$P_t = \exp(tD)$ are strictly positive for all $t>0$, and so it is possible to disperse between any two patches.  Also, we will assume that
the covariance matrix $\Sigma$ has full rank 
(that is, that it is non-singular).
This assumption implies that the randomness in the temporal fluctuations
is genuinely $n$-dimensional.

\section{The stable patch distribution and stochastic growth rate}
\label{S:distribution_growth}

\subsection{Stable patch distribution} 
The key to understanding the asymptotic stochastic
growth rate of the population
is to first examine the dynamics of the spatial distribution of the population.  
Let $S_t := X_t^1 + \cdots + X_t^n$ denote the total population abundance at time $t$
and write  $Y_t^i := X_t^i/S_t$ for the proportion of the total population 
that is in patch $i$.  Set $\bY_t:=(Y^1_t,\dots,Y^n_t)^T$.  
The stochastic process $\bY$ takes values in the 
probability simplex $\Delta:=\{y\in \R^n: \sum_i y_i =1, \, y_i \ge 0\}$.

The following proposition, 
proved in Appendix~\ref{apx:frequencies}, shows that the 
stochastic process $\bY$ is autonomously Markov; that is, that
its evolution dynamics are governed by a stochastic differential equation that does not involve the total population size.  Moreover, 
it says that the law of the random vector $\bY_t$ converges to
a unique  equilibrium as $t \to \infty$. Recall, the \emph{law of a random vector $\bY\in \R^n$} is the probability measure $\mu_\bY$ on $\R^n$ defined by $\mu_\bY(A)=\P\{ \bY\in A\}$ for all Borel sets $A\subseteq \R^n$. Moreover, for any $\mu_\bY$-integrable function $h:\R^n\to\R$, the \emph{expectation of $h(\bY)$} is defined by  
\[
\E[h(\bY)]=\int h(y) \,\mu_\bY(dy).
\]
A sequence of random vectors $\bY_1,\bY_2,\dots$ \emph{converges in law} to a random vector $\bY_\infty$ if 
\[
\lim_{n\to\infty} \E[h(\bY_n)] = \E [ h(\bY_\infty)] 
\]
for every continuous, bounded function $h:\R^n \to \R$. Convergence in law of a sequence of random vectors is also called convergence in distribution of the random vectors  and is equivalent to weak convergence of their laws.

\begin{proposition}\label{prop:frequencies} Suppose that $\bX_0 \ne 0$. Then, the stochastic process $\bY$ satisfies the stochastic differential equation
\begin{equation}\label{eq:freq}
d\bY_t =  \left( \mathrm{diag}(\bY_t) - \bY_t \bY_t^T  \right) \Gamma^T  d\bB_t    + D^T \bY_t  dt + \left( \mathrm{diag}(\bY_t) - \bY_t \bY_t^T  \right) \left( \mu - \Sigma \bY_t\right)  dt.
\end{equation}
Moreover, there exists a random variable 
$\bY_\infty$  taking values in the probability simplex $\Delta$ 
such that $\bY_t$ converges in law to $\bY_\infty$ as $t \to \infty$
and such that the empirical measure $\Pi_t := \frac{1}{t}\int_0^t \delta_{\bY_s} \, ds$ converges almost surely to the law of $\bY_\infty$
as $t \to \infty$.  The law of $\bY_\infty$ does not depend on $X_0$.
\end{proposition}

The empirical probability measure $\Pi_t$ appearing in Proposition~\ref{prop:frequencies}
describes the proportions of the time interval $[0,t]$ that
the process $\bY$ spends in the various subsets of its state space $\Delta$.
Namely, for a Borel set $A\subseteq \Delta$ of patch occupancy states, $\Pi_t(A)$ equals the fraction of time spent in these states over the time interval $[0,t]$. For example, if $A=\{y\in \Delta: y_1>1/2\}$, then $\Pi_t(A)$ equals the fraction of time for which at least 50\% of the population is in patch $1$ during the time interval $[0,t]$.

\subsection{Stochastic growth rates.} 
Recall that $S_t= X_t^1 + \cdots + X_t^n$  
is the total population size at time $t$.  That is,
$S_t = \bone^T \bX_t$, where $\bone = (1, \ldots, 1)^T$.
Because $D \bone = 0$, it follows from \eqref{eq:main2}
that  
\[ dS_t
 = \bX_t^T \Gamma^T d\bB_t + \mu^T \bX_t dt\\
 = S_t\bY_t^T \Gamma^T d\bB_t + S_t \mu^T \bY_t dt.
\]
Therefore, by It\^o's lemma \citep{gardiner-04},
\[
 \log S_t
 = S_0+\int_0^t \bY_t^T \Gamma^T d\bB_t + \int_0^t\mu^T \bY_t dt  - \frac{1}{2} \int_0^t \bY_t^T \Gamma^T \Gamma \bY_t dt.
\]
Dividing by $t$, taking the limit as $t\to\infty$, and applying Proposition~\ref{prop:frequencies} yields the following result.

\begin{theorem}\label{thm:nice}
Suppose that $\bX_0 \ne 0$. 
Then, 
\begin{equation}\label{eq:lyapunov}
\chi:= \lim_{t \to \infty} t^{-1} \log S_t 
=
 \mu^T \mathbb{E}[ \bY_\infty]  - 
\frac{1}{2}\E \left[\bY_\infty^T \Sigma \bY_\infty \right]\quad \text{almost surely},
\end{equation}
where $\bY_\infty$ is described in Proposition~\ref{prop:frequencies}.
\end{theorem}
The limit $\chi$ in \eqref{eq:lyapunov} is generally known as the 
{\em Lyapunov exponent} for the Markov process $\bX$. 
Following \citet{tuljapurkar-90}, we also call $\chi$ the \emph{stochastic growth rate} of the population, 
as it describes the asymptotic growth rate of the population in the presence of stochasticity.  
To interpret \eqref{eq:lyapunov}, notice that 
\begin{equation}
\label{eq:weighted_rate}
\langle \mu \rangle:= \mu^T \E[\bY_\infty] 
= \sum_i \mu_i \E[Y_\infty^i] 
= \lim_{t \to \infty} \sum_i \mu_i \E[Y_t^i]
\end{equation} 
corresponds to weighted average of the per-capita
growth rates  with respect to the long-term spatial distribution $\bY_\infty$ of the population. 
To interpret the other component of \eqref{eq:lyapunov}, let $\Var[X]$ denote the variance of a random variable $X$. 
Since $\sum_i Y_t^i (E_{t+\Delta t}^i-E_t^i)$ for small $\Delta t>0$ is approximately the average environmental change experienced by the population over time interval $[t,t+\Delta]$,  
\begin{equation}
\label{eq:weighted_variance}
\langle \sigma^2 \rangle=\E \left[\bY_\infty^T \Sigma \bY_\infty \right]
=
\lim_{t \to \infty} \frac{1}{\Delta t}\Var\left[\bY_t^T (\bE_{t + \Delta t} - \bE_t)\right]
=\lim_{t \to \infty} \frac{1}{\Delta t}\Var\left[\sum_i\bY_t^i (E^i_{t + \Delta t} - E^i_t)\right]
\mbox{ for any }\Delta t>0
\end{equation}
corresponds to the infinitesimal variance of the  environmental fluctuations weighted by the  long-term spatial distribution. 

\begin{biology}{ of Theorem~\ref{thm:nice}} The stochastic growth rate $\langle \mu \rangle - \langle \sigma^2 \rangle/2 $ for a spatially structured population is just what we  see for an unstructured population
where  $\langle \mu \rangle$ and $\langle \sigma^2\rangle$ are the per-capita growth rate and the infinitesimal
covariances of the temporal fluctuations averaged appropriately with
respect to the equilibrium spatial distribution. Hence, as in a spatially homogeneous environment, environmental fluctuations reduce the population growth rate. However, as we show in greater detail below, interactions between dispersal patterns, spatial heterogeneity and environmental fluctuations may increase the stochastic growth rate by increasing $\langle \mu\rangle$ or decreasing $\langle \sigma^2\rangle$.  \end{biology}

To get a more explicit expression for the stochastic growth rate, we need to determine the distribution of the equilibrium  $\bY_\infty$, or at least
find its first and second moments. 
This problem reduces to solving for the time-invariant solution of the Fokker-Planck equations with appropriate boundary conditions~\citep{gardiner-04},
Namely, the density $\rho:\Delta \to [0,\infty)$ of $\bY_\infty$ satisfies
\begin{eqnarray}\label{PDE}
  -\sum_i \frac{\partial}{\partial y_i} M_i(y)\rho(y)+\frac{1}{2}\sum_{i,j} \frac{\partial^2}{\partial y_i \partial y_j} V_{ij}(y) \rho(y) &=& 0 \quad \mbox{for} \; y \in \Delta ,
% 0&=& \rho(y) \mbox{ for }y\in \left\{y \in \Delta: \prod_i y_i =0\right\}\\
 % \nonumber 1&=& \int_\Delta \rho(y) \,dy ,
\end{eqnarray}
where $M_i$ and $V_{ij}$ are the entries of 
\[
M(y)=  D^T y + \left( \mathrm{diag}(y) - y y^T  \right) \left( \mu - \Sigma y \right)  
\mbox{ and }
V(y)=\left( \mathrm{diag}(y) - y y^T  \right) \Gamma^T \Gamma  \left( \mathrm{diag}(y) - y y^T  \right) ,
\]
respectively,
and $\rho$ is constrained to have $\int_\Delta \rho(y) dy = 1$.
However, the PDE \eqref{PDE} needs
to be supplemented with appropriate boundary conditions.
In principle, these are found by characterizing the domain of the
infinitesimal generator of the Feller diffusion process $Y$ and thence
characterizing the domain of the adjoint of this operator \citep{khasminskii1960ergodic,bhattacharya1978criteria,bogachev2002uniqueness,bogachev2009elliptic}.
This appears to be a quite difficult problem.
However, in the case of two patches, the problem simplifies to solving an ODE on the unit interval.

\begin{figure}[t]
\includegraphics[width=7in]{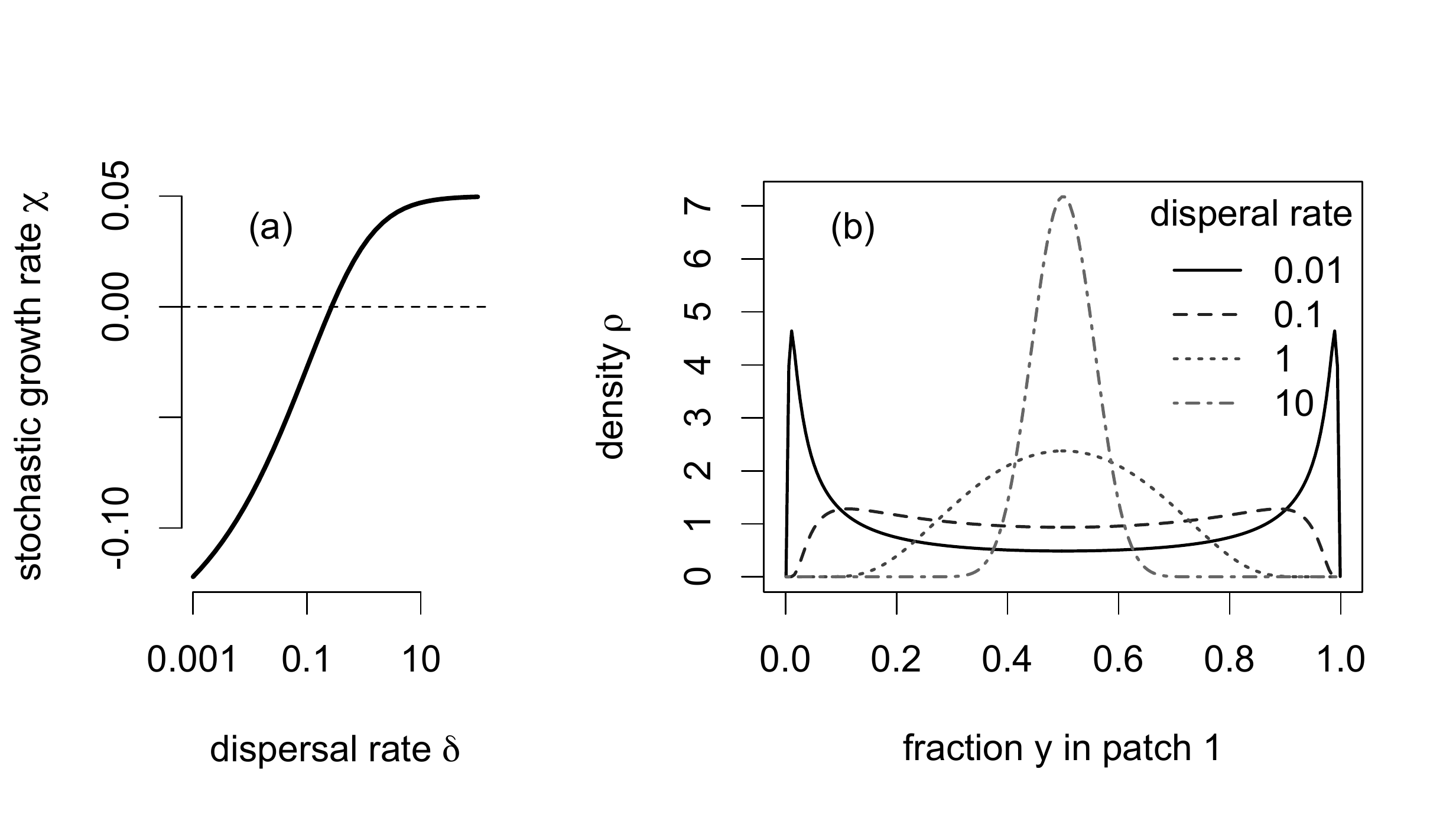}
\caption{Spatial distribution and population growth in a two patch environment. In (a), the stochastic growth rate $\chi$ is plotted as a function of the dispersal rate $\delta$.  In (b), the stationary density of the fraction of individuals in patch $1$ is plotted for different dispersal rates. Parameter values are $\mu_1=\mu_2=0.3$, $\sigma_1=\sigma_2=1$, and $D_{12}=D_{21}=\delta$.}\label{fig:delta}
\end{figure}

\begin{example}[Stochastic growth in two patch environments]\label{ex:2patch}
Assume there are two patches. For simplicity, suppose
 there are no environmental correlations between the patches; that is, that  $\sigma_{ii}=\sigma_i^2$ and $\sigma_{ij}=0$ for $i\neq j$. Proposition~\ref{prop:frequencies} gives that  $Y_t^1 =X_t^1/(X_t^1+X_t^2)$ satisfies the one-dimensional stochastic differential equation
\[
dY_t^1 = M_*(Y_t^1)\,dt+\sqrt{V_*(Y_t^1)}\,dB_t
\]
where 
\[
\begin{split}
M_*(y)&:=y(1-y)(\mu_1-\mu_2-\sigma_1^2 y +\sigma_2^2 (1-y))-D_{12}y+D_{21}(1-y)\\
& \text{and}\\
V_*(y)&:= y^2(1-y)^2 (\sigma_1^2+\sigma_2^2).\\
\end{split}
\]
We can then apply standard tools for one-dimensional diffusions \citep{gardiner-04}
(checking that the boundaries at 0 and 1 are ``entrance'', and hence inaccessible)
to find that the density $\rho(x):[0,1]\to [0,\infty)$ of  $Y^1_\infty$ is given by
%\[
%    \frac12 (V_*\rho)'' - (M_*\rho)' = 0 \mbox{ on }(0,1)\mbox{ with } \rho(0)=\rho(1)=0.
%\]
% If we integrate once and rearrange, we see that
% \[
%     (V_*\rho)' = 2 \frac{M_*}{V_*} (V_*\rho),
% \]
% an equation that is solved by
\begin{eqnarray*}
\rho(y)&=& \frac{C_1}{V_*(y)} \exp\left( 2\int \frac{M_*(y)}{V_*(y)}\,dy\right)\\
    &=&\frac{C_2}{y^2(1-y)^2} \exp\left( \frac{2}{\sigma_1^2+\sigma_2^2}
    \int \frac{\mu_1-\mu_2}{y(1-y)}-\frac{\sigma_1^2}{1-y}+\frac{\sigma_2^2}{y}-\frac{D_{12}}{y(1-y)^2}+\frac{D_{21}}{y^2(1-y)}\,dy \right)\\
    &=& C_3 \,y^{\beta-\alpha_1}(1-y)^{-\beta-\alpha_2}\exp\left(-\frac{2}{\sigma_1^2+\sigma_2^2}\left( \frac{D_{21}}{y}+\frac{D_{12}}{1-y}\right)\right),
\end{eqnarray*}
where the $C_i$  are normalization constants, and 
\begin{eqnarray*}
\alpha_i&:=& \frac{2 \sigma_i^2}{\sigma_1^2+\sigma_2^2}\\
\beta&:=& \frac{2}{\sigma_1^2+\sigma_2^2}\left( \mu_1-\mu_2 +D_{21}-D_{12}\right).
\end{eqnarray*}
Using this expression in \eqref{eq:lyapunov}, we get the following explicit 
expression for the stochastic growth rate 
\begin{eqnarray*}
\chi&=&\mu_1\int_0^1y\rho(y)\,dy+\mu_2\int_0^1 (1-y)\rho(y)\,dy-\frac{\sigma_1^2}{2}\int_0^1 y^2\rho(y)\,dy -\frac{\sigma_2^2}{2}\int_0^1(1- y)^2\rho(y)\,dy \\
&=& \mu_2-\frac{\sigma_2^2}{2}+(\mu_1-\mu_2+\sigma_2^2)\int_0^1 y \rho(y)\,dy-\frac{\sigma_1^2+\sigma_2^2}{2}\int_0^1 y^2 \rho(y)\,dy .
\end{eqnarray*}

Despite its apparent complexity,  this formula provides insights into how dispersal may influence population growth. For example, consider a population dispersing diffusively between statistically similar but uncorrelated patches (that is,  $D_{12}=D_{21}=\delta/2$, $\mu_1=\mu_2=\mu$, and $\sigma_1=\sigma_2=\sigma$). We claim that the stochastic growth rate  $\chi$ is an increasing function of the dispersal rate $\delta$.  Intuitively, this occurs because increasing $\delta$ decreases the variance of the random variable $\bY_\infty$ but has no effect on its expectation. 

To verify our claim that $\chi$ is increasing with $\delta$, 
write $\rho(\cdot;\delta)$ for the density of $Y_\infty^1$
to emphasize its dependence on $\delta$ and
notice that in this case 
\[
\rho(y;\delta) =\frac{1}{C(\delta)} y^{-1}(1-y)^ {-1} \exp \left (  - \frac{\delta}{2\sigma^2y(1-y)} \right ), \quad y \in (0, 1),
\] 
where $C(\delta)=\int_0^1 y^{-1}(1-y)^ {-1} \exp \left (  - \frac{\delta}{2\sigma^2y(1-y)} \right )\,dy$ is the normalization constant and 
\begin{equation}\label{twopatch}
\chi(\delta) = \mu - \sigma^2/2 + \sigma^2 \int_0^1 y(1-y)\,\rho(y;\delta)\,dy.
\end{equation}
It suffices to show that
\begin{eqnarray*}
 \int_0^1 y(1-y)\rho(y; 2\delta\sigma^2)\,dy &=&
 \frac{\int_0^1  \exp \left (  - \frac{\delta}{y(1-y)} \right ) \, dy}{C(2\delta \sigma^2)}\\
&=& \frac{\int_0^1  \exp \left (  - \frac{\delta}{y(1-y)} \right ) \, dy}{ \int_0^1  y^{-1}(1-y)^{-1} \exp \left (  - \frac{\delta}{y(1-y)} \right ) \, dy}\end{eqnarray*}
is an increasing function of $\delta>0$. Differentiating with respect to  $\delta$ and carrying the differentiation inside the integral sign, we obtain
\[ \begin{split} C(2 \sigma^2 \delta)^{-2} \times \left[ \int_0^1  y^{ -2}(1-y)^{-2} \exp \left(  - \frac{\delta}{y(1-y)} \right) \, dy
\times  \int_0^1  \exp \left( - \frac{\delta}{y(1-y)} \right) \, dy \right.  \\
-  \left. \left(\int_0^1  y^{-1}(1-y)^{-1} \exp \left(  - \frac{\delta}{y(1-y)} \right) \, dy \right )^2\right].
\end{split}\]
This quantity is the variance 
of the random variable
$\left( Y_\infty^1(1-Y_\infty^1) \right)^{-1}$
and is thus nonnegative.

For the purpose of comparison with
general asymptotic approximations that we develop later, we
note that after a change of variable
\[
\frac{
\int_0^1 
\exp \left (  - \frac{\delta}{2\sigma^2 y(1-y)} \right )
\, dy
}
{
\int_0^1 
y^{-1}(1-y)^{-1}
\exp \left (  - \frac{\delta}{2\sigma^2 y(1-y)} \right )
\, dy
} 
=
\frac{
\int_0^\infty e^{-z} z^{-\frac{1}{2}} (\frac{2 \sigma^2 z}{\delta} +4)^{-\frac{3}{2}} \, dz
}
{
\int_0^\infty e^{-z} z^{-\frac{1}{2}} (\frac{2 \sigma^2 z}{\delta} +4)^{-\frac{1}{2}} \, dz
}. 
\]
Upon expanding the two functions 
$w \mapsto (w+4)^{-\frac{1}{2}}$
and 
$w \mapsto (w+4)^{-\frac{3}{2}}$
in Taylor series around $0$ and integrating, we find that
the ratio of integrals is of the form
\[
\frac{1}{4} - \frac{1}{\delta} \frac{\sigma^2}{16} 
+ \mathrm{O}\left(\frac{1}{\delta^2}\right)
\]
as $\delta \to \infty$, so that
\begin{equation}
\label{eq:twopatch_asymptotics}
\chi(\delta) \approx \mu - \frac{\sigma^2}{4} - \frac{1}{\delta} \frac{\sigma^4}{16}
\end{equation}
as $\delta \to \infty$. 

Approximation~\eqref{eq:twopatch_asymptotics} implies, as we prove more generally in Proposition~\ref{prop:limit}, that $\lim_{\delta\to\infty} \chi(\delta)=\mu-\sigma^2/4$. % Since $\chi(\delta)$ in this case is an increasing function, we can derive the following

\begin{biology}{ of Example~\ref{ex:2patch}}Even if two patches are unable to sustain a population in the absence of dispersal, connecting the patches by dispersal can permit persistence. This phenomenon occurs only at intermediate levels of environmental stochasticity (i.e.  $2\mu<\sigma^2<4\mu$). Moreover, when this phenomenon occurs, there is a critical dispersal threshold $\delta^*>0$ such that the metapopulation decreases to extinction whenever its dispersal rate is too low (i.e. $\delta\le \delta^*$) and persists otherwise (Fig.~\ref{fig:delta}). 
\end{biology}

\end{example}

Because there do not appear to be closed-form expressions for
the law of the stable patch distribution $\bY_\infty$ when
there are more than two patches, we must seek other routes to 
understanding the stochastic growth rate in such cases.  One approach
would be to solve the PDE \eqref{PDE} numerically.
A second approach would be to simulate the stochastic process $\bY$
for long time intervals and derive approximate values for the
first and second moments of the equilibrium distribution.
To give an indication of the range of phenomena that
can occur in even relatively simple systems where there
is biased movement between patches, we adopt the
even simpler approach of simulating the stochastic process $\bX$
directly for long time intervals to obtain an approximate
value of the stochastic growth rate.
We implemented the simulations in a manner similar to that of \citet{talay1991lyapunov},
and the \texttt{R} code used is provided as supplementary material.

\begin{example}[Spatially heterogeneous environments with biased emigration]  \label{ex:numerics}
For these simulations, we consider a metapopulation with either $n=8$ or $n=40$ patches of which one quarter are higher quality ($\mu_i=10$ in these patches) and the remainder are lower quality ($\mu_i=1$ in the remaining patches). All patches have the same level of spatially uncorrelated environmental noise ( $\sigma_{ii}=16$ for all $i$ and $\sigma_{ij} = 0$ for $i\neq j$). When an organism exits a patch it chooses from the other patches with equal probability, but the
emigration rate from a patch depends on the patch quality. 

\begin{figure}[t!!]
\includegraphics[width=7in]{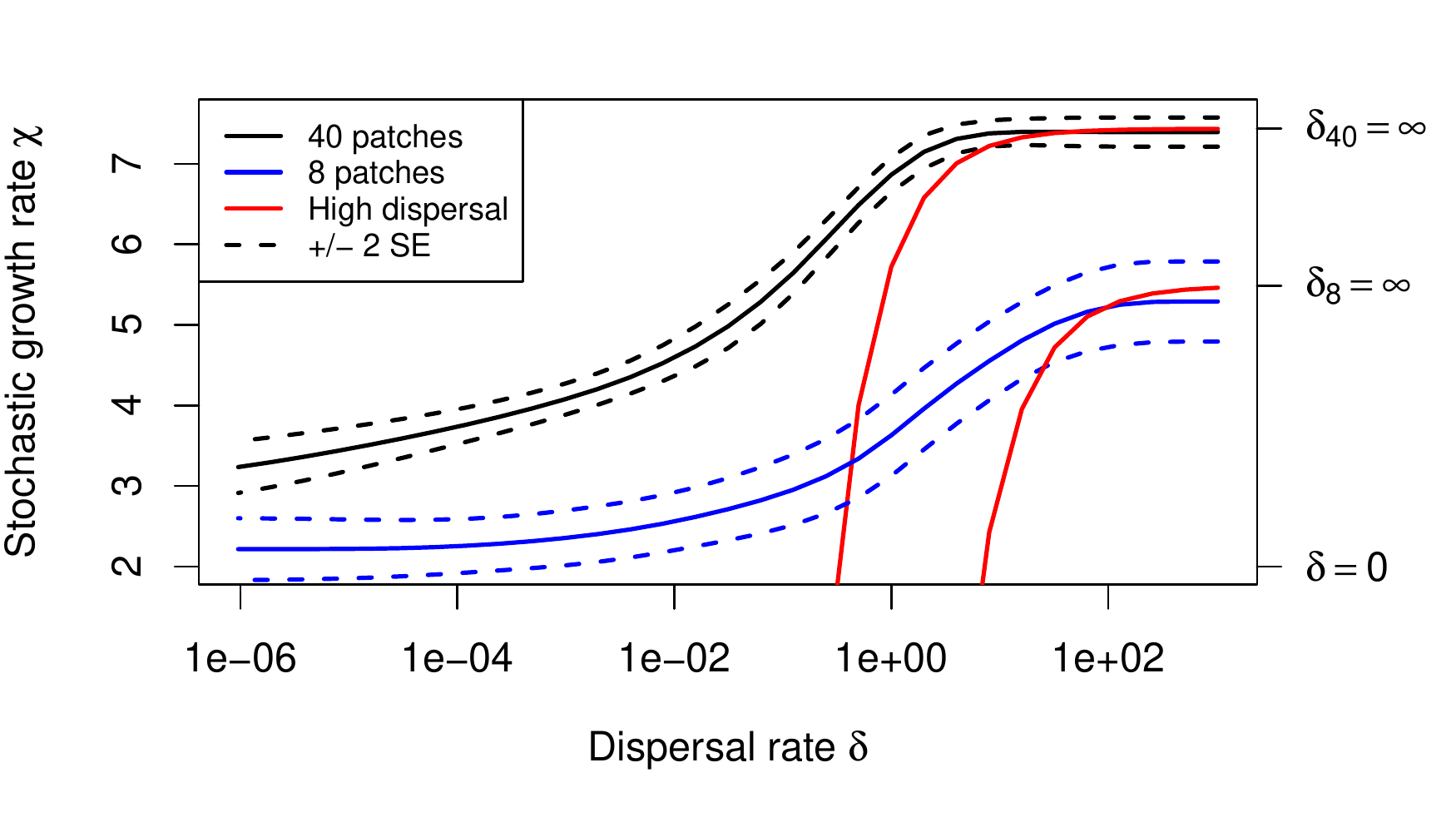}
\caption{The effect of dispersal rate $\delta$ on populations emigrating
more rapidly out of lower quality patches than higher quality patches.
Shown is the stochastic growth rate $\chi$ estimated from simulation of the
SDE for 100 time units, across a range of values of $\delta$, for both a
40-patch and a 8-patch model.  Standard errors are estimated using the
standard deviation of the stochastic growth rates across nonoverlapping time
segments of a given simulation.  Details of the dispersal matrix and
parameter values are described in the main text.
The right-hand axis shows asymptotic values for $\delta = 0$ and
$\delta = \infty$, which are: $\chi(0) = \max_i \mu_i$ and
$\chi(\infty)=\mu^T \pi - \frac{1}{2} \pi^T \Sigma \pi$ (Proposition
\ref{prop:limit}).
``High dispersal'' shows the approximation of the form $\chi(\delta) \approx a + b/\delta$ for large $\delta$
calculated from formula~\eqref{eq:bigone} of Theorem 2. }
\label{fig:numericsadaptive}
\end{figure}

First, we consider the case in which emigration is ``adaptive''
in the sense that individuals emigrate more rapidly out of lower quality patches than higher quality patches:
\[
D_{ij}=\left\{
\begin{array}{cc}
\delta, \quad & \mbox{for }i=1,\dots,n/4\mbox{ and }i\neq j , \\
10\,\delta, \quad & \mbox{for }i=n/4+1,\dots, n \mbox{ and } i\neq j.
\end{array}
\right.
\]
Here, the  parameter $\delta>0$ scales the emigration rate, so that doubling $\delta$ doubles the emigration rate from all patches.  As expected, since in this case dispersal is ``adaptive'', Figure~\ref{fig:numericsadaptive} shows that stochastic growth rate $\chi=\chi(\delta)$ as a function of $\delta$ increases with $\delta$. Moreover, Figure~\ref{fig:numericsadaptive} shows asymptotic values at $\delta=\infty$ for each case, and illustrates that the analytic approximation developed later in Theorem~\ref{thm:bigone} works reasonably well for large values of $\delta$.
The Figure also shows extremely slow convergence as $\delta\to 0$ to $\chi(0)=\max_i \mu_i-(1/2)\sigma_i^2$
(note the logarithmic scale on the horizontal axis),
indicating that although $\chi$ is continuous at $\delta=0$ by Proposition~\ref{prop:cont} below,
it may not be differentiable there.

\begin{figure}[t!!]
\includegraphics[width=7in]{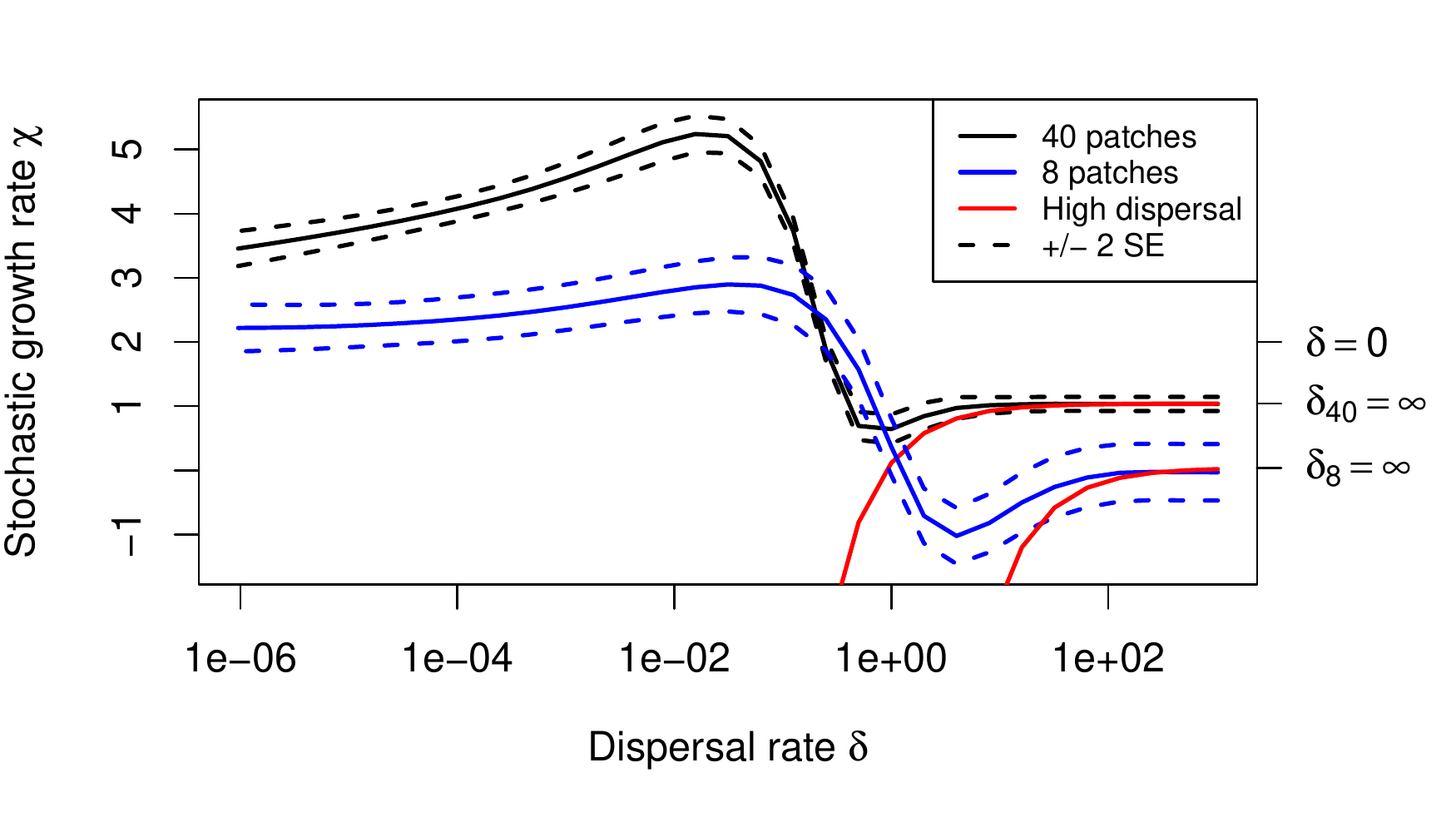}
\caption{The effect of dispersal rate $\delta$ on populations emigrating
more rapidly out of higher quality patches than lower quality patches.
Details are as in Figure \ref{fig:numericsadaptive}, but with different
dispersal scheme; parameter values are described in the main text.}
\label{fig:numericsmal}
\end{figure}

Next we consider a case in which emigration is ``maladaptive'',
in the sense that individuals emigrate more rapidly out of higher quality patches than out of lower quality patches:
\[
D_{ij}=\left\{
\begin{array}{cc}
10\,\delta, \quad & \mbox{for }i=1,\dots,n/4\mbox{ and }i\neq j,\\
\delta, \quad & \mbox{for }i=n/4+1,\dots, n \mbox{ and } i\neq j.
\end{array}
\right.
\]
It is possible to show using the results of Section~\ref{S:constraints} below that in this regime, high dispersal rates lead to a {\em lower} stochastic growth rate than sedentary populations (that is, 
$\lim_{\delta \to \infty} \chi(\delta)$ is dominated by 
$\lim_{\delta \to 0} \chi(\delta)$), and yet $\chi(\delta)$ increases 
with $\delta$ when $\delta$ is large. 
As illustrated in Figure~\ref{fig:numericsmal}, the stochastic growth rate $\chi(\delta)$ exhibits a rather complex
dependence on $\delta$: increasing at low dispersal rates,  declining at higher dispersal rates, and finally increasing  again at the highest dispersal rates.

In a conservation framework, increasing $\delta$ corresponds to facilitating movement between patches by increasing the size or number of dispersal corridors between patches. 
\begin{biology}{ of Example~\ref{ex:numerics}}
For populations exhibiting adaptive movement, increasing the size or number of dispersal corridors between patches enhances metapopulation growth rates. For populations exhibiting maladaptive movement, however, increasing dispersal rates can either increase or decrease metapopulation growth rates. 
\end{biology} 
 
\end{example}

\section{Ideal free dispersal in a stochastic environment}
\label{S:ideal_free}

A basic quandary in evolutionary ecology is, ``For a given set of environmental conditions, what dispersal pattern maximizes fitness?''  Since fitness in our context corresponds to the stochastic growth rate of the population, we can rephrase this question as, ``Given $\mu$ and $\Sigma$,  what form of the dispersal matrix $D$ maximizes $\chi$?'' Following \citet{fretwell-lucas-70}, we call such an optimal dispersal mechanism {\em ideal free dispersal} as individuals  have no constraints on their dispersal (i.e. are ``free'') and have complete knowledge about the distribution of spatial-temporal fluctuations (i.e. are ``ideal'').

 Equation~\eqref{eq:lyapunov} provides a means to answer this question. Because $\Sigma$ has full rank, the function $y\mapsto \frac{1}{2} y^T\Sigma y$ is strictly convex, and so Jensen's inequality implies that 
\[
\E[\bY_\infty^T \Sigma \bY_\infty]
\ge \E[\bY_\infty]^T \Sigma\, \E[\bY_\infty],
\] 
with equality if and only if the random vector $\bY_\infty$ is almost surely constant.  Hence, to maximize the stochastic growth rate $\chi$, we need to eliminate the variability in $\bY_\infty$, so that $\bY_\infty=y$ almost surely for a constant $y$ 
that is chosen to maximize 
\begin{equation}\label{eq:lyap3}
\mu^T y  - 
\frac{1}{2}y^T \Sigma\, y
\end{equation}
subject to the constraint $y\in \Delta$.  Under our standing
non-degeneracy assumptions on $D$ and $\Sigma$, the law
of $\bY_\infty$ is supported on all of $\Delta$, and so we cannot
actually achieve a situation in which $\bY_\infty$ is a constant.
However, the following result, which we prove in  Appendix~\ref{apx:limit}, shows  that we can approach this regime arbitrarily closely.  Recall that the \emph{stationary distribution $\pi$} for an irreducible dispersal matrix $Q$ is a probability vector $\pi \in \Delta$ such that $\pi^T Q=0$.  We note that any vector $\pi$ in the interior of $\Delta$ is the stationary distribution for some irreducible dispersal matrix $Q$. 
For example, given $\pi$, we can define $Q = \bone \pi^T - I$ where $I$ denotes the identity matrix.

\begin{proposition}~\label{prop:limit} Consider a vector $\pi$
in the interior of $\Delta$ and an irreducible dispersal matrix $Q$ that has $\pi$ as its unique stationary distribution.  
Let $\bY_\infty(\delta)$ be the equilibrium patch distribution and $\chi(\delta)$ be the stochastic growth rate for \eqref{eq:main2} with  $D=\delta Q$. Then $\bY_\infty(\delta)$ converges in law to the constant vector $\pi$ as $\delta\to \infty$, and $\chi(\delta)$ converges to 
$\mu^T \pi  - 
\frac{1}{2}\pi^T \Sigma\, \pi$ as $\delta\to\infty$.
\end{proposition}

In the absence of population growth due to deterministic or stochastic effects, each of the dispersal matrices $\delta Q$ in Proposition~\ref{prop:limit} sends the patch distribution
to the vector $\pi$ regardless of the initial conditions, 
and the speed at which this happens increases
with $\delta$, so that it becomes effectively instantaneous
for large $\delta$.  Proposition~\ref{prop:limit} says that this
push towards a deterministic equilibrium overcomes any disruptive
effects introduced by population growth provided $\delta$ is
sufficiently large, and so it is possible to produce random
equilibrium patch distributions that are arbitrarily close
to any given vector $\pi$ in the interior of $\Delta$.  If we further
approximate vectors $\pi$ on the boundary of $\Delta$ by ones
in the interior, we see that it is possible to produce 
equilibrium patch distributions that are arbitrarily close
to any given vector in $\Delta$.

Given that any patch distribution can be approximated arbitrary closely
by the equilibrium patch distribution of a suitable
population of rapidly dispersing individuals, 
the problem of optimizing $\chi$ reduces, as we have already noted, 
to maximizing the strictly concave
function $g(y)=\mu^T y  -  \frac{1}{2}y^T \Sigma\, y$ 
over the compact, convex set $\Delta$.  This concavity implies there exists at most one local maximum. Denote this unique
maximizer by $y^*=(y_1^*,\dots,y_n^*)^T$.  

It is optimal 
for all individuals to remain in the single patch $k$
(that is, $y_k^* = 1$) only if
\[
\frac{\partial g}{\partial y_i}(e_k)- \frac{\partial g}{\partial y_k}(e_k)
=
\mu_i - \sigma_{ik} -\mu_k + \sigma_{kk} <0 \text{ for all } i\neq k,
\]
where $e_k$ is the $k$-th element of the standard basis of $\R^n$,
or, equivalently, 
\begin{equation}\label{eq:stay-one-place}
\mu_k-\mu_i > \sigma_{kk}-\sigma_{ik}\mbox{ for all }i\neq k.
\end{equation}
\begin{biology}{ of equation~\eqref{eq:stay-one-place}} If the variances of environmental fluctuations are sufficiently large in all patches and the spatial covariances in these environmental fluctuations are sufficiently small, then ideal free dispersers occupy multiple patches. \end{biology}

When it is optimal to disperse between several patches, we can solve for the optimal dispersal strategy $y^*$ by using the method of Lagrange multipliers. Without loss of generality, assume that the optimal strategy $y^*$ makes use of all patches,  that is, that $y^*$ is in the interior of $\Delta$. Indeed, if the optimal strategy does not make use of all patches, then we can consider analogous problems 
on the faces of the convex polytope
$\Delta$ of the form $\{y \in \Delta : y_i = 0, \, i \in A\}$,
where $A$ is a subset of $\{1, \ldots, n\}$. 
Because 
\[
\nabla g(y)=\mu -\Sigma y  \mbox{ and } \nabla \left(\sum_i y_i\right) = \bone, 
\]  
the optimal $y^*$ must satisfy 
\begin{equation}\label{eq:idf}
\mu -\Sigma y^*= \lambda \bone,
\end{equation}
where $\lambda$ is a Lagrange multiplier. Notice that 
\[
(\Sigma y)_i 
= 
\frac{1}{\Delta t} \E\left[(E_{t + \Delta t}^i - E_t^i) \sum_j y_j (E_{t + \Delta t}^j - E_t^j)\right].
\] 
Hence, we get the following interpretation. 
\begin{biology}{ of equation~\eqref{eq:idf}}
Ideal free populations using multiple patches are distributed across the patches in such a way that  the differences between the mean per-capita growth rates and the covariances between the within patch noise and the noise experienced on average by an individual are equal in all occupied patches. In particular, the local stochastic growth rates $\mu_i -\sigma_{ii}/2$ need not be equal in all occupied patches. 
\end{biology} 
 
Now,
\begin{equation}\label{eq:opt1}
y^*=\Sigma^{-1} (\mu -\lambda\bone),
\end{equation} 
and the constraint $\bone^T y=1$ yields
\[
1= \bone^T\Sigma^{-1} (\mu -\lambda \bone), \\
\]
so that
\begin{equation}\label{eq:opt2}
\lambda =\frac{ \bone^T \Sigma^{-1} \mu -1}{\bone^T\Sigma^{-1} \bone}
\end{equation}
and
\begin{equation}\label{eq:opt3}
y^* = 
\Sigma^{-1} 
\left(\mu  - \frac{ \bone^T \Sigma^{-1} \mu -1}{\bone^T\Sigma^{-1} \bone} \bone \right).
\end{equation}
The right-hand side of equation \eqref{eq:opt3} is the optimal vector $y^*$ we seek, provided that it belongs to the interior of $\Delta$. 
Otherwise, as we remarked above, we need to perform similar analyses
on the faces of the simplex $\Delta$.

To illustrate the utility of this formula, we examine two special cases: 
when the environmental noise between patches is uncorrelated, and 
when the patches experience the same individual levels of noise but they are spatially correlated.  

\begin{figure}[t]
\begin{tabular}{c}
\includegraphics[width=6in]{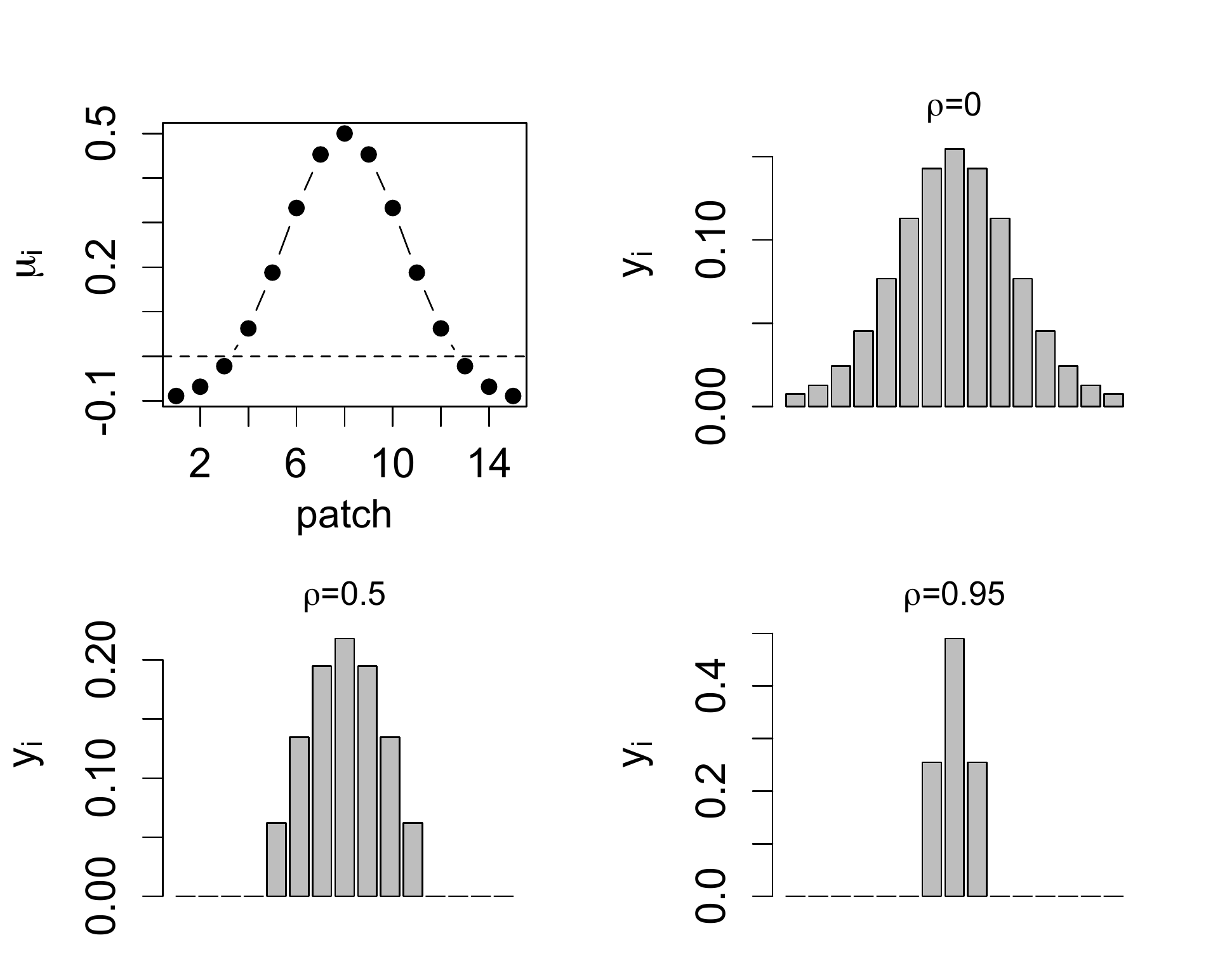}
\end{tabular}
\caption{Effects of spatial correlations on the ideal free patch distribution in a 15 patch environment. Per-capita growth rates $\mu_i$ are plotted in the top left. The ideal free patch distribution $y^*$ is plotted at three levels of spatial correlation $\rho$. Covariances are $\sigma_{ii}=2$ and $\sigma_{ij}=2\rho $ for $i\neq j$. }\label{fig:two}
\end{figure}

\begin{example}[Spatially uncorrelated environments]
Suppose that there are no spatial correlations in 
the environmental noise, so that $\Sigma$ is a diagonal matrix with diagonal entries $\sigma_{ii} = \sigma_i^2$.  It follows from
equation \eqref{eq:opt3}  that the ideal free patch distribution is 
\begin{equation}\label{eq:ifd1}
y_i^* = \frac{1}{\sigma_i^2\sum_j 1/\sigma_j^2} \left[ \sum_j \frac{\mu_i-\mu_j}{\sigma_j^2}+1\right], 
\end{equation}
provided that $\sum_j (\mu_j-\mu_i)/\sigma_j^2<1$  for all $i$. 
\begin{biology}{ of equation~\eqref{eq:ifd1}} In the absence of spatial correlations in environmental fluctuations, ideal free dispersers visit all patches  whenever 
the environmental variation is sufficiently 
great relative to differences in the mean
per-capita growth rates. In particular, if all 
mean per-capita growth rates are equal, then the fraction of individuals in a patch is inversely proportional to the  variation in temporal fluctuations
in the patch; that is, $y_i^* =(1/\sigma_i^2)/(\sum_j 1/\sigma_j^2)$. \end{biology}
\end{example}

\begin{example}[Spatially correlated environments]
Suppose that the infinitesimal variance of the temporal fluctuations in each 
patch is $\sigma^2$ and that the correlation between the fluctuations in
any pair of patches is $\rho$. Thus, $\Sigma = \sigma^2 (1-\rho) I+ \sigma^2 \rho J$, where $J = \bone \bone^T$ is the matrix in which every entry is
$1$. Provided that $-\frac{1}{n-1}< \rho <1$,  the matrix $\Sigma$ 
is non-singular with inverse 
 \[
 \Sigma^{-1}=\frac{1}{(1-\rho)\sigma^2} I - \frac{\rho}{(1-\rho)(1+(n-1)\rho)\sigma^2}J.
 \]
Denoting by $\bar \mu = \frac{1}{n}\sum_i \mu_i$ 
the average across the patches of the mean per-capita growth rates, 
the optimal dispersal strategy is given by 
 \begin{equation}\label{eq:ifd2}
 y_i^* = \frac{\mu_i - \bar \mu}{\sigma^2(1-\rho)}+\frac{1}{n}
 \end{equation}
 provided that $y^i_*>0$ for all $i$. Notice that \eqref{eq:ifd2} agrees with \eqref{eq:ifd1} when $\rho=0$ and $\sigma_i =\sigma$. 
\begin{biology}{ of equation~\eqref{eq:ifd2}}
If environmental fluctuations have a sufficiently large variance $\sigma^2$, then ideal free dispersers visit all patches and spend more time in patches that support higher mean per-capita growth rates. Increasing the common spatial correlation $\rho$ results in ideal free dispersers spending more time in patches whose mean per-capita growth rate is greater than the average of the mean per-capita growth rates and less time  in other patches  (Fig.~\ref{fig:two}). When the spatial correlations are sufficiently large, it is no longer optimal to disperse to the patches with lower mean per-capita growth rates  ($\rho=0.5$ and $\rho=0.95$ in Fig.~\ref{fig:two}). \end{biology}
\end{example}

\section{The effect of constraints on dispersal}
\label{S:constraints}

While the ideal free patch distribution is a useful idealization to investigate how organisms should disperse in the absence of constraints, organisms in the natural world have limits on their ability to disperse and to collect and interpret environmental information.  Recall from
Section~\ref{S:ideal_free} that if the optimal patch distribution $y^*$
for an ideal free disperser is in the interior of the probability simplex
$\Delta$, then, loosely speaking, the ideal free disperser achieves the maximal
stochastic growth rate by using a strategy for which dispersal rate matrix
is of the form $D = \delta Q$, where $Q$ is any irreducible dispersal matrix with $(y^*)^T Q = 0$ and $\delta = \infty$.  At the
opposite extreme, if $y^*$ assigns all of its mass to a single patch,
then an ideal free disperser never leaves that single
most-favored patch.

To get a better understanding of how constraints on dispersal influence population growth, we consider dispersal  matrices of the form 
$D=\delta Q$,
where $\delta \ge 0$ and $Q$ is a fixed irreducible dispersal matrix $Q$ with a stationary distribution $\pi$ that is not necessarily the optimal patch distribution for
an ideal free disperser in the given environmental conditions.
We write $\chi(\delta)$ for the stochastic growth rate of the population as a function of the dispersal parameter $\delta$ 
and ask which choice of $\delta$ maximizes $\chi(\delta)$.  In particular,
we are interested in conditions under which some intermediate $\delta>0$ maximizes the stochastic growth rate $\chi(\delta)$.

We know from Proposition~\ref{prop:limit} that  
$\chi(\delta)$ approaches $\pi^T\mu - \frac{1}{2} \pi^T \Sigma \pi$ as $\delta\to\infty$.  We therefore set
$\chi(\infty) = \pi^T\mu - \frac{1}{2} \pi^T \Sigma \pi$.
On the other hand, if there is no dispersal ($\delta =0$), then 
$
\lim_{t\to\infty} \frac{1}{t} \log X_t^i = \mu_i -\frac{\sigma_i^2}{2}
$
with probability one whenever $X_0^i>0$, and 
so
$
 \lim_{t\to\infty}\frac{1}{t}\log S_t =\max_i \{\mu_i -\frac{\sigma_i^2}{2}\}
$
whenever $X_0^i>0$ for all $i$.  Hence, it is
reasonable to set $\chi(0)= \max_i \{\mu_i -\frac{\sigma_i^2}{2}\}$. The following result, which we prove in Appendix~\ref{apx:continuity},
implies that the function $\delta \mapsto \chi(\delta)$ is continuous
on $[0,\infty)$.

\begin{proposition}\label{prop:cont} The
function $\delta \mapsto \chi(\delta)$ is analytic on the
interval $(0,\infty)$ and continuous at the point $\delta=0$. 
\end{proposition}

One way to establish that $\chi(\delta)$
is maximized for an intermediate value of $\delta$ is to show that
$\chi(0) < \chi(\infty)$ and that $\chi(\delta) > \chi(\infty)$ for
all sufficiently large $\delta$.
The following theorem provides an asymptotic approximation for $\chi(\delta)$ when
$\delta$ is large that  allows us to check when the latter condition holds. 
We prove the theorem under the hypothesis that the 
dispersal matrix $Q$ is \emph{reversible} with respect to its stationary distribution $\pi$;
that is, that  $\pi_i Q_{ij} =\pi_j Q_{ji}$ for all $i,j$. Reversibility implies that at stationarity the Markov chain defined by $Q$ exhibits ``balanced dispersal in the absence of local demography.'' Namely,  if a large number of individuals are 
independently executing the equilibrium movement dynamics, then the rate at which  individuals move from patch $i$ to patch $j$  equals the rate at which  individuals move from patch $j$ to patch $i$.  We note that diffusive movement (that is, the matrix $Q$ is symmetric) 
and any form of movement along a one-dimensional landscape 
(that is, the matrix $Q$ is tridiagonal) are examples of reversible Markov chains. 
We provide a proof of the theorem in Appendix~\ref{apx:bigone}.
Corollary~\ref{cor:diagonalized} below, which we prove in 
Appendix~\ref{apx:diagonalization}, provides a more readily computable
expression for the asymptotics of the
stochastic growth rate under further assumptions.

\begin{theorem}\label{thm:bigone} 
Suppose that $Q$ is reversible with respect to its
stationary distribution $\pi$. Then, 
\begin{eqnarray}\label{eq:bigone}
\nonumber\chi({\delta}) &=&  
  \left( \mu^T \pi -  
    \frac{1}{2}  \pi^T \Sigma \pi \right)  +  
    \frac{1}{\delta}
    \biggl [ 
    ( \mu  - \Sigma \pi)^T \nu \\
  & &\qquad -   \frac{1}{2} \int_0^\infty  \mathrm{Tr}\left(   \exp(Q^Ts) \left( \mathrm{diag}(\pi) - \pi \pi^T  \right) \Sigma \left( \mathrm{diag}(\pi) - \pi \pi^T  \right)   \exp( Q s) \Sigma \right)  \, ds \biggr ]  + 
\mathrm{O}(\delta^{-\frac{5}{4}}) 
\end{eqnarray}
as $\delta \to \infty$, where $\nu$ is the unique vector satisfying $\bone^T\nu =0$ and $Q^T\nu= -\left( \mathrm{diag}(\pi) - \pi \pi^T  \right)(\mu-\Sigma\pi )$.
\end{theorem} 

When the dispersal matrix $D=\delta Q$ is consistent with ideal dispersal in the limit $\delta \to\infty$, equation~\eqref{eq:idf} implies that  $(\mu-\Sigma \pi)^T \nu=\lambda \bone^T \nu =0$. On the other hand, the proof of Theorem~\ref{thm:bigone} shows that 
\[
\int_0^\infty  \mathrm{Tr}\left(   \exp(Q^Ts) \left( \mathrm{diag}(\pi) - \pi \pi^T  \right) \Sigma \left( \mathrm{diag}(\pi) - \pi \pi^T  \right)   \exp( Q s) \Sigma \right)  \, ds=\mathrm{Tr} \left(\E[\bV_\infty \bV_\infty^T]  \Sigma \right)>0
\]
where $\bV_\infty$ is a Gaussian random vector. Hence, as expected, $\chi(\delta)$ is an increasing function for large $\delta$ when $\pi$ corresponds to the ideal free distribution associated with $\mu$ and $\Sigma$. However, when $\pi$ does not correspond to the ideal free distribution, $\chi(\delta)$ may be increasing or decreasing for large $\delta$ as we illustrate below.  

When $Q$ and $\Sigma$ commute, the asymptotic expression \eqref{eq:bigone}
for $\chi(\delta)$ simplifies a great deal.

\begin{corollary}
\label{cor:diagonalized}
Suppose that 
$Q$ is symmetric and $Q \Sigma = \Sigma Q$.
Let $\lambda_1 \le \ldots \le \lambda_{n-1} < \lambda_n=0$ be the eigenvalues of $Q$ 
with corresponding orthonormal eigenvectors $\xi_1, \ldots, \xi_n$.
Then, the eigenvalues $\theta_1, \ldots, \theta_n$ of $\Sigma$ 
can be ordered so that $\Sigma \xi_k = \theta_k \xi_k$, for each $1\le k \le n$,
and the approximation \eqref{eq:bigone} reduces to
\begin{equation} \label{eq:diagonalized}
\chi(\delta) 
  = 
  \left(\bar \mu - \frac{1}{2n} \theta_n\right)
  -
  \frac{1}{\delta n}
  \left[
    \sum_{k=1}^{n-1}
      \frac{1}{\lambda_k}
      \left(
      (\xi_k^T \mu)^2
      -
      \frac{1}{4 n} \theta_k^2
      \right)
  \right]
  + \mathrm{O}(\delta^{-5/4})
\end{equation}
as $\delta \to \infty$, where $\bar \mu = \frac{1}{n}\sum \mu_i$.
\end{corollary}

To illustrate the utility of this latter approximation, we develop more explicit formulas for three scenarios: diffusive movement in a landscape where all patches are equally connected (that is, a classic ``Levins'' style landscape \citep{levins-69b}), diffusive movement in a landscape 
consisting of a ring of patches, and diffusive movement in a landscape with multiple spatial scales (that is, a hierarchical Levins landscape). 

\begin{example}[Fully connected metapopulations with unbiased movement]   \label{ex:fully_connected}
Consider a population in which individuals disperse at the same per-capita rate $\delta/n$ between all pairs of patches.  Let $\sigma^2$ be the variance of the within patch fluctuations and $\rho$ be the correlation in these fluctuations between any pair of patches. Under these assumptions, the dispersal matrix is $Q =  J/n - I$ and the environmental covariance matrix is $\Sigma = (1- \rho)\sigma^2 I + \rho \sigma^2 J$, where recall that
$J = \bone \bone^T$ is the matrix of all ones.  Because $Q$ is symmetric, 
the stationary distribution of $Q$ is uniform; that is,
$\pi_1=\dots= \pi_n = \frac{1}{n}$.  Hence, 
in the absence of population growth there would be
equal numbers of individuals in each patch at large times. 

Because the matrices $I$ and $J$ commute, the matrices $Q$ and $\Sigma$ also commute.  Recall the notation of Corollary~\ref{cor:diagonalized}.
The eigenvector $\xi_n$
is $\frac{1}{\sqrt{n}} \bone$. If $\xi$ is any vector of length one orthogonal to $\xi_n$, 
then $J \xi = 0$, and so $Q \xi = - \xi$ and $\Sigma \xi = (1-\rho) \sigma^2 \xi$.
We may thus take $\xi_1, \ldots, \xi_{n-1}$ to be any orthonormal set of vectors orthogonal to $\xi_n$. 
Moreover,  $\lambda_1 = \cdots = \lambda_{n-1} = -1$
and $\theta_1 = \cdots = \theta_{n-1} = (1-\rho) \sigma^2$.

Now, $(\xi_n^T \mu)^2 = (1/n) \left( \sum_{k=1}^n \mu_k \right)^2 = n (\bar \mu)^2$,  and so Parseval's identity  implies that 
$\sum_{k=1}^{n-1} (\xi_k^T \mu)^2 
= \sum_{k=1}^n \mu_k^2 - n (\bar \mu)^2
= \mu^T \mu - n (\bar \mu)^2$.
Denote the variance of the vector $\mu$ by
\[
\mathrm{Var}[\mu] 
=\frac{1}{n} \mu^T \mu - (\bar \mu)^2 
= \frac{1}{n} \sum_{k=1}^{n-1} (\xi_k^T \mu)^2.
\]

Substituting these observations
 into equation \eqref{eq:diagonalized}, we get that
\begin{equation}\label{eq:cool} 
\chi(\delta) = \bar \mu -   
    \frac{\sigma^2}{2n} \left (1 + (n-1)\rho \right) + 
    \frac{1}{\delta} \left [ \mathrm{Var}[\mu]  - \frac{ (n-1)((1-\rho)\sigma^2)^2}{4n^2}  \right ]  + 
\mathrm{O}(\delta^{-\frac{5}{4}}). 
\end{equation}
Recall that for the special case of two uncorrelated patches with  $D_{12}=D_{21}=\delta/2$, $\mu_1=\mu_2=\mu$, and $\sigma_1=\sigma_2=\sigma$, we showed from our exact formula for $\chi(\delta)$  in the two patch case that
\[
\chi(\delta) \approx \mu - \frac{\sigma^2}{4} - \frac{1}{\delta} \frac{\sigma^4}{16}
\]
as $\delta \to \infty$, see \eqref{eq:twopatch_asymptotics}.
Hence, this approximation agrees with \eqref{eq:cool}.

Approximation~\eqref{eq:cool} implies that $\chi(\delta)$ is decreasing for large $\delta$ whenever
\begin{equation}\label{eq:dec}
  \frac{n}{\sqrt{n-1}} \sqrt{\mathrm{Var}[\mu]} > \frac{(1-\rho)\sigma^2}{2},
\end{equation}
and that $\chi(\delta)$ is increasing if the opposite inequality holds. We have remarked that, in general, an intermediate dispersal rate is
optimal when  $\chi(0) < \chi(\infty)$ and  
$\chi(\delta) > \chi(\infty)$ for
all sufficiently large $\delta$.  This will occur for
individuals in this diffusive dispersal regime when 
\begin{equation}\label{eq:biginfinity}
\frac{(1-\rho)\sigma^2}{2} > \frac{\max_i \mu_i -\bar \mu}{1-1/n}
\end{equation}
and \eqref{eq:dec} holds. In particular, when there are many patches (that is, $n\to\infty$), inequalities \eqref{eq:biginfinity} and \eqref{eq:dec} are both satisfied if  
\[
(1-\rho)\sigma^2/{2} > \max_i \mu_i -\bar \mu > 0.
\]
\begin{biology}{ of equations~\eqref{eq:dec} and ~\eqref{eq:biginfinity}}
Highly diffusive movement has a negative impact on population growth whenever there are sufficiently many patches and there is sufficient spatial variation in the mean per-capita growth rates. Alternatively, if there is no spatial variation in the mean per-capita rates and stochastic fluctuations are not perfectly correlated, then the population growth rate continually increases with higher dispersal rates. This latter observation is consistent with individuals being distributed equally across the landscape is the optimal patch distribution. In contrast, if there is some spatial variation in the mean per-capita growth rates and there are sufficiently large, but not perfectly correlated environmental fluctuations, then an intermediate dispersal rate maximizes the stochastic growth rate for diffusively dispersing populations. 
\end{biology}
\end{example}

%\subsection{Dispersal on groups}

In order to apply Corollary~\ref{cor:diagonalized}, we need to 
to simultaneously diagonalize the matrices $Q$ and $\Sigma$.
A situation in which this is possible and the resulting formulas
provide insight into biologically relevant scenarios is when
the dispersal mechanism and the covariance structure of the noise
both exhibit the symmetries of an underlying group.
Example~\ref{ex:fully_connected} above is a particular instance of
this situation.

More specifically, we suppose that the patches 
can be labeled with the elements of a finite group $G$ in such a way
that the migration rate $Q_{g,h}$ and environmental covariance $\Sigma_{g,h}$ between patches $g$ and $h$ 
both only depend on the ``displacement'' $gh^{-1}$
from $g$ to $h$ in $G$. That is, we assume
there exist functions $q$ and $s$ on $G$ such that
$Q_{gh} = q(gh^{-1})$ and $\Sigma_{gh} = s(gh^{-1})$.
For instance, if $G$ is the group of integers modulo $n$,
then the habitat has $n$ patches arranged in a circle, 
and the dispersal rate and environmental covariance between two patches
only depends on the distance between them, measured in steps around the circle.
We do not require
that the vector $\mu$ of mean per-capita growth rates
satisfies any symmetry conditions. 

The matrices $Q$ and $\Sigma$ will commute if  
$q$ and $s$ are {\em class functions},
that is, if $q(gh) = q(hg)$ and 
$s(gh) = s(hg)$
for all $g, h \in G$.  We assume this condition holds
from now on.  
Note that if $G$ is Abelian (that is,
the group operation is commutative), then
any function is a class function.  
%The group structure allows us to use representation theory
%to find the eigenvalues and eigenvectors of $Q$ and $\Sigma$
%and then to use Corollary~\ref{cor:diagonalized}.

\subsection{Background on group representations}

We now record a few facts about representation theory,
the tool that will enable us to find the eigenvalues and eigenvectors
of $Q$ and $\Sigma$,
resulting in Theorem~\ref{thm:character_rep}.
We refer readers interested in more detail to \citep{serre-77,diaconis-88},
while readers interested in less mathematical detail 
may skip directly to Examples \ref{ex:circle} and \ref{ex:multiscale} without loss of continuity.

A {\em unitary representation} of a group $G$ is a homomorphism $\rho$
from $G$ into the group of $d_\rho \times d_\rho$ unitary matrices,
where $d_\rho$ is called the {\em degree} of the
representation.
Two representations $\rho'$ and $\rho''$ are {\em equivalent} 
if there exists a unitary matrix $U$ such that
$\rho''(g) = U \rho'(g) U^{-1}$ for all $g \in G$.  
A representation $\rho'$ is {\em irreducible} if it is not equivalent
to some representation $\rho''$ for which $\rho''(g)$ is
of the same block diagonal form for all $g \in G$.
A finite group has a finite set of inequivalent, irreducible, unitary representations, which we denote by $\hat G$.  
The simplest representation is the {\em trivial representation} 
$\rho_{\mbox{tr}}$ of degree one, for which $\rho_{\mbox{tr}}(g) = 1$ for all $g$.

For a simple example that we will return to, let $G=\Z_n$, the group of integers modulo $n$.
Since $\Z_n$ is Abelian, all the irreducible representations are one-dimensional 
($d_\rho=1$ for all $\rho \in \hat G$),
and are of the form $\rho^{(m)}(j) = \exp(2\pi i m j/n)$, 
so that $\hat G  = \{ \rho^{(m)}_m : 0 \le m \le n-1\}$.

The matrix entries of irreducible representations are {\em orthogonal}: for $\rho', \rho'' \in \hat G$,  
\begin{equation}
\label{eq:matrix_entry_orthog}
\sum_{g \in G} \rho'_{ij}(g) \rho''_{k \ell}(g)^* 
= 
\begin{cases}
\frac{\# G}{ d_\rho}, & \text{if $\rho' = \rho''$ and $(i,j) = (k,\ell)$},\\
0, & \text{otherwise,}
\end{cases}
\end{equation}
where $z^*$ denotes the complex conjugate of 
a complex number $z$, and $\# G$ is the number of elements of $G$.
%(This implies that the degrees satisfy $\sum_\kappa d_\rho^2 = \# G$.)

The {\em Fourier transform} of a function $f: G \to \mathbb{C}$ is a function $\hat {f}$ on $\hat G$ defined by
\begin{equation} \label{eqn:ft_defn}
    \hat f(\rho) := \sum_{g\in G} f(g) \rho(g) \quad \mbox{for} \; \rho \in \hat G.
\end{equation}
Note that $\hat f(\rho)$ is a $d_\rho \times d_\rho$ matrix.
It follows from the orthogonality properties of the matrix entries of
the irreducible representations recorded above that
the Fourier transform may be inverted, giving $f$ explicitly as the linear combination of matrix entries of $\hat f$.
The inversion formula is
\[
    f(g) = \frac{1}{\# G} \sum_{\rho \in \hat G} d_\rho \Tr \left( \rho(g^{-1}) \hat f(\rho) \right).
\]

For $G=\Z_n$, this is the familiar {\em discrete Fourier transform},
for which orthogonality of matrix entries is the fact that $(1/n) \sum_{j=0}^{n-1} \exp(2\pi i j(\ell-m)/n) = \delta_{\ell m}$.
The transform is given by
$\hat f(\rho^{(m)}) = \sum_{k=0}^{n-1} f(k) \exp( 2\pi i m k/n)$ for $0 \le m \le n-1$,
and
$f(k) = (1/n) \sum_{m=0}^{n-1} \hat f(\rho^{(m)}) \exp( -2\pi i m k/n)$.
The trivial character is $\trchar=\rho^{(0)}$.

Associated with a representation $\rho \in \hat G$ is its {\em character} $\kappa$, defined by $\kappa(g) := \Tr \rho(g)$.  We write $\tilde G$
for the set of characters of irreducible representations.  
The characters are class functions, and
form an orthogonal basis for the subspace of class functions on $G$
and all have the same norm: $\sum_{g \in G} |\kappa(g)|^2 = \# G$, 
where $|z| = \sqrt{z z^*}$ is the modulus of the complex number $z$.
For $\rho \in \hat G$ with character $\kappa \in \tilde G$,
the Fourier transform of a class function $f$ satisfies
\[
\hat f(\rho)= \frac{1}{d_\rho} \tilde f(\kappa) I
\]
where $I$ is the $d_\rho \times d_\rho$ identity matrix and
\begin{equation} \label{eqn:clas_ft_defn}
    \tilde f(\kappa) := \sum_{g \in G} f(g) \kappa(g).
\end{equation}
% Note that $\hat f(\kappa) = \Tr( \ft f(\rho) )$. 
% Since we mostly use $\hat f$, this should not lead to confusion.
%For every $\kappa \in G^*$, $\kappa(\id) = d_\kappa$ and $\kappa(g^{-1}) = \kappa(g)^*$,
%and the class functions $\kappa \in G^*$ form an orthogonal basis of the space of class functions on $G$.
%The inverse Fourier transform for class functions is
Consequently,
\begin{equation}
\label{eq:expansion_class_function}
  f(g) = \frac{1}{\# G} \sum_{\kappa \in \tilde G} \kappa(g)^* \tilde f(\kappa).
\end{equation}

As noted above, if $G=\Z_n$ then all irreducible representations are one-dimensional,
so in this case we may identify the characters with the irreducible representations, $\hat G =\tilde G$.
Since $\Z_n$ is Abelian, all functions on $\Z_n$ are class functions,
so that the two Fourier transforms \eqref{eqn:ft_defn} and \eqref{eqn:clas_ft_defn} are equal.

Finally, given a function $f$ on $G$ and character $\kappa$, define
% that corresponds to a representation $\rho$,
% denote by $\|f\|_\kappa$ the norm of the projection of $f$ onto the subspace $\mathbb{C}^G$ spanned by the matrix entries of $\rho$.
% It is shown in Appendix~\ref{apx:character_rep} that this can be written as
% \|f\|_\kappa^2 := \frac{d_\rho}{\# G} \sum_{i,j=1}^{d_\rho} \left| \sum_{g \in G} \rho_{ij}(g) f(g) \right|^2.
\[
\|f\|_\kappa^2 := \frac{d_\rho}{\# G} \sum_{g,h\in G} \kappa(gh^{-1}) f(g) f(h)^* .
\]

The following theorem is proved in Appendix~\ref{apx:character_rep}.

\begin{theorem}\label{thm:character_rep}
Suppose that the $n$ 
patches are labeled by a finite group $G$ in such a way that
$Q_{gh} = q(gh^{-1})$ and $\Sigma_{gh} = s(gh^{-1})$,
where $q$ and $s$ are class functions.
Suppose further that $q(g) = q(g^{-1})$, $g \in G$, so that
the matrix $Q$ is symmetric. Let $\bar \mu = \frac{1}{\# G}\sum_{g\in G} \mu(g)$ and $\bar s = \frac{1}{\# G}\sum_{g\in G} s(g)$.
Then,
\begin{equation}\label{eq:character_rep}
\chi(\delta) =
\left(\bar \mu - \frac12 \bar s \right)
- \frac{1}{\delta n} \sum_{\kappa \in \tilde G \setminus \{\trchar\} } 
    \frac{ d_\kappa }{ \tilde{q}(\kappa) } 
      \left( \| \mu \|_\kappa^2 
         - \frac{1}{4 n} \tilde{s}(\kappa)^2 \right) 
  + O(\delta^{-5/4})
\end{equation}
as $\delta \to \infty$.
Furthermore, $\tilde{q}(\kappa) < 0$ for all 
$\kappa \in \tilde G \setminus \{\trchar\}$.
\end{theorem}

Roughly speaking, this expression tells us about the respective roles of variance of patch quality ($\mu$) and covariance of environmental noise ($s$).
The fact that $\tilde{q}(\kappa)$ is negative for all $\kappa$ leads to the following.
\begin{biology}{ of equation~\eqref{eq:character_rep}}
    If variability in patch quality at a certain scale is larger 
    than the correlation in environmental noise at that scale,
    in a sense made precise above,
    then the stochastic growth rate decreases with increasing dispersal rates  at that scale.
    Conversely, if environmental noise is strongly correlated between patches
    and the mean patch quality is similar, then more dispersal is expected to be better.
    The relevant sense of ``at that scale'' is in the sense of the Fourier transform,
    analogous to the ``frequency domain'' in Fourier analysis.
\end{biology}

\begin{example}[Circle of Patches]
  \label{ex:circle}
Suppose that the $n$ patches of a habitat are arranged in a circle
and are labeled by $\Z_n = \{0,1,\ldots,n-1\}$, the group of
integers modulo $n$ with identity element $0$.
As reviewed above, the Fourier transform is the familiar discrete Fourier transform.
% Because $\Z_n$ is Abelian, any function is a class function.
% The irreducible representations of $\Z_n$ are all one-dimensional 
% (that is, $d_\rho=1$ for all $\rho \in \hat G$), and hence
% an irreducible representation can be identified with its character.
% The characters are of the form 
% $j \mapsto \kappa_m(j) = \exp\left(2 \pi i  m j  / n \right)$ for $0 \le m \le n-1$,
% and $\kappa_0$ is the trivial character $\trchar$. Given a function $\Z_n \to \C$, its Fourier transform 
% is given by 
% \[ 
% \tilde 
% f(\kappa_m) 
% = \sum_{j=0}^{n-1} f(j) \exp\left(2 \pi i m j \pi / n \right) \quad \mbox{for} \; 0 \le m \le n-1. 
% \]

If we assume that individuals disperse only to neighboring patches and these dispersal rates are equal, then  
$q(1) = q(n-1) = 1/2$, $q(0) = -1$ and $q(2) = \ldots = q(n-2) = 0$. Assume the environmental noise is independent between patches and has variance $\sigma^2$ i.e. $s(0)=\sigma^2$ and $0 = s(1) = \ldots = s(n-1)$.
% here variance is $a$ , the covariance of environmental noise in neighboring patches is $s(1)=s(n-1)=b/2$ with $|b|/2<a$, and the environmental noise in 
%non-neighboring patches are independent i.e. $0 = s(2) = \ldots = s(n-2)$.  
Finally, suppose  that patch quality as measured by the average per-capita growth rates is spatially periodic, so that $\mu(k) = \bar \mu + c \; \cos(2\pi k \ell / n)$ for some $c>0$,  $\bar \mu$, and $1\le \ell <n/2$.

Under this set of assumptions, we can compute that for $m \ne 0$,
$\tilde q(m) = \cos(2 \pi m/n) - 1 $ and $\tilde s(m) = \sigma^2$.
%$\tilde q(m) = \cos(2 \pi m/n) - 1 $ and $\tilde s(m) = a + b \cos(2 \pi m/n)$.
Furthermore,  $\|\mu\|^2_{\kappa_\ell} = \|\mu\|^2_{\kappa_{n-\ell}} = n c^2/4$  and $\|\mu\|^2_{\kappa_m} = 0$ otherwise.
From these computations, Theorem~\ref{thm:character_rep} implies that
\[
\chi(\delta) \approx \bar \mu -\sigma^2/2 
%\left(\bar \mu - \frac{1}{2} \frac{a + b}{n}\right)
-
\frac{1}{\delta n}
\left(
\frac{ n c^2 }{ 2 (\cos(2\pi \ell/n) - 1) }
-
\sum_{m=1}^{n-1}
\frac{\sigma^2}{4 n (\cos(2 \pi m/n) - 1) } 
%\frac{(a + b \cos(2 \pi m/n))^2}{4 n (\cos(2 \pi m/n) - 1) } 
\right)
\]
for large $\delta$. Using the identity 
$\sum_{k=1}^{n-1} ( 1-\cos(2\pi k / n) )^{-1} = (n^2-1)/6$
(see equation 1.381.1 in \citet{gradshteyn-07}'s table of integrals and series), this approximation simplifies to
\begin{equation}\label{eq:circle}
\chi(\delta) \approx 
\bar \mu - \sigma^2/2
%\left(\bar \mu - \frac{1}{2} \frac{a + b}{n}\right)
        + \frac{1}{4 \delta n^2}
        \left(
        \frac{ 2 n^2 c^2 }{ 1 - \cos(2\pi \ell/n) } -
         \frac{1}{6} (n^2-1) \sigma^4  
		%\left\{ \frac{1}{6} (n^2-1) a^2 + \frac{1}{3} n (n-5) ab + (n^2 -6n + 11) b^2 
        %\right\} 
        \right).
\end{equation}

Since $\chi(0)=\bar \mu +c-\sigma^2/2$, high dispersal is better than no dispersal if 
%$\chi(\infty)-\chi(0) =  a (1-n^{-1})/2 - b/2n - c > 0$.
$\chi(\infty)-\chi(0) =  \sigma^2 (1-1/n)/2 - c > 0$. 
When the number of patches is sufficiently large,  this inequality implies that highly dispersive populations grow faster than sedentary populations provided that the temporal variation is sufficiently greater than the spatial variation in per-capita growth rates i.e. $\sigma^2>2c$. On the other hand, $\chi(\delta)$ is decreasing for large $\delta$ 
if the coefficient of $1/\delta$ is positive i.e. 
\[
 4c^2 > 
\frac{1}{3} (1-\cos(2\pi \ell/n)) (1-n^{-2}) \sigma^4 .
% \frac{1}{3} (1-\cos(2\pi \ell/n)) \left( (1-n^{-2}) a^2 + 2(1-5n^{-1}) ab + 6(1 -6n^{-1} + 11n^{-2}) b^2 \right) .
\]
Hence, if $\ell/n$ is small enough, then $\chi(\delta)$ is decreasing for large $\delta$. 
\begin{biology}{ of equation~\eqref{eq:circle}}
    In a circular habitat with nearest-neighbor dispersal and sinusoidally varying patch quality,
intermediate dispersal rates maximize the stochastic growth rate provided that spatial heterogeneity occurs on a short scale (i.e. $\ell/n$ sufficiently small) 
and temporal variability is sufficiently large.
\end{biology}
\end{example}

\begin{example}[Multi-scale patches]
  \label{ex:multiscale}
Suppose now that our organism lives in a hierarchically structured habitat.
For example, individuals might live on bushes, the bushes grow around the edges of clearings, 
and the clearings are scattered across an archipelago of islands.
We label each bush with an ordered triple recording on which island, in which clearing, and in what 
bush around the clearing it lives,
so that for instance $(2,1,4)$ denotes the fourth bush in the first clearing of the second island.
To make the mathematical picture a pretty one, 
we suppose that each of the $I$ islands 
has the same number $C$ of clearings
and each clearing has the same number $B$ of bushes.
This enables us identify the habitat structure with the group $\Z_I \otimes \Z_C \otimes \Z_B$,
where, as above, $\Z_m$ is the group of integers modulo $m$.
We will get particularly simple and interpretable results if
we also assume that dispersal rates and environmental covariances
only depend on the scale at which the movement occurs --
between bushes, clearings, or islands.

Although it requires imaginative work to find examples 
with many more scales than this
(do the organism's fleas have fleas?)
it does not cost us anything to work in greater generality.
Suppose, then, that the patches in the habitat 
are labeled with the group $G = G_1 \otimes \cdots \otimes G_k$,
where $G_j = \Z_{n_j}$ for $1 \le j \le k$.

Thus, one patch is labeled with the identity element
$\id_G = (\id_1, \ldots, \id_k)$
and every other patch is labeled by the displacement 
required to get there from $\id_G$.
The later coordinates are understood to be at finer ``scales'',
so that if $g_i = h_i$ for all $1 \le i \le j-1$,
then $g$ and $h$ represent patches in the same \emph{metapatch at scale $j$}.
For instance, in our example above, the archipelago of islands
is the single metapatch at scale $1$ and the metapatches
at scales $2$ and $3$ are, respectively, 
the islands and the clearings. 
We label the metapatches at scale $r$ 
with the set $Z_r := \{g \in G : g_r = \id_r, \ldots g_k = \id_k\}$,
with the convention that $Z_{k+1} := G$.
Because a label $g=(g_1,\ldots,g_k) \in G$ represents displacement,
the coordinate of the leftmost non-identity element of $g$,
denoted by
\[
    \ell(g) := \min\{ j : g_j \neq \id_j \}  \mbox{ and }\ell(\id_G)=k+1,
\]
tells us the scale on which the motion occurs: 
$g \in G$ corresponds
to a displacement that moves between patches within
the same metapatch at scale $\ell(g)$ but moves from a patch
within a metapatch at scale $\ell(g)+1$ to a patch within some
other metapatch at that scale.  Note that $1 \le \ell(g) \le k+1$.

We assume that the dispersal rate 
and the environmental covariance between two patches 
only depends on the scale of the
displacement necessary to move between the two patches.
That is, we suppose there are numbers 
$q_1, \ldots, q_{k+1}$
and
$s_1, \ldots, s_{k+1}$
such that 
$q(g) = q_{\ell(g)}$
and
$s(g) = s_{\ell(g)}$.

In Appendix~\ref{apx:multiscale} we show that the Fourier transforms appearing in Theorem~\ref{thm:character_rep}
depend on  the following quantities. 
Let $N_r  := \# Z_r = \prod_{j=1}^{r-1} n_j$ 
be the number of metapatches at scale 
$r$.  Write $\bar Z_r := \{ g \in G : g_j = \id_j, \; j \le r \}$
for the subgroup of displacements that move from one patch to another
within the same metapatch at scale $r+1$ and set
$\bar N_r := \# \bar Z_r = \prod_{j={r+1}}^{k} n_j$.
%We can identify $Z_r$ with $G_1 \otimes \cdots \otimes G_{r-1}$ and
%$\bar Z_r$ with $G_{r+1} \otimes \cdots \otimes G_k$.
Set
\[
    v_\mu(r) := \frac{1}{N_r}  \sum_{g \in Z_r} \left(
        \frac{1}{n_r} \sum_{h \in G_r} \left( \frac{1}{\bar N_r} \sum_{z \in \bar Z_r} \mu(g \gtimes h \gtimes z) \right)^2
                    - \left( \frac{1}{n_r} \sum_{h \in G_r} \frac{1}{\bar N_r} \sum_{z \in \bar Z_r} \mu(g \gtimes h \gtimes z) \right)^2 \right).
\]
We can interpret this quantity as follows. There are
$N_r$ metapatches at scale $r$.  
Each one has within it $n_r$ metapatches at scale $r+1$.
First, compute the average of $\mu$ over all the patches
within each metapatch
at scale $r+1$, then compute the variance 
of these averages within each metapatch at scale $r$,
and finally average these variances across all the
 metapatches at scale $r$ to produce $v_\mu(r)$.  Thus,
$v_\mu(r)$ measures the variability in $\mu$ that can
be attributed to scale $r+1$.
Set
\[
    \tilde s(r) = \sum_{\ell=r}^k (s_{\ell+1} - s_\ell) \bar N_\ell
\]
and
\[
    \tilde q(r) = - \sum_{\ell=1}^{r} q_\ell (\bar N_{\ell-1} - \bar N_\ell) - q_{r} \bar N_{r}. \\
\]
%Because there are $\bar N_{\ell-1} - \bar N_\ell$ patches that differ from a given patch at the $\ell^\mathrm{th}$ level,
%and $q_\ell$ is the migration rate to such patches,
%the total rate of migration to patches differing at least on the $r^\mathrm{th}$ level is
%$\sum_{\ell=1}^{r} q_\ell (\bar N_{\ell-1} - \bar N_\ell)$.

The following result agrees with equation \eqref{eq:cool}, which
describes the special case where there is a single scale.
\begin{theorem} \label{thm:multiscale}
For a habitat with the above multi-scale structure, 
equation \eqref{eq:bigone} reduces to
\begin{equation}\label{eq:multiscale}
    \chi(\delta) 
    = \left(\bar \mu - \frac12 \bar s \right)
        - \frac{1}{\delta} \sum_{r=1}^k \frac{ 1 }{ \tilde{q}(r) } \left( v_\mu(r) - \frac{N_{r+1}-N_r}{4 N_{k+1}^2} \tilde{s}(r)^2 \right) 
                + O(\delta^{-5/4})
\end{equation}
as $\delta \to \infty$. Furthermore, $\tilde{q}(r)<0$ for all $1\le r\le k$. 
\end{theorem}

Note that if $s_\ell$ increases with $\ell$ 
(that is, two patches within the same 
metapatch have a higher environmental covariance than two
patches in different metapatches at that scale), then $\tilde s(r)$ decreases 
with $r$.  Also,  if $q_\ell$ increases with $\ell$ 
(that is, there is a higher rate for dispersing to a patch within the same metapatch at
some scale than to a patch in another metapatch at that scale), 
then $\tilde q(r)$ is negative and decreases with $r$.
Using these observations, we may read off several things from 
\eqref{eq:multiscale}.

First, consider a simple example
with a fixed, large number $n$ of patches distributed among a variable number of islands.
Now $k=2$, and let the number of islands $n_1=1/\alpha$, with $\alpha \ge 1$, 
so that the number of patches on each island is $n_2=\alpha n$.
In this case, $N_1=1$, $N_2=1/\alpha$, and $N_3=n$, while $\bar N_0 = n$, $\bar N_1=\alpha n$, and $\bar N_2=1$,
so \eqref{eq:multiscale} reads
\begin{eqnarray}\label{eq:multiscaleA}
\chi(\delta) &\approx &(\bar \mu - \frac{1}{2} \bar s) - \frac{1}{\delta} \left\{ 
  - \frac{v_\mu(1)}{q_1 n} + \frac{(1-\alpha)((s_3-s_2)+\alpha n (s_2-s_1))^2}{\alpha q_1 n^2} 
  - \frac{ v_\mu(2) - (\alpha n - 1) (s_3-s_2)^2 }{ \alpha n^2 (q_2 \alpha + q_1 (1-\alpha)) }
  \right\} \\
 \nonumber &= &(\bar \mu - \frac{1}{2} \bar s) - \frac{ \alpha(1-\alpha) (s_2-s_1)^2 }{ \delta q_1 } + O(n^{-1}) .
\end{eqnarray}
The effect of higher dispersal depends on 
the difference in covariances between patches on the same island and on different islands,
and on the number of islands.
\begin{biology}{ of equation~\eqref{eq:multiscaleA}}
    If a sufficiently large number of patches are distributed equally across a number of islands,
    then for a given dispersal pattern,
    the stochastic growth rate increases with the dispersal rate (at high levels of dispersal).
    This effect is strongest if there are only two islands (i.e. $\alpha=1/2$).
\end{biology}

Secondly, imagine a fixed ensemble of patches with varying mean per-capita
growth rates and consider the following two possibilities
for assignment of these patches to metapatches at scale $2$
(the islands in our bush-clearing-island example). 
One possibility is that some islands are assigned patches
that are primarily of high quality, whereas other islands
are mostly assigned poor patches.
The other possibility is that 
patches of different quality are evenly spread across the islands,
with the range of quality within an island similar to the range of quality between islands.
In the first case, the variance across islands of within-island means is comparable to the variance across all patches,
so $v_\mu(1) \approx v_\mu(k)$. 
In the second case, the within-island means are approximately constant,
so that $v_\mu(1)$ will be small.
Therefore, since $\tilde q(r)$ is negative for all $r$, 
having local positive association of $\mu$ at nearby patches leads to higher stochastic growth rates,
at least for large enough values of the dispersal parameter $\delta$.
\begin{biology}{ of equation~\eqref{eq:multiscale}}
    All other things being equal,
    the species will do better if the good habitat is concentrated on particular islands, 
    rather than spread out across many.
\end{biology}

Finally, we can observe that adding new scales of metapatch 
may change the situation from one in which $\chi(\delta)$ is maximal at high values of the dispersal parameter $\delta$ 
to one in which $\chi(\delta)$ is maximal at intermediate values of $\delta$,
or vice-versa.
If $n_1 = 1$, then $\tilde s(1)$ and $v_\mu(1)$ are both zero,
and changing $n_1$ (for example, going from one to several
islands in our example) will increase $\tilde s(1)$.
Changing $n_1$ will also add the quantity $-q_1 (n_1-1)\bar N_1$ to all values of $\tilde q(r)$.
The result of this could be to change the sign of the coefficient of $\frac{1}{\delta}$ in \eqref{eq:bigone}.
\begin{biology}{ of equation~\eqref{eq:multiscale}}
    The optimal level of dispersal for a subpopulation, and the growth rate at that level of dispersal,
    may differ drastically depending on whether it is connected (or connectable) by dispersal to other subpopulations.
\end{biology}

\end{example}

\section{Discussion}
\label{S:discussion}

Classical ecology theory predicts that environmental stochasticity increases extinction risk by reducing the long term per-capita growth rate of populations~\citep{may-75,turelli-78}. For sedentary populations in a spatially homogeneous yet temporally variable environment, a simple model of their growth is given by the stochastic differential equation
$dZ_t = \mu Z_t dt+ \sigma Z_t dB_t$, where $B$ is a standard Brownian motion. The stochastic  growth rate for such populations equals $\mu -\frac{\sigma^2}{2}$;  the reduction in the growth rate is proportional to the infinitesimal variance of the noise. Here, we show  that a similar expression describes the growth of populations dispersing in spatially and temporally heterogeneous environments. More specifically, 
if average per-capita growth rate in patch $i$ is $\mu_i$
and the infinitesimal spatial covariance between environmental noise
in patches $i$ and $j$ is $\sigma_{ij}$, then
the stochastic growth rate equals the average of the mean per-capita growth rate $\langle \mu \rangle=\sum_i \mu_i \E[Y_\infty^i]$ experienced by the population when the 
proportions of the population in the various patches have
reached equilibrium
minus half of the average temporal variation 
$\langle \sigma^2 \rangle=\E[\sum_{i,j} \sigma_{ij} Y_\infty^i Y_\infty^j]$ 
experienced by the population in equilibrium. 
The law of $\bY_\infty$, the random equilibrium spatial distribution of the population which provides the weights in these averages,  is determined by interactions between spatial heterogeneity in mean per-capita growth rates, the infinitesimal spatial covariances of the environmental noise, and population movement patterns.  To investigate how these interactions
effect the stochastic growth rate, we derived analytic expressions for the law of $\bY_\infty$,  determined what choice of  dispersal mechanisms resulted
in optimal stochastic growth rates for a freely dispersing population, 
and considered the consequences on the stochastic growth rate
of limiting the population to a fixed dispersal mechanism. As we now discuss, these analytic results provide fundamental insights into ``ideal free'' movement in the face of uncertainty,  the persistence of coupled sink populations, the evolution of dispersal rates, and the single large or several small (SLOSS) debate in conservation biology. 

In spatially heterogeneous environments, ``ideal free'' individuals disperse to the patch or patches that maximize their long term per-capita growth rate~\citep{fretwell-lucas-70,harper-82,oksanen-etal-95,vanbaalen-sabelis-99,amnat-00,prsb-06,siap-06,cantrell-etal-07}.  In the absence of environmental stochasticity and density-dependent feedbacks, ideal free dispersers only select the patches supporting the highest per-capita growth rate.  Here, we show that uncertainty due to environmental stochasticity can overturn this prediction. Provided environmental stochasticity is sufficiently strong and spatial correlations are sufficiently weak, equation \eqref{eq:opt3} implies that ideal free dispersers occupy all patches despite spatial variation in the local stochastic growth rates $\mu_i - \sigma_i^2/2$.  Intuitively, by spending time in multiple patches, including  those  that in isolation support lower stochastic growth rates, individuals reduce the net environmental variation $\langle \sigma^2 \rangle$ they experience and, thereby, increase their stochastic growth rate. Hence, dispersing to lower quality patches is a form of bet-hedging against environmental uncertainty~\citep{slatkin-74,philippi-seger-89,wilbur-rudolf-06}. When environmental fluctuations in higher quality patches are sufficiently strong, this spatial bet-hedging can result in ideal free dispersers occupying sink patches; patches that are unable in the absence of immigration to sustain a population. This latter prediction is consistent with Holt's analysis of a discrete-time two patch model~\citep{holt-97}. Spatial correlations in environmental fluctuations, however, can disrupt  spatial bet-hedging. Movement between patches exhibiting strongly covarying environmental fluctuations has little effect on the net environmental variation $\langle \sigma^2 \rangle$ experienced by individuals and, therefore, movement to lower quality patches may confer little or no advantage to individuals. Indeed, when the spatial covariation is sufficiently strong, ideal free dispersers only occupy patches with the highest local stochastic growth rates $\mu_i -\sigma_i^2/2$, similar to the case of deterministic environments.  In deterministic environments, density dependent feedbacks can result in ideal-free dispersers occupying multiple patches including sink patches~\citep{fretwell-lucas-70,cantrell-etal-07,holt-mcpeek-96}. Our results show that even density-independent processes can result in populations occupying multiple patches. However, both of these processes are likely to play important roles in the evolution of patch selection.

A sink population is a local population that is sustained by immigration~\citep{holt-85,pulliam-88,dias-96}. Removing immigration results in a steady decline to extinction. In contrast, source populations persist in the absence of immigration. Empirical studies have shown that landscapes often partition into mosaics of source and sink populations~\citep{murphy-01,kreuzer-huntly-03,re-05}. For discrete-time two-patch models, \citet{jansen-yoshimura-98} showed, quite surprisingly, that sink populations coupled by dispersal can persist, a prediction supported by recent empirical studies with protozoan populations~\citep{matthews-gonzalez-07} and extended to discrete-time multi-patch models~\citep{roy-etal-05,prsb-10}. Here, we show a similar phenomena occurs for populations experiencing continuous temporal fluctuations. For example, if the stochastic growth rates in all patches equal $\mu-\sigma^2/2$ and the spatial correlation between patches is $\rho$, then equations \eqref{eq:lyapunov} and \eqref{eq:ifd2} imply that populations dispersing freely between $n$ patches persist whenever $\mu - ((n-1)\rho+1)\sigma^2/2n>0$. Hence, ideal free movement mediates persistence whenever local environmental fluctuations produce sink populations (i.e., $\sigma^2/2>\mu>0$), environmental fluctuations aren't fully spatially correlated (i.e. $\rho<2\mu/\sigma^2$) and there are sufficiently many patches (i.e., $n>((1-\rho)\sigma^2)/(2\mu-\rho\sigma^2)$). This latter expression for the necessary number of patches to mediate persistence is an exact, continuous time counterpart to an approximation by \citet{bascompte-etal-02} for discrete time models.  When two patches are sufficient to mediate persistence,  equation \eqref{twopatch}  reveals that there is a critical dispersal rate below which the population is extinction prone and above which it persists. Our high dispersal approximation (see equation \eqref{eq:cool} with $\mathrm{Var}[\mu]=0$) suggests this dispersal threshold also exists for an arbitrary number of patches.

 While ideal free movement corresponds to the optimal dispersal strategy for species without any constraints on their movement or their ability to collect information, many organisms experience these constraints. For instance, in the absence of information about environmental conditions in other patches, individuals may move randomly between patches,  in which case the rate of movement (rather than the pattern of movement) is subject to natural selection~\citep{hastings-83,levin-etal-84,mcpeek-holt-92,holt-mcpeek-96,dockery-etal-98,hutson-etal-03,siap-06}. When density-dependent feedbacks are weak and certain symmetry assumptions are met, our high dispersal approximation in \eqref{eq:diagonalized} implies there is selection for higher dispersal rates whenever
\begin{equation}\label{fast}
\sum_{k=1}^{n-1}
\frac{1}{|\lambda_k|}\frac{1}{4 n} \theta_k^2
>
\sum_{k=1}^{n-1}
\frac{1}{|\lambda_k|}
(\xi_k^T \mu)^2
\end{equation}
 where, recall, $\lambda_k<0$, $\xi_k$ are the eigenvalues/vectors of the dispersal matrix, $\mu$ is the vector of per-capita growth rates, and $\theta_k$ are the eigenvalues of the covariance matrix for the environmental noise. Roughly speaking, equation~\eqref{fast} asserts that if temporal variation (averaged in the appropriate manner) exceeds spatial variation, then there is selection for faster dispersers; a prediction consistent with the general consensus of earlier  studies~\citep{levin-etal-84,mcpeek-holt-92,hutson-etal-03}. More specifically, in the highly symmetric case where the temporal variation in all patches equals $\sigma^2$ and the spatial correlation between patches is $\rho$, equation~\eqref{fast} simplifies to 
 \begin{equation}\label{fast2}
 \frac{(1-\rho)\sigma^2}{2}>  \frac{n}{\sqrt{n-1}} \sqrt{\mathrm{Var}[\mu]} ,
\end{equation} 
in which case lower spatial correlations and larger number of patches also facilitate selection for faster dispersers. Another important constraint influencing the evolution of dispersal are travel costs that reduce fitness of dispersing individuals. While the effect of these costs have been investigated for deterministic models~\citep{deangelis-etal-11}, it remains to be seen how these traveling costs interact with environmental stochasticity in determining optimal dispersal strategies.  
  
Previous studies have shown that  spatial heterogeneity in per-capita growth rates increases the net population growth rate for deterministic models with diffusive movement~\citep{adler-92,amnat-09b}. Intuitively, spatial heterogeneity provides patches with higher per-capita growth rates that boost the population growth rate, a boost that gets diluted at higher dispersal rates. Our high dispersal approximation \eqref{eq:diagonalized} shows that this boost also occurs in temporally heterogeneous environments, i.e. the correction term $-\sum_{k=1}^{n-1}\frac{1}{\lambda_k}(\xi_k^T \mu)^2/\delta$ is positive. More importantly, the multiscale version of this correction term \eqref{eq:multiscale} implies this boost is larger when the variation in the per-capita growth rates occurs at multiple spatial scales. For example, for insects living on plants in meadows on islands, the largest boost  occurs when the higher quality plants (i.e. the plants supporting the largest $\mu_i$ values) occur on the same island in the same meadow.  This analytic conclusion is consistent with numerical simulations showing that habitat fragmentation (e.g. distributing high quality plants more evenly across islands and meadows)  increases extinction risk~\citep{fahrig-97,fahrig-02}. Intuitively, spatial aggregation of higher quality patches increases the chance of individuals dispersing away from a high quality patch arriving in another high quality patch. Even without spatial variation in per-capita growth rates, equation \eqref{eq:multiscale} implies that strong spatial aggregation of patches maximizes stochastic growth rates for dispersive populations living in environments where temporal correlations decrease with spatial scale. This finding promotes the view that a single large (SL) reserve is a better for conservation  than several small (SS)  reserves. This finding is consistent with many arguments in the SLOSS debate~\citep{diamond-75,wilcox-murphy-85,gilpin-88,cantrell-cosner-89,cantrell-cosner-91}. For example, using reaction-diffusion equations, \citet{cantrell-cosner-91} found that even in deterministic environments ``[it] is  better  for  a population to have  a  few  large  regions  of  favorable  habitat  than  a  great  many  small  ones 
closely  intermingled  with  unfavorable  regions.'' However, our results run contrary to a numerical simulation study of \citet{quinn-hastings-87} that, unlike ours, applies to sedentary populations experiencing independent environments~\citep{gilpin-88}.  

 While our work provides a diversity of analytical insights into the interactive effects of temporal variability, spatial heterogeneity, and movement on long-term population growth, many challenges remain. Most notably, are there analytic approximations for relatively sedentary populations? What effect do correlations in the temporal fluctuations have on the stochastic growth rate? Can the explicit formulas for stochastic growth rates in two patch environments be extended to special classes of higher dimensional models? Can one extend the analysis to account for density-dependent feedbacks?  Answers to these questions are likely to provide important insights into the evolution of dispersal and metapopulation persistence.

% \bibliography{linearSDE}

\begin{thebibliography}{75}
\providecommand{\natexlab}[1]{#1}
\providecommand{\url}[1]{\texttt{#1}}
\expandafter\ifx\csname urlstyle\endcsname\relax
  \providecommand{\doi}[1]{doi: #1}\else
  \providecommand{\doi}{doi: \begingroup \urlstyle{rm}\Url}\fi

\bibitem[Adler(1992)]{adler-92}
F.R. Adler.
\newblock The effects of averaging on the basic reproduction ratio.
\newblock \emph{Mathematical Biosciences}, 111:\penalty0 89--98, 1992.

\bibitem[Bascompte et~al.(2002)Bascompte, Possingham, and
  Roughgarden]{bascompte-etal-02}
J.~Bascompte, H.~Possingham, and J.~Roughgarden.
\newblock Patchy populations in stochastic environments: Critical number of
  patches for persistence.
\newblock \emph{American Naturalist}, 159:\penalty0 128?--137, 2002.

\bibitem[Bena\"{i}m and Schreiber(2009)]{tpb-09}
M.~Bena\"{i}m and S.~J. Schreiber.
\newblock Persistence of structured populations in random environments.
\newblock \emph{Theoretical Population Biology}, 76:\penalty0 19--34, 2009.

\bibitem[Bhattacharya(1978)]{bhattacharya1978criteria}
R.~N. Bhattacharya.
\newblock Criteria for recurrence and existence of invariant measures for
  multidimensional diffusions.
\newblock \emph{The Annals of Probability}, 6\penalty0 (4):\penalty0 pp.
  541--553, 1978.
\newblock ISSN 00911798.
\newblock URL \url{http://www.jstor.org/stable/2243121}.

\bibitem[Bogachev et~al.(2002)Bogachev, R\"ockner, and
  Stannat]{bogachev2002uniqueness}
V~I Bogachev, M~R\"ockner, and W~Stannat.
\newblock Uniqueness of solutions of elliptic equations and uniqueness of
  invariant measures of diffusions.
\newblock \emph{Sbornik: Mathematics}, 193\penalty0 (7):\penalty0 945, 2002.
\newblock URL \url{http://stacks.iop.org/1064-5616/193/i=7/a=A01}.

\bibitem[Bogachev et~al.(2009)Bogachev, Krylov, and
  R\"ockner]{bogachev2009elliptic}
Vladimir~I Bogachev, Nikolai~V Krylov, and Michael R\"ockner.
\newblock Elliptic and parabolic equations for measures.
\newblock \emph{Russian Mathematical Surveys}, 64\penalty0 (6):\penalty0 973,
  2009.
\newblock URL \url{http://stacks.iop.org/0036-0279/64/i=6/a=R02}.

\bibitem[Boyce et~al.(2006)Boyce, Haridas, Lee, and the NCEAS Stochastic
  Demography Working~Group]{boyce-etal-06}
M.S. Boyce, C.V. Haridas, C.T. Lee, and the NCEAS Stochastic Demography
  Working~Group.
\newblock Demography in an increasingly variable world.
\newblock \emph{Trends in Ecology \& Evolution}, 21:\penalty0 141 -- 148, 2006.

\bibitem[Cantrell and Cosner(1991)]{cantrell-cosner-91}
R.~S. Cantrell and C.~Cosner.
\newblock The effects of spatial heterogeneity in population dynamics.
\newblock \emph{Journal of Mathematical Biology}, 29:\penalty0 315--338, 1991.

\bibitem[Cantrell et~al.(2006)Cantrell, Cosner, and Lou]{cantrell-etal-06}
R.~S. Cantrell, C.~Cosner, and Y.~Lou.
\newblock Movement toward better environments and the evolution of rapid
  diffusion.
\newblock \emph{Mathematical Biosciences}, 204\penalty0 (2):\penalty0 199--214,
  2006.

\bibitem[Cantrell and Cosner(1989)]{cantrell-cosner-89}
R.S. Cantrell and C.~Cosner.
\newblock Diffusive logistic equations with indefinite weights: population
  models in disrupted environments.
\newblock \emph{Proceedings of the Royal Society of Edinburgh. Section A.
  Mathematics}, 112\penalty0 (3-4):\penalty0 293--318, 1989.

\bibitem[Cantrell et~al.(2007)Cantrell, Cosner, Deangelis, and
  Padron]{cantrell-etal-07}
R.S. Cantrell, C.~Cosner, D.~L. Deangelis, and V.~Padron.
\newblock The ideal free distribution as an evolutionarily stable strategy.
\newblock \emph{Journal of Biological Dynamics}, 1:\penalty0 249--271, 2007.

\bibitem[Chesson(2000)]{chesson-00}
P.L. Chesson.
\newblock General theory of competitive coexistence in spatially-varying
  environments.
\newblock \emph{Theoretical Population Biology}, 58:\penalty0 211--237, 2000.

\bibitem[Da~Prato and Zabczyk(1996)]{MR1417491}
G.~Da~Prato and J.~Zabczyk.
\newblock \emph{Ergodicity for infinite-dimensional systems}, volume 229 of
  \emph{London Mathematical Society Lecture Note Series}.
\newblock Cambridge University Press, Cambridge, 1996.

\bibitem[DeAngelis et~al.(2011)DeAngelis, Wolkowicz, Lou, Jiang, Novak,
  Svanb{\"a}ck, Ara{\'u}jo, Jo, and Cleary]{deangelis-etal-11}
D.L. DeAngelis, G.S.K. Wolkowicz, Y.~Lou, Y.~Jiang, M.~Novak, R.~Svanb{\"a}ck,
  M.S. Ara{\'u}jo, Y.S. Jo, and E.A. Cleary.
\newblock The effect of travel loss on evolutionarily stable distributions of
  populations in space.
\newblock \emph{American Naturalist}, 178:\penalty0 15--29, 2011.

\bibitem[Delibes et~al.(2001)Delibes, Gaona, and P.]{delibes-etal-01}
M.~Delibes, P.~Gaona, and Ferreras P.
\newblock Effects of an attractive sink leading into maladaptive habitat
  selection.
\newblock \emph{American Naturalist}, 158:\penalty0 277--285, 2001.

\bibitem[Dennis et~al.(1991)Dennis, Munholland, and Scott]{dennis-etal-91}
B.~Dennis, P.L. Munholland, and J.M. Scott.
\newblock Estimation of growth and extinction parameters for endangered
  species.
\newblock \emph{Ecological monographs}, 61:\penalty0 115--143, 1991.

\bibitem[Diaconis(1988)]{diaconis-88}
P.~Diaconis.
\newblock \emph{Group representations in probability and statistics}.
\newblock Institute of Mathematical Statistics Lecture Notes---Monograph
  Series, 11. Institute of Mathematical Statistics, Hayward, CA, 1988.

\bibitem[Diamond(1975)]{diamond-75}
J.M. Diamond.
\newblock The island dilemma: lessons of modern biogeographic studies for the
  design of natural reserves.
\newblock \emph{Biological Conservation}, 7:\penalty0 129--146, 1975.

\bibitem[Dias(1996)]{dias-96}
P.C. Dias.
\newblock Sources and sinks in population biology.
\newblock \emph{Trends Ecol. Evol.}, pages 326--330, 1996.

\bibitem[Dockery et~al.(1998)Dockery, Hutson, Mischaikow, and
  Pernarowski]{dockery-etal-98}
J.~Dockery, V.~Hutson, K.~Mischaikow, and M.~Pernarowski.
\newblock The evolution of slow dispersal rates: a reaction diffusion model.
\newblock \emph{Journal of Mathematical Biology}, 37:\penalty0 61--83, 1998.

\bibitem[Durrett and Remenik(in press)]{durrett-remenik-11}
R.~Durrett and D.~Remenik.
\newblock Evolution of the dispersal distance.
\newblock \emph{Journal of Mathematical Biology}, in press.

\bibitem[Fahrig(1997)]{fahrig-97}
L.~Fahrig.
\newblock Relative effects of habitat loss and fragmentation on population
  extinction.
\newblock \emph{The Journal of Wildlife Management}, 61:\penalty0 603--610,
  1997.

\bibitem[Fahrig(2002)]{fahrig-02}
L.~Fahrig.
\newblock Effect of habitat fragmentation on the extinction threshold: a
  synthesis.
\newblock \emph{Ecological Applications}, 12:\penalty0 346--353, 2002.

\bibitem[Foley(1994)]{foley-94}
P.~Foley.
\newblock Predicting extinction times from environmental stochasticity and
  carrying capacity.
\newblock \emph{Conservation Biology}, pages 124--137, 1994.

\bibitem[Fretwell and Lucas(1970)]{fretwell-lucas-70}
S.D. Fretwell and H.L.~Jr. Lucas.
\newblock On territorial behavior and other factors influencing habitat
  distribution in birds.
\newblock \emph{Acta Biotheoretica}, 19:\penalty0 16--36, 1970.

\bibitem[Gardiner(2004)]{gardiner-04}
C.W. Gardiner.
\newblock \emph{Handbook of stochastic methods: for physics, chemistry \& the
  natural sciences}, volume~13 of \emph{Series in synergetics}.
\newblock Springer, 4th edition, 2004.

\bibitem[Gei{\ss} and Manthey(1994)]{MR1290705}
C.~Gei{\ss} and R.~Manthey.
\newblock Comparison theorems for stochastic differential equations in finite
  and infinite dimensions.
\newblock \emph{Stochastic Processes and Applications}, 53:\penalty0 23--35,
  1994.

\bibitem[Gilpin(1988)]{gilpin-88}
M.E. Gilpin.
\newblock A comment on quinn and hastings: extinction in subdivided habitats.
\newblock \emph{Conservation Biology}, 2:\penalty0 290--292, 1988.

\bibitem[Gonzalez and Holt(2002)]{gonzalez-holt-02}
A.~Gonzalez and R.D. Holt.
\newblock The inflationary effects of environmental fluctuations in source-sink
  systems.
\newblock \emph{Proceedings of the National Academy of Sciences}, 99:\penalty0
  14872--14877, 2002.

\bibitem[Gradshteyn and Ryzhik(2007)]{gradshteyn-07}
I.~S. Gradshteyn and I.~M. Ryzhik.
\newblock \emph{Table of integrals, series, and products}.
\newblock Elsevier/Academic Press, Amsterdam, seventh edition, 2007.
\newblock Translated from the Russian, Translation edited and with a preface by
  Alan Jeffrey and Daniel Zwillinger.

\bibitem[Harper(1982)]{harper-82}
D.G.C. Harper.
\newblock Competitive foraging in mallards: ``ideal free'' ducks.
\newblock \emph{Animal Behaviour}, 30:\penalty0 575--584, 1982.

\bibitem[Harrison and Quinn(1989)]{harrison-quinn-89}
S.~Harrison and J.~F. Quinn.
\newblock Correlated environments and the persistence of metapopulations.
\newblock \emph{Oikos}, 56:\penalty0 293--298, 1989.

\bibitem[Hastings(1983)]{hastings-83}
A.~Hastings.
\newblock Can spatial variation alone lead to selection for dispersal?
\newblock \emph{Theoretical Population Biology}, 24:\penalty0 244--251, 1983.

\bibitem[Holt(1985)]{holt-85}
R.D. Holt.
\newblock Patch dynamics in two-patch environments: {S}ome anomalous
  consequences of an optimal habitat distribtuion.
\newblock \emph{Theor. Pop. Biol.}, 28:\penalty0 181--208, 1985.

\bibitem[Holt(1997)]{holt-97}
R.D. Holt.
\newblock On the evolutionary stability of sink populations.
\newblock \emph{Evolutionary Ecology}, 11:\penalty0 723--731, 1997.

\bibitem[Holt and McPeek(1996)]{holt-mcpeek-96}
R.D. Holt and M.A. McPeek.
\newblock Chaotic population dynamics favors the evolution of dispersal.
\newblock \emph{American Naturalist}, 148\penalty0 (44):\penalty0 709--718,
  1996.

\bibitem[Hutson et~al.(2001)Hutson, Mischaikow, and
  Pol{\'a}{\v{c}}ik]{hutson-etal-03}
V.~Hutson, K.~Mischaikow, and P.~Pol{\'a}{\v{c}}ik.
\newblock The evolution of dispersal rates in a heterogeneous time-periodic
  environment.
\newblock \emph{Journal of Mathematical Biology}, 43:\penalty0 501--533, 2001.

\bibitem[Ikeda and Watanabe(1989)]{MR1011252}
N.~Ikeda and S.~Watanabe.
\newblock \emph{Stochastic differential equations and diffusion processes},
  volume~24 of \emph{North-Holland Mathematical Library}.
\newblock North-Holland Publishing Co., Amsterdam, second edition, 1989.

\bibitem[Jansen and Yoshimura(1998)]{jansen-yoshimura-98}
V.~A.~A. Jansen and J.~Yoshimura.
\newblock Populations can persist in an environment consisting of sink habitats
  only.
\newblock \emph{Proceedings of the National Academy of Sciences USA},
  95:\penalty0 3696--3698, 1998.

\bibitem[Keagy et~al.(2005)Keagy, Schreiber, and Cristol]{re-05}
J.~Keagy, S.~J. Schreiber, and D.~A. Cristol.
\newblock Replacing sources with sinks: When do populations go down the drain?
\newblock \emph{Restoration Ecology}, 13:\penalty0 529--535, 2005.

\bibitem[Khas$'$minskii(1960)]{khasminskii1960ergodic}
R.~Z. Khas$'$minskii.
\newblock Ergodic properties of recurrent diffusion processes and stabilization
  of the solution to the cauchy problem for parabolic equations.
\newblock \emph{Theory of Probability and its Applications}, 5\penalty0
  (2):\penalty0 179--196, 1960.
\newblock ISSN 0040585X.
\newblock \doi{DOI:10.1137/1105016}.
\newblock URL \url{http://dx.doi.org/10.1137/1105016}.

\bibitem[Kirkland et~al.(2006)Kirkland, Li, and Schreiber]{siap-06}
S.~Kirkland, C.K. Li, and S.~J. Schreiber.
\newblock On the evolution of dispersal in patchy landscapes.
\newblock \emph{SIAM Journal on Applied Mathematics}, 66:\penalty0 1366--1382,
  2006.

\bibitem[Kreuzer and Huntly(2003)]{kreuzer-huntly-03}
M.~P. Kreuzer and N.~J. Huntly.
\newblock Habitat-specific demography: evidence for source-sink population
  structure in a mammal, the pika.
\newblock \emph{Oecologia}, 134:\penalty0 343--349, 2003.

\bibitem[Lande et~al.(2003)Lande, Engen, and S{\ae}ther]{lande-etal-03}
R.~Lande, S.~Engen, and B.E. S{\ae}ther.
\newblock Stochastic population dynamics in ecology and conservation: an
  introduction.
\newblock 2003.

\bibitem[Levin et~al.(1984)Levin, Cohen, and Hastings]{levin-etal-84}
S.~A. Levin, D.~Cohen, and A.~Hastings.
\newblock Dispersal strategies in patchy environments.
\newblock \emph{Theoretical Population Biology}, 26:\penalty0 165 -- 191, 1984.

\bibitem[Levins(1969)]{levins-69b}
R.~Levins.
\newblock Some demographic and genetic consequences of environmental
  heterogeneity for biological control.
\newblock \emph{Bulletin of the ESA}, 15:\penalty0 237--240, 1969.

\bibitem[Lonsdale(1993)]{lonsdale-93}
W.~M. Lonsdale.
\newblock Rates of spread of an invading species- {\emph{ {m}imosa pigra}} in
  northern {A}ustralia.
\newblock \emph{Journal of Ecology}, 81:\penalty0 513--521, 1993.

\bibitem[Lundberg et~al.(2000)Lundberg, Ranta, Ripa, and
  Kaitala]{lundberg-etal-00}
P.~Lundberg, E.~Ranta, J.~Ripa, and V.~Kaitala.
\newblock Population variability in space and time.
\newblock \emph{Trends in Ecology and Evolution}, 15:\penalty0 460--464, 2000.

\bibitem[Matthews and Gonzalez(2007)]{matthews-gonzalez-07}
D.~P. Matthews and A.~Gonzalez.
\newblock The inflationary effects of environmental fluctuations ensure the
  persistence of sink metapopulations.
\newblock \emph{Ecology}, 88:\penalty0 2848--2856, 2007.

\bibitem[May(1975)]{may-75}
R.~M. May.
\newblock \emph{Stability and Complexity in Model Ecosystems, 2nd edn.}
\newblock Princeton University Press, Princeton, 1975.

\bibitem[McPeek and Holt(1992)]{mcpeek-holt-92}
M.A. McPeek and R.D. Holt.
\newblock The evolution of dispersal in spatially and temporally varying
  environments.
\newblock \emph{American Naturalist}, 6:\penalty0 1010--1027, 1992.

\bibitem[Metz et~al.(1983)Metz, de~Jong, and Klinkhamer]{metz-etal-83}
J.~A.~J. Metz, T.~J. de~Jong, and P.~G.~L. Klinkhamer.
\newblock What are the advantages of dispersing; a paper by {K}uno extended.
\newblock \emph{Oecologia}, 57:\penalty0 166--169, 1983.

\bibitem[Murphy(2001)]{murphy-01}
M.~T. Murphy.
\newblock Source-sink dynamics of a declining eastern kingbird population and
  the value of sink habitats.
\newblock \emph{Conserv. Biol.}, 15:\penalty0 737--748, 2001.

\bibitem[Oksanen et~al.(1995)Oksanen, Power, and Oksanen]{oksanen-etal-95}
T.~Oksanen, M.E. Power, and L.~Oksanen.
\newblock Ideal free habitat selection and consumer-resource dynamics.
\newblock \emph{American Naturalist}, 146:\penalty0 565--585, 1995.

\bibitem[Petchey et~al.(1997)Petchey, Gonzalez, and Wilson]{petchy-etal-97}
O.~L. Petchey, A.~Gonzalez, and H.~B. Wilson.
\newblock Effects on population persistence: The interaction between
  environmental noise colour, intraspecific competition and space.
\newblock \emph{Proceedings: Biological Sciences}, 264:\penalty0 1841--1847,
  1997.

\bibitem[Philippi and Seger(1989)]{philippi-seger-89}
T.~Philippi and J.~Seger.
\newblock Hedging one's evolutionary bets, revisited.
\newblock \emph{Trends Ecol. Evol.}, 4:\penalty0 41--44, 1989.

\bibitem[Pulliam(1988)]{pulliam-88}
H.~R. Pulliam.
\newblock Sources, sinks, and population regulation.
\newblock \emph{Amer. Nat.}, 132:\penalty0 652--661, 1988.

\bibitem[Quinn and Hastings(1987)]{quinn-hastings-87}
J.F. Quinn and A.~Hastings.
\newblock Extinction in subdivided habitats.
\newblock \emph{Conservation Biology}, 1:\penalty0 198--209, 1987.

\bibitem[Reme\v{s}(2000)]{remes-00}
V.~Reme\v{s}.
\newblock How can maladaptive habitat choice generate source-sink population
  dynamics?
\newblock \emph{Oikos}, 91:\penalty0 579--582, 2000.

\bibitem[Roy et~al.(2005)Roy, Holt, and Barfield]{roy-etal-05}
M.~Roy, R.D. Holt, and M.~Barfield.
\newblock Temporal autocorrelation can enhance the persistence and abundance of
  metapopulations comprised of coupled sinks.
\newblock \emph{American Naturalist}, 166:\penalty0 246--261, 2005.

\bibitem[Ruelle(1979)]{ruelle-79}
D.~Ruelle.
\newblock Analycity properties of the characteristic exponents of random matrix
  products.
\newblock \emph{Adv. in Math.}, 32:\penalty0 68--80, 1979.

\bibitem[Schmidt(2004)]{schmidt-04}
K.~A. Schmidt.
\newblock Site fidelity in temporally correlated environments enhances
  population persistence.
\newblock \emph{Ecology Letters}, 7:\penalty0 176?--184, 2004.

\bibitem[Schreiber(2010)]{prsb-10}
S.~J. Schreiber.
\newblock Interactive effects of temporal correlations, spatial heterogeneity,
  and dispersal on population persistence.
\newblock \emph{Proceedings of the Royal Society: Biological Sciences},
  277:\penalty0 1907--1914, 2010.

\bibitem[Schreiber and Lloyd-Smith(2009)]{amnat-09b}
S.~J. Schreiber and J.~O. Lloyd-Smith.
\newblock Invasion dynamics in spatially heterogenous environments.
\newblock \emph{American Naturalist}, 174:\penalty0 490--505, 2009.

\bibitem[Schreiber and Saltzman(2009)]{amnat-09a}
S.~J. Schreiber and E.~Saltzman.
\newblock Evolution of predator and prey movement into sink habitats.
\newblock \emph{American Naturalist}, 174:\penalty0 68--81, 2009.

\bibitem[Schreiber and Vejdani(2006)]{prsb-06}
S.~J. Schreiber and M.~Vejdani.
\newblock Handling time promotes the coevolution of aggregation in
  predator-prey systems.
\newblock \emph{Proceedings of the Royal Society: Biological Sciences},
  273:\penalty0 185--191, 2006.

\bibitem[Schreiber et~al.(2000)Schreiber, Fox, and Getz]{amnat-00}
S.~J. Schreiber, L.~R. Fox, and W.~M. Getz.
\newblock Coevolution of contrary choices in host-parasitoid systems.
\newblock \emph{American Naturalist}, pages 637--648, 2000.

\bibitem[Serre(1977)]{serre-77}
J.P. Serre.
\newblock \emph{Linear representations of finite groups}.
\newblock Springer-Verlag, New York, 1977.
\newblock Translated from the second French edition by Leonard L. Scott,
  Graduate Texts in Mathematics, Vol. 42.

\bibitem[Slatkin(1974)]{slatkin-74}
M.~Slatkin.
\newblock Hedging one's evolutionary bets.
\newblock \emph{Nature}, 250:\penalty0 704--705, 1974.

\bibitem[Talay(1991)]{talay1991lyapunov}
Denis Talay.
\newblock Approximation of upper {Lyapunov} exponents of bilinear stochastic
  differential systems.
\newblock \emph{SIAM Journal on Numerical Analysis}, 28\penalty0 (4):\penalty0
  1141--1164, 1991.
\newblock ISSN 00361429.
\newblock URL \url{http://www.jstor.org/stable/2157791}.

\bibitem[Tuljapurkar(1990)]{tuljapurkar-90}
S.~Tuljapurkar.
\newblock \emph{Population Dynamics in Variable Environments}.
\newblock Springer-Verlag, New York, 1990.

\bibitem[Turelli(1978)]{turelli-78}
M.~Turelli.
\newblock Random environments and stochastic calculus.
\newblock \emph{Theoretical Population Biology}, 12:\penalty0 140--178, 1978.

\bibitem[van Baalen and Sabelis(1999)]{vanbaalen-sabelis-99}
M.~van Baalen and M.~W. Sabelis.
\newblock Nonequilibrium population dynamics of ``ideal and free'' prey and
  predators.
\newblock \emph{The American Naturalist}, 154:\penalty0 69--88, 1999.

\bibitem[Wilbur and Rudolf(2006)]{wilbur-rudolf-06}
H.~M. Wilbur and V.~H.~W. Rudolf.
\newblock Life-history evolution in uncertain environments: Bet hedging in
  time.
\newblock \emph{The American Naturalist}, 168:\penalty0 398--411, 2006.

\bibitem[Wilcox and Murphy(1985)]{wilcox-murphy-85}
B.A. Wilcox and D.D. Murphy.
\newblock Conservation strategy: the effects of fragmentation on extinction.
\newblock \emph{American Naturalist}, 125:\penalty0 879--887, 1985.

\end{thebibliography}

\appendix 

\section{Proof of Proposition~\ref{prop:frequencies}}
\label{apx:frequencies}

Define the matrix $R$ by 
\[
R := \mathrm{diag}(\mu) + D.
\]
Equation \eqref{eq:main2} becomes
\[ d\bX_t =  \mathrm{diag}( \bX_t) \Gamma^T d\bB_t + R^T \bX_t dt. \]
Recall that $Y_t^j = X_t^j/(X_t^1 + \cdots + X_t^n)$ for each $1 \le j \le n$
and $\bY_t = (Y_t^1, \ldots, Y_t^n)^T$.
Fix $j$ and define $f_j(x_1, \ldots, x_n) := x_j/(x_1 + \cdots + x_n)$, so that $Y^j = f_j(\bX)$. 
Using $\partial_k$ to denote differentiation with respect to $x_k$,
observe that 
\[ \partial_j f_j(x_1, \ldots, x_n) = \left(\sum_{\ell \ne j} x_\ell\right) \bigg / \left(\sum_\ell x_\ell\right)^2, \quad 
\partial_k f_j(x_1, \ldots, x_n) = - x_j \bigg / \left(\sum_\ell x_\ell\right)^2, \; k \ne j.
\]
Moreover,
\[
\partial_{jj} f_j(x_1, \ldots, x_n) = -2 \left(\sum_{\ell \ne j} x_\ell\right) \bigg / \left(\sum_\ell x_\ell\right)^3,\]
\[ \partial_{jk}f_j(x_1, \ldots, x_n) =  - 1 \bigg / \left(\sum_{\ell} x_{\ell}\right)^2 + 2x_j \bigg / \left(\sum_{\ell} x_{\ell}\right)^3, k \ne j
\]
and 
\[
\partial_{k m} f_j(x_1, \ldots, x_n) = 2 x_j \bigg / \left(\sum_\ell x_\ell\right)^3, \; k, m \ne j.
\]

It follows from It\^o's lemma~\citep{gardiner-04} that for each $1 \le j \le n$, 
\begin{eqnarray*}
dY_t^j
& =&
\sum_{k=1}^n  \partial_k f_j(\bX_t)  X^k_t \Gamma_{*k}^T d\bB_t 
 +
\sum_{k =1}^n \partial_k f_j(\bX_t) \bX_t^T  R_{*k} dt \\
&& \qquad {} + 
(1/2) \sum_{k, m=1}^n \partial_{k m}f_j(\bX_t)  X_t^k X_t^m  (\Sigma)_{km} dt,
\end{eqnarray*}
 where $\Gamma_{*k}$ and $R_{*k}$ denote the $k^\mathrm{th}$ columns of the matrices $\Gamma$ and  $R$ respectively. 
Substituting in the derivatives of $f_j$ gives
\[
\begin{split}
dY_t^j
 =
&- \sum_{k \ne j} Y_t^j Y_t^k  \Gamma_{*k}^T d\bB_t 
+  \sum_{k \ne j} Y_t^j Y_t^k  \Gamma_{*j}^T d\bB_t \\
& \qquad {} - \sum_{k \ne j} Y_t^j  Y_t^T R_{*k} dt
+  \sum_{k \ne j} Y_t^k  Y_t^T R_{*j}  dt \\
& \qquad {} + 
(1/2) \sum_{k, m \ne j} 2 Y_t^j  Y_t^k Y_t^m  \Sigma_{km}  dt  - 
(1/2) \sum_{ k  \ne j} 2 Y_t^k (Y_t^j)^2 \Sigma_{jj}  dt \\
&\qquad {} +   (1/2) \times 2 \sum_{k \ne j} \Big( -  Y_t^j  Y_t^k  + 2 Y_t^k (Y_t^j)^2  \Big) \Sigma_{kj} dt \\
 =
&-  Y_t^j  \sum_{k} Y_t^k  \Gamma_{*k}^T d\bB_t 
+  Y_t^j   \Gamma_{*j}^T d\bB_t 
-  Y_t^j  \sum_{k }  \bY_t^T R_{*k} dt
+    \bY_t^T R_{*j}  dt \\
& \qquad {} + 
Y_t^j  \sum_{k, m }  Y_t^k Y_t^m \Sigma_{km} dt   - Y_t^j   \sum_{k }   Y_t^k  \Sigma_{kj} dt.
% =
%&-  Y_t^j  \sum_{k} Y_t^k  d\langle \sigma^{(k)}, B_t \rangle 
%+  Y_t^j   d\langle \sigma^{(j)}, B_t \rangle 
%-  Y_t^j   \langle Y_t, \mu  \rangle dt
%+    \langle Y_t , R^{(j)}  \rangle dt \\
%& + 
%Y_t^j  \sum_{k, m }  Y_t^k Y_t^m  \langle \sigma^{(k)},  \sigma^{(m)}  \rangle dt   - Y_t^j   \sum_{k }   Y_t^k   \langle \sigma^{(k)},  \sigma^{(j)}  \rangle dt.
\end{split}
\]
Since $D \bone = 0$, we have 
$\sum_k R_{*k} = R \bone = \diag(\mu) \bone = \mu$,
and the above system of SDEs can be written in the following compact way
\[ \begin{split}
d\bY_t =&  - \bY_t \bY_t^T \Gamma^T d\bB_t + \mathrm{diag}(\bY_t)  \Gamma^T d\bB_t  \\
&- \bY_t \bY_t^T \mu dt + R^T \bY_t dt  + \bY_t  \bY_t^T \Sigma \bY_t dt  - \mathrm{diag}(\bY_t) \Sigma \bY_t dt\\
&=  \left( \mathrm{diag}(\bY_t) - \bY_t \bY_t^T  \right) \Gamma^T d\bB_t    + D^T \bY_t dt \\
&\quad + \left( \mathrm{diag}(\bY_t) - \bY_t \bY_t^T  \right) \left( \mu - \Sigma \bY_t\right)dt.\\ 
\end{split}\]

Now that the SDE \eqref{eq:freq} is established,
we will prove the ergodicity of the Markov process $(\bY_t)_{ t \ge 0}$ defined in \eqref{eq:freq}.

\bigskip
\noindent
\textbf{Existence.} 
Clearly $(\bY_t)_{ t \ge 0}$ is a Feller process. 
Since for each $t \ge 0$, the random vector $\bY_t$ takes values in the compact state space $\Delta$,  
it trivially follows that the family of probability measures $\{\P^y \{\bY_t \in \cdot\}: t>0\}$ is uniformly tight for any fixed $y \in \Delta$,
where $\P^y$ denotes the law of the process with $\bY_0=y$. 
Hence, by the Krylov-Bogolyubov theorem (see, for example, \citep[Corollary 3.1.2]{MR1417491}), there exists at least one probability measure $\mu$ on $\Delta$ which is an invariant measure for the process $(\bY_t)_{ t \ge 0}$, 
that is,
\[ \int_{ \Delta} \mu(dy) \P^y \{\bY_t \in \cdot\} = \mu\{\cdot\}.\]

\bigskip
\noindent
\textbf{Uniqueness.}
The uniqueness of the invariant measure for $(\bY_t)_{ t \ge 0}$ is ensured by the Doob-Khasminskii theorem  (see, for example, \citep[Chapter 7]{MR1417491}
), 
provided this process satisfies the following two properties:
\begin{itemize}
\item $(\bY_t)_{ t \ge 0}$ is \emph{irreducible}, that is, $\P^y \{\bY_t \in V\} > 0$ for any $t>0$ and any open set $V$ in the simplex $\Delta$.
\item $(\bY_t)_{ t \ge 0}$ is \emph{strong Feller}, that is, $\Delta \ni y \mapsto \int_\Delta \P^y \{\bY_t \in dz\} f(z) $ is continuous for any bounded measurable function $f: \Delta \to \mathbb R$.
\end{itemize}
These conditions also ensure that $(\bY_t)_{ t \ge 0}$ converges in law to the unique invariant measure. 
We next establish irreducibility and the strong Feller property of  $(\bY_t)_{ t \ge 0}$ separately.

\medskip
\noindent
\textbf{(a) Irreducibility.} It clearly suffices to show that the process $(\bX_t)_{ t \ge 0}$ as defined by \eqref{eq:main2} is irreducible,
that is, that $\P^x \{ \bX_t \in U\} > 0$ for each $t > 0$, $x \in \mathbb R_+^n \setminus \{0\}$ and open set $U \subseteq \mathbb R_+^n$.

We will first prove that $\P^x\{ X^i_t > 0 \ \forall i\} =1$ for all $t > 0$ and all $x \in \mathbb R_+^n\setminus \{0\}$,
by induction on the size of the set $G: = \{ 1 \le i \le n: x_i = 0\}$. 
First consider the case $\#G = 0$. By a suitable comparison theorem for SDEs
\citep[Theorem 1.1]{MR1290705},
$\P^x\{ \bX_t \ge \widehat \bX_t \text { for all } t \ge 0 \} =1$, where $\hat X$ is defined by
\[ d\hat X^i_t  = \mu_i \hat X^i_t  dt + \hat X^i_t dE^i_t + D_{ii} \hat X^i_t dt, \quad 1 \le i \le n.  \]
This SDE has the unique solution $\hat X^i_t = x^i \exp( E^i_t + (\mu + D_{ii} - \frac{1}{2}\Sigma_{ii}) t)>0$, so
\begin{equation}\label{eq:strictly_positive}
   \P^x\{ X^i_t > 0 \ \forall i \text{ for all }  t>0 \}, \quad  x \in (0, \infty)^n.  
\end{equation}

Now suppose $\# G = k< n$. By the irreducibility of the
infinitesimal generator matrix $D$, there exist $i_0 \in G, j_0 \not \in G$ such that $D_{j_0, i_0} > 0$.  
Consider the new SDE
\[ d\tilde X^i_t  = \mu_i \tilde X^i_t  dt + \tilde X^i_t dE^i_t + D_{ii} \tilde X^i_t dt, \quad i \ne i_0,  \]
and 
\[ d\tilde X^{i_0}_t  = \mu_{i_0} \tilde X^{i_0}_t  dt +\tilde X^{i_0}_t dE^i_t +  (D_{j_0i_0} \tilde X^{j_0}_t + D_{i_0i_0} \tilde X^{i_0}_t)  dt. \]
By the same comparison theorem, $\P^x\{ \bX_t \ge \widetilde \bX_t \text { for all } t \ge 0 \} =1$. 
Clearly, $\P^x\{\tilde X^i_t > 0\} = 1$ for all $i \not \in G$ and for all $t >0$. 
Since $\tilde X^{i_0}_0 = 0$ and $\tilde X^{j_0}_0>0$,  at time $t= 0$ the diffusion component of $\tilde X_t^{i_0}$ vanishes but its drift coefficient is strictly positive. It follows that $\P^x\{ \tilde X^{i_0}_t > 0\} = 1$ for all $t > 0$. Hence, at any positive time $t$, almost surely $ \widetilde \bX_t$ has at most $k-1$ zero coordinates, 
and, by the comparison theorem, so does $\bX_t$. 
Using the Markov property and the induction hypothesis, we deduce that $\P^x\{ X^i_t > 0 \ \forall i \}=1$ for all $t > 0$. 
This proves that each component of $\bX$ is strictly positive with probability 1 for each $t>0$.

Let $\varphi: (0, \infty)^n \to \mathbb R^n$ be the homeomorphism given by $ \varphi(x) =  (\log x_1, \ldots, \log x_n)$.  Set $\bH_t = \varphi(\bX_t)$, with $\bH_t = (H^1_t, \ldots, H^n_t)^T$. By \eqref{eq:strictly_positive}, this stochastic process is well defined provided $\bX_0 \in (0, \infty)^n$.
Note that  $(\bH_t)_{ t \ge 0}$ satisfies the following SDE,
\[ dH^i_t = (\mu_i - \frac{1}{2} \Sigma_{ii}) dt + dE^i_t + e^{ - H^i_t} \sum_{j=1}^n D_{ji}  e^{H^j_t} dt, \quad 1 \le i \le n.\]
 By Girsanov's theorem 
 (see \citep[Section 4 of Chapter IV]{MR1011252}), the law of $(\Gamma^T)^{-1} \bH_t$ (and hence the law of $\bH_t$) 
is absolutely continuous with respect to the law of $\bB_t$ for any $t>0$. 
Thus, $\P^x\{ \bH_t \in V\} > 0$ for any open set $V \subseteq \mathbb R^n$. 
Finally, for any $x \in \mathbb R^n \setminus \{0\}$,
\begin{eqnarray*}
\P^x\{ \bX_t \in U\} &=& \int_{ \mathbb R_+^n}  \P^x\{ \bX_{t/2} \in dy \}  \P^y\{ \bX_{t/2} \in U\} \\
&= & \int_{ (0, \infty)^n} \P^x\{ \bX_{t/2} \in dy \} \P^y\{ \bX_{t/2} \in U\}\\
&= & \int_{ (0, \infty)^n} \P^x\{ \bX_{t/2} \in dy \} \P^{\varphi(y)}\{ \bH_{t/2} \in \varphi(U)\} >0.
\end{eqnarray*}

\medskip
\noindent
\textbf{(b) Strong Feller property.}  
Note that $\bH$ satisfies a SDE of the form
$d\bH_t = \Gamma^T d\bB_t + b(\bH_t) dt$
for some smooth function $b: \mathbb R^n \to \mathbb R^n$.
For each $K \ge 1$, consider a new SDE
\[ d\bH^K_t = \Gamma^T d\bB_t + b^K(\bH_t)dt,\]
where $b^K:\mathbb R^n \to \mathbb R^n$ is a smooth bounded function with bounded derivative such that $b^K(x) = b(x)$ on $[ - K, K]^n$.
Since the matrix $\Gamma$ is nonsingular, the associated Fisk-Stratonovich type generator of  $(\bH^K_t )_{t\ge 0}$ is trivially hypoelliptic, which in turn implies that  $(\bH^K_t )_{t\ge 0}$ is strong Feller for every $K \ge 1$
%for any $t > 0, bounded measurable $f:\mathbb R^n \to \mathbb R$,  the map $x %\mapsto \E^x f(\bH^K_t)$ is continuous 
(see \citep[Section 8 of Chapter V]{MR1011252}). If we define a sequence of stopping times $\tau_K  := \inf \{ t : \|X_t\|_\infty \ge K\}$, then  $\bH^K_0 = \bH_0 = x \in [-K, K]^n$ implies $\bH^K_t = \bH_t$ for $t \in [0, \tau_K]$. Let $t> 0$ and  $f$ be a bounded measurable function. Fix $\epsilon>0$.  Then for any $x \in \mathbb R^n$,
 \[ \big|\E^x [f(\bH_t)] - \E^x [f(\bH^K_t)]\big |  \le   2 \|f\|_\infty  \P^x \{ \tau_K < t\}.\] 
 Hence, for any open neighborhood $U(x)$ of $x$,
 \[ \big|\E^y [f(\bH_t)] - \E^x [f(\bH_t)]\big |  
 \le   \big|\E^y [f(\bH^K_t)] - \E^x [f(\bH^K_t)]\big | + 4 \|f\|_\infty  \sup_{ z\in U(x)} \P^x \{ \tau_K < t\} \quad \text{ for all } y \in U(x).\]
 
 Since almost surely $\tau_K \uparrow \infty$, we can choose $K$ large enough such that $\P^x \{ \tau_K < t\} < \epsilon(8 \|f\|_\infty)^{-1}$. 
Moreover, by the Feller property of $(\bH_t)_{ t \ge 0}$, there exists a neighborhood $U^1(x)$ of $x$  
such that $\sup_{z \in U^1(x)} \P^z \{ \tau_K < t\} < \epsilon(8 \|f\|_\infty)^{-1}$. 
From the strong Feller property of $(\bH_t^K)_{ t \ge 0}$, 
there exists a neighborhood $U^2(x)$ of $x$ such that  $\big|\E^y [f(\bH^K_t)] - \E^x [f(\bH^K_t)]\big | < \epsilon/2$ for all $y \in U^2(x)$. Thus, 
$\big|\E^y [f(\bH_t)] - \E^x [f(\bH_t)]\big | < \epsilon$ for all $y \in  U^1(x) \cap U^2(x)$. 
Hence, $x \mapsto \E^x [f(\bH_t)]$ is continuous.
Now, for  $t> 0$ and a bounded measurable function $g: \mathbb R_+^n \to \mathbb R$,
\[ \E^x [g(\bX_t)] = \int_{ (0, \infty)^n} \P^x\{\bX_{t/2} \in dy\} \E^{\varphi(y)} [g (\varphi^{-1}( \bH_{t/2}))], \quad x \in \mathbb R_+^n. \]
Therefore, the map $x \mapsto \E^x [g(\bX_t)]$ is continuous, and so $(\bX_t)_{ t \ge 0}$ is a strong Feller process.
It follows easily  that $(\bY_t)_{ t \ge 0}$ is also a strong Feller process.
\qed

\section{Proof of Proposition~\ref{prop:limit}}
\label{apx:limit}

By rescaling time $\tau := \delta t$ and setting $\epsilon := 1/\delta$, \eqref{eq:freq} becomes 
\begin{equation}\label{eq:rescaled}
d\bY_\tau^\epsilon =\sqrt{\epsilon} f(\bY_\tau^\epsilon)d\bB_\tau + \epsilon g(\bY_\tau^\epsilon)dt + Q^T \bY_\tau^\epsilon dt
\end{equation}
where 
$
f(y) :=\left( \mathrm{diag}(y) - y y^T  \right) \Gamma^T
$,
$g(y) := \left( \mathrm{diag}(y) - y y^T  \right) \left( \mu - \Sigma y\right)
$,
and
$\bY_\tau^\epsilon := \bY_{\tau/\epsilon}$.

 For $\epsilon>0$, let $\nu_\epsilon$ be the unique invariant probability measure for \eqref{eq:rescaled} guaranteed by Proposition~\ref{prop:frequencies}. The irreducibility of $Q$ implies that
$\pi$ is the unique stable point for the ODE
\[
\frac{d}{d\tau} y_\tau^x = Q^T y_\tau^x, \quad y_0^x = x \in \Delta,
\]
and that $\lim_{\tau \to \infty} y_\tau^x = \pi$ for any $x \in \Delta$.
Write $\nu_0$
for the Dirac measure at the point $\pi \in \Delta$.   
By the compactness of Borel probability measures on $\Delta$
in the topology of weak convergence, it suffices to show if $\nu_{\epsilon_k}$ converges weakly to $\nu$ for some sequence $\epsilon_k \downarrow 0$, then $\nu=\nu_0$, and hence it is sufficient to check that 
\[
\int_\Delta  h(y_\tau^x) \, \nu(dx) = \int_\Delta h(x) \, \nu(dx) 
\]
for every $\tau \ge 0$ and Lipschitz function $h:\Delta \to \R$.

Set $\bY^k_\tau=\bY^{\epsilon_k}_\tau$ 
and $\nu_k=\nu_{\epsilon_k}$ for ease of notation.
Let $L$ be the Lipschitz constant for the function $h$. Then, 
\begin{eqnarray*}
\left| \int_\Delta \left(  h(y_\tau^x) - h(x) \right) \, \nu(dx) \right| 
 &=& \lim_{k\to\infty} \left| \int_\Delta \left(  h(y_\tau^x) -h(x) \right) \, \nu_k(dx)\right| \\
 &\le & \limsup_{k\to\infty}  \underbrace{\left| \int_\Delta \left( \E^x\left[ h(\bY_\tau^k)\right]-h(x) \right) \, \nu_k(dx)\right|}_{=0 \mbox{ by invariance of $\nu_k$}}\\
 &&+\limsup_{k \to \infty} \left| \int_\Delta \E^x\left[ h(y_\tau^x) - h(\bY_\tau^k)\right] \,\nu_k(dx)\right|\\
 &\le& \limsup_{k\to\infty} L\, \int_\Delta \E^x \left[ \| y_\tau^x - \bY_\tau^k\| \right] \, \nu_k(dx),
\end{eqnarray*}
where $\| \cdot \|$ is the usual Euclidean norm on $\mathbb{R}^n$.

It remains to show that 
$ \lim_{k\to\infty}  \sup_{x \in \Delta} \E^x \left[ \| y_\tau^x - \bY_\tau^k\| \right] = 0$. 
Fix $x \in \Delta$ and set $Z_\tau^k:= y_\tau^x - \bY_\tau^k$.
By It\^o's formula, 
\begin{eqnarray*}
\E^x \left[ \|\bZ^k_\tau\|^2 \right] &=&  \E \left[\int_0^\tau 2\langle \bZ_s^k, Q^T \bZ_s^k \rangle - 2\epsilon_k \langle \bZ_s^k, g(\bY_s^k)\rangle+\epsilon_k \Tr(f(\bY^k_s) f(\bY^k_s)^T) ds \right]\\
&\le &  2\|Q^T\| \int_0^\tau \E^x \left[\| \bZ_s^k\|^2\right]ds + \epsilon_k C\tau,
\end{eqnarray*} 
for some constant $C$ that does not depend on $x$ or $\tau$, 
where we write $\langle \cdot, \cdot \rangle$
for the usual Euclidean inner product
on $\mathbb{R}^n$, and 
$\|Q^T\| = \sup_{\|z\| = 1} |\langle z, Q^T z \rangle|$. 
Gronwall's inequality implies that
\[
\E^x \left[ \|Z_\tau^k\|^2 \right] \le \epsilon_k C e^{2\|Q^T\| \tau },
\] 
and so, by Jensen's inequality,
\[
\E^x \left [ \|Z_\tau^k \|\right] \le \sqrt{\epsilon_k C} e^{\|Q^T\| \tau}.
\]
It follows that $ \lim_{k\to\infty}  \sup_{x \in \Delta} \E^x \left[ \| y_\tau^x - \bY_\tau^k\| \right] = 0$, and hence $\nu = \nu_0$, as required.

In particular,
\[
\begin{split}
\chi(\delta) 
& = 
\int_\Delta \mu^T y \, \nu_{1/\delta}(dy)
-
\frac{1}{2} \int_\Delta  y^T \Sigma y \, \nu_{1/\delta}(dy) \\
& \to \mu^T \pi - \frac{1}{2} \pi^T \Sigma \pi \\
\end{split}
\]
as $\delta \to \infty$. \qed

\section{Proof of Proposition~\ref{prop:cont}}
\label{apx:continuity}

Fix $\delta \in [0, \infty)$, and denote our underlying
probability space by $(\Omega,\mathcal{F},\P)$.
Define
\[ 
\Phi^\delta_{s, t} : \mathbb R^n \times \Omega \to \mathbb R^n, 
\quad  0 \le s \le t,
\] 
by
$\Phi^\delta_{s,t}(\bx, \omega) = \bX^\delta_t(\omega)$, where $(\bX^\delta_u)_{ u \ge s}$ is the unique solution of
\[ 
\bX^\delta_u =  \bx + \int_s^u \mathrm{diag}( \bX^\delta_v) \Gamma^T d\bB_v + \int_s^u  (R_\delta)^T \bX^\delta_v dv
\]
with $R_\delta := \mathrm{diag}(\mu) + \delta Q$.

Note that for all $0 \le s \le w \le t$, 
\begin{equation}\label{eq:flow}
\Phi_{s, t}(\cdot, \omega) =  \Phi_{w, t}^\delta(\cdot, \omega)\circ \Phi_{s, w}^\delta(\cdot, \omega).
\end{equation}
It is easy to see that $ \Phi^\delta_{s, t}( \cdot, \omega)$ is a linear map from $\mathbb R^n$ to $\mathbb R^n$ and thus can be represented by a matrix $\bM_{s,t}^\delta(\omega)$. From \eqref{eq:flow}, it follows that
\[ \bM_{s, t}^\delta(\omega) =  \bM_{w, t}^\delta(\omega) \bM_{s, w}^\delta( \omega) \quad \text{for all } 0 \le s \le w \le t .
 \]

Since $\bM^\delta_{s,t}$ is constructed from $ (\bB_u-\bB_s)_{ u \in [s, t]}$, 
the matrices $\{\bM^\delta_{k, k+1}\}_{ k \in \mathbb N }$ are independent. Moreover, since the drift and the diffusion coefficients do not depend on time,  $\{\bM^\delta_{k, k+1}\}_{ k \in \mathbb N }$  is a stationary sequence.

We note that the  Lyapunov exponent $\chi(\delta)$ of $(\bX^\delta_t)_{t \ge 0}$ is the same as 
\[
\lim_{k \to \infty } \E \left[ k^{-1} \log \| \bM^\delta_{0, k}\| \right] 
= 
\inf_{k \ge 1} \E \left[ k^{-1} \log \| \bM^\delta_{0, k}\| \right],
\]
where we set 
\[
\|A\| := \sup\left\{\sum_{i,j} A_{ij} x_j : \sum_k x_k = 1, \, x_k \ge 0 \, \forall k\right\}
\]
for a matrix $A$ with nonnegative entries.

Set $\mathbb R^n_+ := \{ \bx \in \mathbb R^n : x \ge 0\}$.  If $\delta>0$, then it follows from the irreducibility of $Q$ that 
\begin{equation} \label{eq:ruelle}
\bM_{s,t}^\delta (\mathbb R^n_+) \subseteq  \{ x \in \mathbb R^n : x_i > 0 \; \mbox{for all} 1 \le i \le n \} \cup \{0\}
\end{equation}
and hence $\chi(\delta)$ is analytic on $(0, \infty)$ by \citep[Theorem~3.1]{ruelle-79}.

The condition \eqref{eq:ruelle} fails to hold when $\delta=0$
and so we must proceed differently.
We first claim that for fixed $t>0$ the map 
$\delta \mapsto t^{-1} \E [\log \| \bM^\delta_{0, t}\|]$ is upper semicontinuous on $[0, \infty)$. 
To see this, fix $\delta \in [0, \infty)$. 
Set $\log^+ x = \max(0,\log x)$ and $\log^- x = \min(0,\log x)$.
It follows from the continuous dependence of the solution of a SDE
on its parameters \citep[4.3.2]{gardiner-04}, 
that $\bX^{\delta'}_t \to \bX^{\delta}_t$ almost surely as $\delta' \to \delta$, 
which implies that $ \| \bM^{\delta'}_{0, t}\| \to \| \bM^\delta_{0, t}\|$ almost surely as $\delta' \to \delta$.
An application of Gronwall's lemma gives that
$\E[\sup_{0 \le  \delta \le c} \| \bX_t^\delta\|] < \infty$ for each $c>0$.
Hence,
\[ \E \left[ \log^+ \| \bM^{\delta'}_{0, t}\| \right] \to \E \left[ \log^+ \| \bM^{\delta}_{0, t}\| \right] \mbox{ as }\delta' \to \delta.\]
% \marginnote{Should we work out the Gronwall's lemma application here?}
On the other hand, by Fatou's lemma,
\[ \E \left[ -  \log^- \| \bM^{\delta}_{0, t}\| \right] \le \liminf_{\delta' \to \delta} \E \left[- \log^- \| \bM^{\delta'}_{0, t}\|\right] .\]
Combining these two inequalities gives 
\[  \limsup_{\delta' \to \delta} \E \left[ \log \| \bM^{\delta'}_{0, t}\| \right] \le \E \left[ \log \| \bM^{\delta}_{0, t}\| \right],\] 
and the claim follows.

Since $\chi(\delta) = \inf_{ t > 0} t^{-1}  \E \log \| \bM^\delta_{0, t}\|$ 
is the infimum of a family of upper semicontinuous functions, 
it is itself upper semicontinuous, or equivalently,
$\limsup_{\delta' \to \delta} \chi(\delta') \le \chi(\delta)$. In particular, $\limsup_{\delta \to 0} \chi(\delta)\le \chi(0)$.

We now prove the opposite inequality that $\liminf_{\delta \to 0} \chi(\delta)\ge \chi(0)$.
Fix $\delta > 0$, and without loss of generality suppose that $\max_i -Q_{ii}=1$, 
so that if $x_i \ge z_i \ge 0$ for $1 \le i \le n$, 
then $(Qx)_i \ge - z_i$ for   $1 \le i \le n$.
% Note that for  $\bx, \by \in \mathbb R^n_+$ such that  $y_i \ge x_i, y_j = x_j,$
% \[ \mu_j y_j+  \sum_{ i} \delta Q_{ij} y_i  = (\mu_j - \delta) y_j+   \delta \sum_{ i}  P_{ij} y_i  \ge  (\mu_j - \delta) x_j   \]
Consider the two SDEs
\[ d\bX^\delta_t =  \mathrm{diag}( \bX^\delta_t) \Gamma^T d\bB_t +  (\mathrm{diag}(\mu)+ \delta Q^T ) \bX^\delta_t dt\]
and 
\[ d\bZ^\delta_t =  \mathrm{diag}( \bZ^\delta_t) \Gamma^T d\bB_t +  \mathrm{diag}(\mu -\delta)  \bZ^\delta_t dt.
\]
If $\bX^\delta_0 = \bZ^\delta_0$, then, by the comparison theorem,  
\[ \bX^\delta_t  \ge \bZ^\delta_t  \quad \text{for all } t \ge 0\] 
almost surely.

Thus, the  Lyapunov exponent of $(\bX^\delta_t)_{t \ge 0}$ dominates that of $ (\bZ^\delta_t)_{t \ge 0}$. 
Note that the coordinates of $\bZ^\delta$ are decoupled 
and hence the  Lyapunov exponent of this process
is the maximum of the stochastic growth rates for 
the individual coordinate processes. Therefore,
\[ \chi(\delta) \ge  \max_{j} \left(\mu_j - \frac12 \sum_k \sigma_{kj}^2\right) - \delta.    \]
In particular,
\begin{equation}\label{eq:lbound}
\liminf_{\delta \to 0+} \chi(\delta)  \ge  \max_{j} \left(\mu_j - \frac12 \sum_k \sigma_{kj}^2\right) = \chi(0),
\end{equation}
as required.
\qed

\section{Proof of Theorem~\ref{thm:bigone}}
\label{apx:bigone}

Recall that 
\[
d\bY_t =  
\left( \mathrm{diag}(\bY_t) - \bY_t \bY_t^T  \right) \Gamma^T d\bB_t + D^T \bY_t dt 
+ \left( \mathrm{diag}(\bY_t) - \bY_t \bY_t^T  \right) \left( \mu - \Sigma \bY_t\right)dt,
\]
where $D$ is of the form $\delta Q$, with $Q$ an irreducible 
infinitesimal generator matrix and $\delta > 0$. 
Moreover, $Q$ is assumed to be reversible with respect to the
unique probability vector $\pi$ satisfying $Q^T \pi = 0$;  
that is, that $\pi_i Q_{ij} = \pi_j Q_{ji}$
for all $i,j$.  

Define an inner product on $\mathbb{R}^n$
by $\langle u, v \rangle_\pi := \sum_i \frac{1}{\pi_i} u_i v_i = u^T \mathrm{diag}(\pi)^{-1} v$. 
It follows from reversibility that 
the linear operator $v \mapsto Q^T v$ is self-adjoint
with respect to this inner product;
that is, that 
$\langle u, Q^T v \rangle_\pi = \langle Q^T u, v \rangle_\pi$ for
all $u,v$.

From the spectral theorem and the Perron-Frobenius theorem, the linear operator
$v \mapsto Q^T v$ has eigenvalues 
$\lambda_1 \le \lambda_2 \le \ldots \le \lambda_{n-1} < \lambda_n = 0$
and corresponding orthonormal eigenvectors $\xi_1, \ldots, \xi_n$ with $\xi_n = \pi$
such that
\[
Q^T v = \sum_{k=1}^{n-1} \lambda_k \xi_k \langle v, \xi_k \rangle_\pi, 
\quad v \in \mathbb{R}^n.
\]
Note that
\begin{equation}
\label{elliptic}
\bone^T v = 
\langle v, \pi \rangle_\pi = 0 
\Longrightarrow 
\langle v, Q^T v \rangle_\pi \le - \kappa \|v\|_\pi^2,
\end{equation}
where $\kappa := - \lambda_{n-1} > 0$ and
$\|\cdot \|_\pi$ is the norm
associated with the inner product $\langle \cdot, \cdot \rangle_\pi$.

Note also that if $\bone^T v = 0$, then
\[
w := \sum_{k=1}^{n-1} \lambda_k^{-1} \xi_k \langle v, \xi_k \rangle_\pi
\]
is the unique vector with the properties
\[
\langle w, \pi \rangle_\pi = 0 \quad \text{and} \quad Q^T w = v.
\]
In particular,
\[
\bone^T \left( \mathrm{diag}(\pi) - \pi \pi^T  \right) \left( \mu - \Sigma \pi\right)
=
\left(\pi^T - \pi^T  \right) \left( \mu - \Sigma \pi\right)
=
0,
\]
and so there is a unique vector we denote $\nu$ such that
\begin{equation}
\label{eq:nu}
\bone^T \nu = \langle \nu, \pi \rangle_\pi = 0
\quad \text{and} \quad
Q^T \nu = - \left( \mathrm{diag}(\pi) - \pi \pi^T  \right) \left( \mu - \Sigma \pi\right).
\end{equation}
We emphasize that $\nu$ does not depend on $\delta$.

Consider the stochastic process
\[
\bU_t := \delta^{\frac{1}{2}}\left(\bY_{t/\delta} - \pi - \delta^{-1} \nu\right),
\]
so that
\[
\bY_t = \delta^{-\frac{1}{2}} \bU_{\delta t} + \pi + \delta^{-1} \nu.
\]
Observe
that $\pi + \delta^{-1} \nu$ is indeed a probability vector for
$\delta$ sufficiently large.
Because we are only interested in the equilibrium law of the process $\bY$, we
assume that $\bY_0 = \pi + \delta^{-1} \nu$ 
and hence $\bU_0 = 0$.    
Note that $0 = \bone^T \bU_t 
= \langle \bU_t, \pi \rangle_\pi$ for all $t \ge 0$.

We have for the standard Brownian motion $\tilde \bB_t := \delta^{\frac{1}{2}} \bB_{t/\delta}$ that
\[
\begin{split}
d\bU_t
& = 
\left( 
\mathrm{diag}(\delta^{-\frac{1}{2}} \bU_t + \pi + \delta^{-1} \nu) 
- (\delta^{-\frac{1}{2}} \bU_t + \pi + \delta^{-1} \nu) (\delta^{-\frac{1}{2}} \bU_t + \pi + \delta^{-1} \nu)^T  \right) 
\Gamma^T \, d \tilde \bB_t \\
& \quad +
\delta^{-\frac{1}{2}} 
\delta Q^T (\delta^{-\frac{1}{2}} \bU_t + \pi + \delta^{-1} \nu) \, dt \\
& \quad +
\delta^{-\frac{1}{2}}
\left( \mathrm{diag}(\delta^{-\frac{1}{2}} \bU_t + \pi + \delta^{-1} \nu) 
- (\delta^{-\frac{1}{2}} \bU_t + \pi + \delta^{-1} \nu)
(\delta^{-\frac{1}{2}} \bU_t + \pi + \delta^{-1} \nu)^T  \right) \\
& \qquad \times \left( \mu - \Sigma (\delta^{-\frac{1}{2}} \bU_t + \pi + \delta^{-1} \nu)\right) \, dt. \\
\end{split}
\]
Using $Q^T \pi = 0$ and \eqref{eq:nu}, we get
\[
\begin{split}
d\bU_t
& =
\left[ \mathrm{diag}(\pi) - \pi \pi ^T  \right] \Gamma^T \, d \tilde \bB_t 
+
Q^T \bU_t \, dt \\
& \quad +
\left[
\delta^{-\frac{1}{2}} A_{\frac{1}{2}}(\bU_t) 
+ \delta^{-1} A_{1}(\bU_t) 
+ \delta^{-\frac{3}{2}} A_{\frac{3}{2}}(\bU_t) 
+ \delta^{-2} A_{2}(\bU_t) 
\right]
\, d \tilde \bB_t \\
& \quad +
\left[
\delta^{-1} b_{1}(\bU_t) 
+ \delta^{-\frac{3}{2}} b_{\frac{3}{2}}(\bU_t) 
+ \delta^{-2} b_{2}(\bU_t) 
+ \delta^{-\frac{5}{2}} b_{\frac{5}{2}}(\bU_t)
+ \delta^{-3} b_{3}(\bU_t) 
+ \delta^{-\frac{7}{2}} b_{\frac{7}{2}}(\bU_t) 
\right]
\, dt, \\
\end{split}
\] 
where
\[
A_{\frac{1}{2}}(u)
:=
\left[\mathrm{diag}(u) - u \pi^T - \pi u^T\right] \Gamma^T
\]
\[
A_{1}(u)
:=
\left[-u u^T + \mathrm{diag}(\nu) - \pi \nu^T - \nu \pi^T\right] \Gamma^T 
\]
\[
A_{\frac{3}{2}}(u)
:=
\left[-u \nu^T - \nu u^T\right] \Gamma^T 
\]
\[
A_{2}(u)
:=
- \nu \nu^T \Gamma^T 
\]
and
\[
\begin{split}
b_{1}(u)
& :=
- \pi  u^T \mu 
- u \pi^T \mu
+ \pi u^T \Sigma \pi
+ u \pi^T \Sigma \pi
+ \pi \pi^T \Sigma u \\
& \quad + \mathrm{diag}(u) \mu 
-\mathrm{diag}(\pi) \Sigma u 
-\mathrm{diag}(u)  \Sigma  \pi \\
\end{split}
\]
\[
\begin{split}
b_{\frac{3}{2}}(u)
& :=
-u u^T \mu
- \pi \nu^T \mu
- \nu \pi^T \mu \\
& \quad
+ \pi u^T \Sigma u
+ u \pi^T \Sigma u
+ u u^T \Sigma \pi
+ \pi \nu^T \Sigma \pi
+ \nu \pi^T \Sigma \pi
+ \pi \pi^T \Sigma \nu \\
& \quad -\mathrm{diag}(\pi) \Sigma \nu
-\mathrm{diag}(u) \Sigma u
+\mathrm{diag}(\nu) \mu
-\mathrm{diag}(\nu) \Sigma \pi \\
\end{split}
\]
\[
\begin{split}
b_{2}(u)
& :=
-u \nu^T \mu
-\nu u^T \mu \\
& \quad
+u u^T \Sigma u
+u \pi^T \Sigma \nu
+u \nu^T \Sigma \pi
+\pi u^T \Sigma \nu
+\pi \nu^T \Sigma u
+\nu u^T \Sigma \pi
+\nu \pi^T \Sigma u \\
& \quad -\mathrm{diag}(u) \Sigma \nu
-\mathrm{diag}(\nu) \Sigma u \\
\end{split}
\]
\[  
\begin{split}
b_{\frac{5}{2}}(u)
& :=
-\nu \nu^T \mu
+u u^T \Sigma \nu
+u \nu^T \Sigma u
+\nu u^T \Sigma u 
+\pi \nu^T \Sigma \nu
+\nu \pi^T \Sigma \nu
+\nu \nu^T \Sigma \pi \\
& \quad
-\mathrm{diag}(\nu) \Sigma \nu \\
\end{split}
\]
\[
b_{3}(u)
:=
u \nu^T \Sigma  \nu
+\nu u^T \Sigma  \nu
+\nu \nu^T \Sigma  u
\]
\[
b_{\frac{7}{2}}(u)
:=
\nu \nu^T \Sigma  \nu.
\]

%By It\^o's lemma,
%\[
%\begin{split}
%d \|\bU_t\|_\pi^2
%& =
%2 \bU_t^T \mathrm{diag}(\pi)^{-1} \left[ \mathrm{diag}(\pi) - \pi \pi ^T  \right] \Gamma^T \, d \tilde \bB_t 
%+
%2 \langle \bU_t, Q^T \bU_t\rangle_\pi \, dt \\
%& \quad +
%2 \sum_{\ell = 1}^4 \delta^{-\frac{\ell}{2}} \bU_t^T \mathrm{diag}(\pi)^{-1} A_{\frac{\ell}{2}}(\bU_t) \, d \tilde \bB_t \\
%& \quad +
%2 \sum_{\ell=2}^7 \delta^{-\frac{\ell}{2}} \bU_t^T \mathrm{diag}(\pi)^{-1} b_{\frac{\ell}{2}}(\bU_t) \, dt \\
%& \quad +
%\bone^T 
%\mathrm{diag}(\pi)^{-1} 
%\left[ \mathrm{diag}(\pi) - \pi \pi ^T  \right] 
%\Gamma^T
%\Gamma 
%\left[ \mathrm{diag}(\pi) - \pi \pi ^T  \right] 
%\mathrm{diag}(\pi)^{-1} 
%\bone \, dt \\
%& \quad +
%\left(
%\sum_{\ell = 1}^4 \delta^{-\frac{\ell}{2}} 
%\bone^T \mathrm{diag}(\pi)^{-1} A_{\frac{\ell}{2}}(\bU_t)
%\right)
%\left(
%\sum_{\ell = 1}^4 \delta^{-\frac{\ell}{2}} 
%A_{\frac{\ell}{2}}(\bU_t)^T \mathrm{diag}(\pi)^{-1} \bone
%\right) \, dt. \\
%\end{split}
%\]

By It\^o's lemma,
\[
\begin{split}
d \|\bU_t\|_\pi^2
& =
2 \bU_t^T \mathrm{diag}(\pi)^{-1} \left[ \mathrm{diag}(\pi) - \pi \pi ^T  \right] \Gamma^T \, d \tilde \bB_t 
+
2 \langle \bU_t, Q^T \bU_t\rangle_\pi \, dt \\
& \quad +
2 \sum_{\ell = 1}^4 \delta^{-\frac{\ell}{2}} \bU_t^T \mathrm{diag}(\pi)^{-1} A_{\frac{\ell}{2}}(\bU_t) \, d \tilde \bB_t \\
& \quad +
2 \sum_{\ell=2}^7 \delta^{-\frac{\ell}{2}} \bU_t^T \mathrm{diag}(\pi)^{-1} b_{\frac{\ell}{2}}(\bU_t) \, dt \\
& \quad +
\Tr \Big(
\mathrm{diag}(\pi)^{-1} 
\left[ \mathrm{diag}(\pi) - \pi \pi ^T  \right] 
\Gamma^T
\Gamma 
\left[ \mathrm{diag}(\pi) - \pi \pi ^T  \right] \Big) dt \\
& \quad +
\Tr \left(  \mathrm{diag}(\pi)^{-1}
\sum_{\ell = 1}^4 \delta^{-\frac{\ell}{2}} 
 A_{\frac{\ell}{2}}(\bU_t) \times 
\sum_{\ell = 1}^4 \delta^{-\frac{\ell}{2}} 
A_{\frac{\ell}{2}}(\bU_t)^T
\right) \, dt. \\
\end{split}
\]

%Observe that
%\[
%\bone^T 
%\mathrm{diag}(\pi)^{-1} 
%\left[ \mathrm{diag}(\pi) - \pi \pi^T  \right] 
%\Gamma^T
%\Gamma 
%\left[ \mathrm{diag}(\pi) - \pi \pi^T  \right] 
%\mathrm{diag}(\pi)^{-1} 
%\bone
%=
%\left[\bone^T - n \pi^T  \right] 
%\Sigma
%\left[\bone - n \pi \right].
%\]
Note also that 
\begin{equation}
\label{apriori}
|U_t^i| \le C \delta^{\frac{1}{2}} \quad 1 \le i \le n,
\end{equation}
for an appropriate
constant $C$ because $0 \le Y_t^i \le 1$, $1 \le i \le n$.  
Each function 
\[
u \mapsto u^T \mathrm{diag}(\pi)^{-1} b_{\frac{\ell}{2}}(u),
\quad 2 \le \ell \le 7,
\]
is a polynomial in $u$ with total degree at most $\ell$,
and each function
%\[
%u \mapsto
%\left(
%\bone^T \mathrm{diag}(\pi)^{-1} A_{\frac{\ell'}{2}}(u)
%\right)
%\left(
%A_{\frac{\ell''}{2}}(u)^T \mathrm{diag}(\pi)^{-1} \bone
%\right),
%\quad
%1 \le \ell', \ell'' \le 4,
%\]
\[
u \mapsto
\Tr \left(
 \mathrm{diag}(\pi)^{-1} A_{\frac{\ell'}{2}}(u) A_{\frac{\ell''}{2}}(u)^T 
\right),
\quad
1 \le \ell', \ell'' \le 4,
\]
is a polynomial in $u$ with total degree at most $\ell'+\ell''$.

It follows that
\begin{equation}
\label{second_moment_U_differential_bound}
\frac{d}{dt} \E\left[\|\bU_t\|_\pi^2\right] 
\le 
-2 \kappa \E\left[\|\bU_t\|_\pi^2\right] + C'
\end{equation}
for all $t \ge 0$ for a suitable constant $C'$ that does
not depend on $\delta$. Hence, 
\begin{equation}
\label{second_moment_U}
\sup_{t \ge 0} \E\left[\|\bU_t\|_\pi^2\right] \le \frac{C'}{2\kappa}
\end{equation}
(recall that $\bU_0 = 0$).

Let $(\bV_t)_{t \ge 0}$ be the solution of the stochastic differential equation
\[
d \bV_t
=
\left[ \mathrm{diag}(\pi) - \pi \pi ^T  \right] \Gamma^T \, d \tilde \bB_t 
+
Q^T \bV_t \, dt \\
\]
with $\bV_0 = \bU_0 = 0$. Note that 
$d (\bone^T \bV_t) = 0$ for all $t \ge 0$,
and so $\langle \bV_t, \pi \rangle_\pi = \bone^T \bV_t = 0$ for all $t \ge 0$.
It is readily checked that
\[
\bV_t = \int_0^t \exp( Q^T (t-s)) \left[ \mathrm{diag}(\pi) - \pi \pi ^T  \right] \Gamma^T \, d \tilde \bB_s.
\]
So $\bV$ is a Gaussian process for which $\bE[\bV_t]=0$ and
\begin{equation}
\label{eq:V_var_covar}
\E[\bV_t \bV_t ^T] = 
\int_0^t 
\exp(Q^T s) 
\left( \mathrm{diag}(\pi) - \pi \pi^T  \right) 
\Sigma
\left( \mathrm{diag}(\pi) - \pi \pi^T  \right)   
\exp( Q s)   \, ds
\end{equation}
for all $t \ge 0$.
%We know \citep[e.g.][4.4.6]{gardiner-04}\marginnote{Is this reference to 4.4.6 correct? If so for which edition? Also check the next time it is cited. -Sebastian} that $(\bV_t)_{t \ge 0}$ is a centered Gaussian process with
%(e.g., {\em Handbook of Stochastic Methods for Physics, Chemistry and the Natural Sciences} by C. W. Gardiner)
Consequently,
\begin{equation}
\label{all_moments_V}
\sup_{t \ge 0} \E\left[|V_t^i|^p\right] < \infty
\end{equation}
for $1 \le i \le n$ and $p \ge 0$.

In the notation above,
\[
d(\bU_t - \bV_t)
=
Q^T (\bU_t - \bV_t) \, dt
+
\left[\sum_{\ell=1}^4 \delta^{-\frac{\ell}{2}} A_{\frac{\ell}{2}}(\bU_t)\right] \, d\tilde \bB_t
+
\left[\sum_{\ell=2}^7 \delta^{-\frac{\ell}{2}} b_{\frac{3}{2}}(\bU_t)\right] \, dt.
\]
Applying It\^o's lemma and a combination of
\eqref{apriori}, \eqref{second_moment_U} and \eqref{all_moments_V},
we can argue along the lines we followed to establish 
\eqref{second_moment_U_differential_bound} to see that
\[
\frac{d}{dt} \E\left[\|\bU_t - \bV_t\|_\pi^2\right] 
\le 
-2 \kappa \E\left[\|\bU_t - \bV_t\|_\pi^2\right] + \delta^{-1} C''
\]
for all $t \ge 0$ for a suitable constant $C''$ that does
not depend on $\delta$. Hence,
\begin{equation}
\label{compare_U_V}
\sup_{t \ge 0} \E\left[\|\bU_t - \bV_t\|_\pi^2\right] \le \delta^{-1} \frac{C''}{2 \kappa}.
\end{equation}

Now let $\bY_\infty$, $\bU_\infty$ and $\bV_\infty$ be random vectors that are distributed
according to the equilibrium laws of $(\bY_t)_{t \ge 0}$,
$(\bU_t)_{t \ge 0}$
and $(\bV_t)_{t \ge 0}$, respectively.  Also  let $\hat U^i$ and $\hat V^i$ be the $i$-th component of the vectors $\bU_\infty$ and $\bV_\infty$ respectively.

% Thus, $\bY_\infty$ has the same distribution as the random vector we called $\bY_\infty$ (and also $\bZ$) in the paper.  
% {\tt This new notation may be better.}

From \eqref{apriori}, \eqref{second_moment_U} and the linearity of the
function $b_1$,
\[
0 
= 
Q^T \E[\bU_\infty] + \delta^{-1} b_1(\E[\bU_\infty]) 
+ \mathrm{O}(\delta^{-\frac{3}{2}}).
\]
Noting that $\langle \E[\bU_\infty], \pi \rangle_\pi = 0$ because
$\langle \bU_t, \pi \rangle_\pi = 0$ for all $t \ge 0$, we have
from \eqref{elliptic} that
\[
\begin{split}
\kappa \|\E[\bU_\infty]\|_\pi^2 
& \le
- \langle \E[\bU_\infty], Q^T \E[\bU_\infty] \rangle_\pi \\
& = \delta^{-1} \langle \E[\bU_\infty], b_1(\E[\bU_\infty]) \rangle_\pi
+ \mathrm{O}(\delta^{-\frac{3}{2}}) \\
& \le
C^{'''} \delta^{-1} \|\E[\bU_\infty]\|_\pi^2 + \mathrm{O}(\delta^{-\frac{3}{2}}) \\
\end{split}
\]
for a suitable constant $C^{'''}$, and hence,
 \begin{equation}
\label{first_moment_U}
\E[\hat U^i] = \mathrm{O}(\delta^{-\frac{3}{4}}), \quad 1 \le i \le n.
\end{equation}

From \eqref{second_moment_U}, \eqref{all_moments_V}  and \eqref{compare_U_V},
\begin{equation}
\label{second_U_V}
\left| 
\E \left[\hat U^i \hat U^j \right] 
- \E \left[\hat V^i \hat V^j \right] 
\right|
=
\mathrm{O}(\delta^{-\frac{1}{2}}), \quad 1 \le i,j \le n.
\end{equation}

Recall that $\chi(\delta)$ is the Lyapunov exponent,
and that
\[
\begin{split}
\chi(\delta)
& =
\mu^T \E\left[\bY_\infty \right]  - 
\frac{1}{2}\E \left[\bY_\infty^T \Sigma \bY_\infty \right] \\
& =
\mu^T \E\left[\delta^{-\frac{1}{2}} \bU_\infty + \pi + \delta^{-1} \nu \right] \\
& \quad - 
\frac{1}{2}
\E \left[
\left(\delta^{-\frac{1}{2}} \bU_\infty + \pi + \delta^{-1} \nu\right)^T 
\Sigma
\left(\delta^{-\frac{1}{2}} \bU_\infty + \pi + \delta^{-1} \nu\right)\right] \\
& =
\delta^{-\frac{1}{2}} \mu^T \E\left[ \bU_\infty\right]
+
\mu^T \left(\pi + \delta^{-1} \nu \right) \\
& \quad -
\delta^{-1}
\frac{1}{2}
\E \left[\bU_\infty^T \Sigma \bU_\infty\right]
- 
2 \delta^{-\frac{1}{2}} \frac{1}{2} \E \left[\bU_\infty^T \right] \Sigma \left(\pi + \delta^{-1} \nu \right) \\
& \quad -
\frac{1}{2} \left(\pi + \delta^{-1} \nu \right)^T \Sigma \left(\pi + \delta^{-1} \nu \right). \\
\end{split}
\]
Substituting in \eqref{first_moment_U} and \eqref{second_U_V}, and noting 
from \eqref{eq:V_var_covar} that
the random vector $\bV_\infty$ is Gaussian with mean vector $0$ and covariance matrix
\[
\int_0^\infty  
\exp(Q^Ts) 
\left( \mathrm{diag}(\pi) - \pi \pi^T  \right) 
\Sigma
\left( \mathrm{diag}(\pi) - \pi \pi^T  \right)   
\exp( Q s)   \, ds ,
\]
we conclude that
\[
\begin{split}
\chi(\delta)
& =
\left( \mu^T \pi -  
\frac{1}{2}  \pi^T \Sigma \pi \right) \\
& \quad +   
\delta^{-1} 
\Big [ 
( \mu  - \Sigma \pi)^T \nu 
- \frac{1}{2} \mathrm{Tr} \left(\E[\bV_\infty \bV_\infty^T]  \Sigma \right) \Big ] \\
& \quad +  
\mathrm{O}(\delta^{-\frac{5}{4}}) \\
& =
\left( \mu^T \pi -  
\frac{1}{2}  \pi^T \Sigma \pi \right) \\
& \quad +  
\delta^{-1} 
\biggl [ 
( \mu  - \Sigma \pi)^T \nu \\
& \qquad -   \frac{1}{2} \int_0^\infty  \mathrm{Tr}(   \exp(Q^Ts) \left( \mathrm{diag}(\pi) - \pi \pi^T  \right) \Sigma \left( \mathrm{diag}(\pi) - \pi \pi^T  \right)   \exp( Q s) \Sigma)  \, ds \biggr ] \\
& \quad + 
\mathrm{O}(\delta^{-\frac{5}{4}}) \\
\end{split}
\]
as $\delta \to \infty$. \qed

\section{Proof of Corollary~\ref{cor:diagonalized}}
\label{apx:diagonalization}

We now assume that the matrices $Q$ and $\Sigma$ are both real symmetric 
($\Sigma$ is, of course, always symmetric) and that they commute.
Hence, as noted in the statement of the corollary, if 
$\lambda_1 \le \ldots \le \lambda_{n-1} < \lambda_n=0$ 
are the eigenvalues of $Q$ with corresponding orthonormal eigenvectors $\xi_1, \ldots, \xi_n$, where
$\xi_n = \frac{1}{\sqrt{n}} \mathbf 1$, then
\[
Q = \sum_{k=1}^n \lambda_k \xi_k \xi_k^T
\]
and it is possible to write the eigenvalues $\theta_1, \ldots, \theta_n$
of $\Sigma$ in some order so that
\[
\Sigma = \sum_{k=1}^n \theta_k \xi_k \xi_k^T.
\]

By the assumption that $Q$ is
symmetric, $\pi = \frac{1}{n} \mathbf 1 = \frac{1}{\sqrt{n}} \xi_n$.
Therefore,
\[
\mu^T \pi -  \frac{1}{2}  \pi^T \Sigma \pi
=
\bar \mu - \frac{1}{2n} \theta_n
\]
where $\bar\mu =\frac{1}{n}\sum_i \mu_i$.
 
To find the unique vector $\nu$ that solves
\[
{\mathbf 1}^T \nu = 0
\quad \text{and} \quad
Q^T \nu = - \left( \mathrm{diag}(\pi) - \pi \pi^T  \right) \left( \mu - \Sigma \pi\right),
\]
write $\nu = \sum_{k=1}^n a_k \xi_k$.  The condition 
${\mathbf 1}^T \nu = 0$ dictates that $a_n = 0$. The second condition becomes
\[
\begin{split}
\sum_{k=1}^{n-1} a_k \lambda_k \xi_k
  & =
    - \frac{1}{n} 
    \left(I - \xi_n \xi_n^T \right)
    \left(\mu - \frac{1}{\sqrt{n}} \theta_n \xi_n \right) \\
  & =
    -\frac{1}{n} 
    \left(\sum_{k=1}^{n-1} \xi_k \xi_k^T \right)
    \left(\mu - \frac{1}{\sqrt{n}} \theta_n \xi_n \right) \\
  & =
    -\frac{1}{n} 
    \sum_{k=1}^{n-1} (\xi_k^T \mu) \xi_k, \\
\end{split}
\]
so that $a_k = -(\xi_k^T \mu)/(n \lambda_k)$ for $1 \le k \le n-1$.
It follows that 
\[
\begin{split}
( \mu  - \Sigma \pi)^T \nu
& =
- 
\left(\mu - \frac{1}{\sqrt{n}} \theta_n \xi_n \right)^T 
\left(\sum_{k=1}^{n-1} \frac{\xi_k^T \mu}{n \lambda_k} \xi_k \right)\\
& =
- 
\sum_{k=1}^{n-1} \frac{(\xi_k^T \mu)^2}{n \lambda_k}. \\
\end{split}
\]
%Note that this term is always positive because the eigenvalues $\lambda_k$, $1 \le k \le n-1$, are
%negative.

Lastly, the matrices inside the trace in the integral
\[
\int_0^\infty  \mathrm{Tr}
\left(   \exp(Q^Ts) \left( \mathrm{diag}(\pi) - \pi \pi^T  \right) \Sigma \left( \mathrm{diag}(\pi) - \pi \pi^T  \right)   \exp( Q s) \Sigma
\right)  \, ds
\]
commute and so the integral is
\[
\begin{split}
& \int_0^\infty  \mathrm{Tr}
\left(
\left( \mathrm{diag}(\pi) - \pi \pi^T  \right)^2 
\Sigma^2   \exp(2 Q s)
\right)  \, ds \\
& \quad =
\frac{1}{n^2} 
\int_0^\infty  \mathrm{Tr}
\left(
\left(I - \xi_n \xi_n^T \right) 
\left(\sum_{k=1}^n \theta_k^2 \xi_k \xi_k^T\right)
\left(\sum_{k=1}^n \exp(2 s \lambda_k) \xi_k \xi_k^T \right)  
\right) \, ds \\
& \quad =
\frac{1}{n^2} 
\int_0^\infty  \mathrm{Tr}
\left(
\sum_{k=1}^{n-1} \theta_k^2  \exp(2 s \lambda_k) \xi_k \xi_k^T
\right)  \, ds \\
& \quad =
\frac{1}{n^2} 
\int_0^\infty  
\left(\sum_{k=1}^{n-1} \theta_k^2  \exp(2 s \lambda_k)\right)  \, ds \\
& \quad =
-\frac{1}{n^2} \sum_{k=1}^{n-1} \frac{\theta_k^2}{2 \lambda_k}. \\
\end{split}
\]
%Note that this term is also always positive.

Therefore, our asymptotic approximation of $\chi(\delta)$ is
\[
\left(\bar \mu - \frac{1}{2n} \theta_n\right)
-
\frac{1}{\delta}
\left[
  \sum_{k=1}^{n-1}
    \frac{1}{n\lambda_k}
    \left(
    (\xi_k^T \mu)^2
    -
    \frac{1}{4 n} \theta_k^2
  \right)
\right]
+ O(\delta^{-5/4})
\]
as $\delta \to 0$. \qed

\section{Proof of Theorem~\ref{thm:character_rep}}
\label{apx:character_rep}

To show that Theorem \ref{thm:character_rep} follows from Corollary \ref{cor:diagonalized},
we show that the matrix entries of each irreducible representation belong to a common eigenspace of $Q$ and $\Sigma$.
Suppose that $c$ is a class function and the matrix $C$ is given by $C_{g,h} = c(gh^{-1})$.
Recall from \eqref{eq:expansion_class_function} that
\[
c(g) 
= 
\frac{1}{\# G} 
\sum_{\kappa \in \tilde G} 
\tilde c(\kappa)
\kappa(g)^*.
\]
Therefore,
\[
C_{g,h} 
= 
\frac{1}{\# G} 
\sum_{\kappa \in \tilde G} 
\tilde c(\kappa)
\kappa(g h^{-1})^*.
\]
If $\kappa$ is associated with 
the irreducible representation $\rho \in \hat G$, then
\[
\kappa(g h^{-1}) 
= \Tr(\rho(g h^{-1}))
= \Tr(\rho(g) \rho(h)^\dag)
= \sum_{i,j = 1}^{d_\rho} \rho_{ij}(g) \rho_{ij}(h)^*
=: (\Xi(\kappa))_{gh},
\]
where $\dag$ denotes the Hermitian conjugate of a matrix.
Set $\Pi_\kappa := (d_\kappa/\#G) \Xi(\kappa)$.
The $\#G \times \#G$ matrix $\Pi_\kappa$ is Hermitian, 
and it follows from \eqref{eq:matrix_entry_orthog}
that 
$\Pi_\kappa^2 = \Pi_\kappa$,
so that $\Pi_\kappa$ is the projection onto
a $d_\kappa^2$-dimensional subspace.  
Again by \eqref{eq:matrix_entry_orthog}, the matrices
$\Pi_{\kappa'}$ and $\Pi_{\kappa''}$ are orthogonal for
distinct $\kappa', \kappa''$.  Thus,
\[
C 
= 
\sum_{\kappa \in \tilde G} 
\frac{\tilde c(\kappa)}{d_\kappa}
\Pi_\kappa.
\]
This expression is nothing other than the spectral decomposition
of the matrix $C$.  
It shows that $\tilde c(\kappa) / d_\kappa$ is an
eigenvalue of $C$ with multiplicity $d_\kappa^2$.
In summary, for each $\kappa \in \tilde G$ there are eigenvalues $\tilde{q}(\kappa)/d_\kappa$ of $Q$ and 
$\tilde s(\kappa)/d_\kappa$ of $\Sigma$, 
each with multiplicity $d_\kappa^2$.

Therefore, in the notation of Corollary~\ref{cor:diagonalized},
\[
\sum_{k=1}^{n-1} 
\frac{\theta_k^2}{\lambda_k}
=
\sum_{\kappa \ne \trchar} 
d_\kappa^2 
\left(\frac{\tilde s(\kappa)}{d_\kappa}\right)^2
\frac{d_\kappa}{\tilde q(\kappa)}
=
\sum_{\kappa \ne \trchar} 
d_\kappa 
\frac{\tilde s(\kappa)^2} {\tilde q(\kappa)}.
\]
Similarly, we can split the sum
\[
\sum_{k=1}^{n-1}
      \frac{1}{\lambda_k}
      (\xi_k^T \mu)^2
\]
up into contributions from each non-trivial character
$\kappa$ that are of the form
\[
\frac{d_\kappa}{\tilde q(\kappa)}
\sum_k (\xi_k^T \mu)^2,
\]
where the sum is over the indices that correspond to
eigenvectors in the range of the projection $\Pi_\kappa$.
By pairwise orthogonality of the matrices $\Pi_\kappa$ and the fact the $\mu$ is real,
this last quantity is equal to
\begin{align*}
\frac{d_\kappa}{\tilde q(\kappa)}
\left \| \Pi_\kappa \mu \right \|^2
 &=
\frac{d_\kappa}{\tilde q(\kappa)}
\left( \frac{d_\kappa}{\# G} \right)
\sum_{g,h\in G} \mu(g) \kappa(gh^{-1}) \mu(h) \\
 &=
\frac{d_\kappa}{\tilde q(\kappa)} \|\mu\|_\kappa^2,
\end{align*}
by definition of $\|\mu\|_\kappa$.
\qed

%
%Finally, note that if $\kappa$ is the character of $\rho$, then
%\begin{align*}
%  \|\mu\|_\kappa^2 &= \sum_{i,j=1}^{d_\rho} \left| \sum_{g \in G} \xi^\rho_{ij}(g) \mu(g) \right|^2 \\
%      &= \frac{\# G}{d_\rho} \sum_{i,j=1}^{d_\rho} \ft \mu(\rho)_{ij} \ft \mu(\rho)_{ij}^* \\
%      &= \frac{\# G}{d_\rho} \Tr\left( \ft \mu (\rho) \ft \mu(\rho)^\conj \right)
%\end{align*}
%where $A^\conj$ denotes the Hermitian conjugate of $A$.
%\qed

\section{Proof of Theorem~\ref{thm:multiscale}}
\label{apx:multiscale}

We first recall some notation. For $0 \le r, \ell \le k+1$, 
\[
Z_r = G_1 \otimes \cdots \otimes G_{r-1}
\otimes \{\id_r\} \otimes \cdots \otimes \{\id_k\},
\]
\[
\bar Z_\ell = \{\id_1\} \otimes \cdots \otimes \{\id_\ell\} 
\otimes G_{\ell+1} \otimes \cdots \otimes G_{k}
\]
and
\[
    \ell(g) := \min\{ j : g_j \neq \id_j \}.
\]
The displacement associated with $g \in G$ moves between two
patches that are in the same metapatch at scale $\ell(g)$
but different metapatches at scales 
$\ell(g)+1, \ell(g)+2, \ldots$
%We also denote the location of the rightmost non-identity element in $g \in G$ by
%\[
%    r(g) := \max\{ j : g_j \neq \id_j \} .
%\]
%Note that $1 \le \ell(g) \le k+1$, $\ell(\id_G) = k+1$,
%$0 \le r(g) \le k$, and $r(\id_G) = 0$.
Recall also that $\# G_r = n_r$, 
$N_r = \#Z_r  = \prod_{j=1}^{r-1} n_j$
and $\bar N_\ell = \# \bar Z_\ell = \prod_{j={l+1}}^{k} n_j$.
%In particular, $\bar N_{k+1} = N_0 = 0$, 
%$\bar Z_k = Z_1 = \{\id_G\}$,
%and $\bar Z_0 = Z_{k+1} = G$.
%Note that $\{ h \in G : r(h) = r \} = Z_{r+1} \setminus Z_r$,
%and that $\{ g \in G : \ell(g) = \ell\} = \bar Z_{\ell-1} \setminus \bar Z_\ell$.

Writing $\bone_j$ for the trivial character on $G_j$, put
\[
\begin{split}
\tilde{Z_r}
& :=
\tilde G_1 \otimes \cdots \otimes \tilde G_{r-1}
\otimes \{\bone_r\} \otimes \cdots \otimes \{\bone_k\} \\
& =
 \{ \kappa \in \tilde G : \kappa(g) = 1 \; \forall g \in \bar Z_{r-1} \} \\
\end{split}
\] 
and
\[
r(\kappa) := \max\{ j : \kappa \notin \tilde Z_j \}.
\]
%Note that $\tilde {Z_1} = \{ \trchar \}$, and $r(\trchar) = 0$.  
%This can be thought of as the set of characters whose elements are the trivial character from position $r$ onwards,
%and so $\tilde {Z_r}$ for larger $r$ corresponds to ``higher frequencies,'' in some sense.

The following orthogonality property of characters:
\[ \sum_{g \in G} \kappa'(g) \kappa''(g)^*  = \left \{ \begin{array}{cc}
 \# G & \text{ if } \kappa'  = \kappa''  \\
0 & \text{otherwise.}
\end{array}  \right.\]
leads to the relation
\[
    \sum_{ g \in \bar Z_r } \kappa(g) = \begin{cases}
      \bar N_r, \quad & \mbox{if} \; \kappa \in \tilde Z_{r+1}, \\
            0, \quad & \mbox{otherwise.}
        \end{cases}
\]
We denote this quantity, as a function of $\kappa$, by $\bar N_r \delta_{\tilde Z_{r+1}}(\kappa)$.

Define the function $f_\ell: G \to \mathbb{C}$ 
by setting $f_\ell(g)=1$ if $\ell(g) = \ell$ and $f_\ell(g)=0$ otherwise.
Then,
\begin{eqnarray*}
    \tilde f_\ell (\kappa) &= &\sum_{g : \ell(g) = \ell} \kappa(g) \\
            &= &\sum_{g \in \bar Z_{\ell-1}} \kappa(g) - \sum_{g \in \bar Z_\ell} \kappa(g) \\
            &= &\bar N_{\ell-1} \delta_{\tilde Z_\ell}(\kappa) - \bar N_\ell \delta_{\tilde Z_{\ell+1}}(\kappa) .
\end{eqnarray*}

Our assumption that $s(g) = s_{\ell(g)}$
implies that $s(g) = \sum_{\ell=1}^{k+1} s_\ell f_\ell(g)$. 
Since $\kappa \in \tilde Z_\ell$ if and only if $r(\kappa) +1 \le  \ell$, it follows by linearity that
\[
\begin{split}
    \tilde s(\kappa) 
    & = \sum_{\ell=1}^{k+1} s_\ell \left( \bar N_{\ell-1} \delta_{\tilde Z_\ell}(\kappa) - \bar N_\ell \delta_{\tilde Z_{\ell+1}}(\kappa) \right) \\
    & = \sum_{\ell= r(\kappa) + 1}^{k+1} s_\ell \bar N_{\ell-1} 
    - \sum_{\ell= r(\kappa)}^{k+1} s_\ell \bar N_{\ell} \\
    & = \sum_{\ell= r(\kappa)}^{k} s_{\ell+1} \bar N_{\ell}
    - \sum_{\ell= r(\kappa)}^{k} s_{\ell} \bar N_{\ell} \\
    & = \sum_{\ell= r(\kappa)}^{k} (s_{\ell+1} - s_{\ell}) \bar N_{\ell}, \\
\end{split}
\]
where we used the convention $\bar N_{k+1} = 0$.

%     \\
%            &= \sum_{\ell=r(\kappa)}^{k+1} s_\ell \left( \bar N_{\ell-1} - \bar N_\ell \right) - s_{r(\kappa)} \bar N_{r(\kappa)-1} \\ 
%            &= \sum_{\ell=r(\kappa)}^k (s_{\ell+1} - s_\ell) \bar N_\ell .
%\end{align*}

Turning to $q$, we have $q(g) = q_{\ell(g)}$ for $g \neq \id_G$
and $q(\id_G) = q_{k+1} = - \sum_{\ell=1}^k q_\ell (\bar N_{\ell-1} - \bar N_\ell)$.
%Thus, since $\bar N_{k+1} = 0$ and $\bar N_k = 1$,
%\begin{align*}
%    \tilde q(\kappa) &= \sum_{\ell=r(\kappa)}^{k+1} q_\ell (\bar N_{\ell-1} - \bar N_\ell) - q_{r(\kappa)} \bar N_{r(\kappa)-1} \\
%            &= \sum_{\ell=r(\kappa)}^{k} q_\ell (\bar N_{\ell-1} - \bar N_\ell) + q_{k+1} - q_{r(\kappa)} \bar N_{r(\kappa)-1} \\
%            &= - \sum_{\ell=1}^{r(\kappa)-1} q_\ell (\bar N_{\ell-1} - \bar N_\ell) - q_{r(\kappa)} \bar N_{r(\kappa)-1}. \\
%\end{align*}
By the same argument as above,
\[
\begin{split}
    \tilde q(\kappa) 
    & = \sum_{\ell= r(\kappa) + 1}^{k+1} q_\ell \bar N_{\ell-1} 
    - \sum_{\ell= r(\kappa)}^{k+1} q_\ell \bar N_{\ell} \\
    & = \sum_{\ell= r(\kappa) + 1}^{k+1} q_\ell (\bar N_{\ell-1} - \bar N_\ell) - q_{r(k)} \bar N_{r(k)} \\
    & = \sum_{\ell= r(\kappa) + 1}^{k} q_\ell (\bar N_{\ell-1} - \bar N_\ell) 
    - \sum_{\ell=1}^k q_\ell (\bar N_{\ell-1} - \bar N_\ell)
    - q_{r(k)} \bar N_{r(k)} \\
    & = - \sum_{\ell= 1}^{r(\kappa)} q_\ell (\bar N_{\ell-1} - \bar N_\ell)
    - q_{r(k)} \bar N_{r(k)} \\
    &= - \sum_{\ell= 1}^{r(\kappa)-1} q_\ell (\bar N_{\ell-1} - \bar N_\ell)
    - q_{r(k)} \bar N_{r(k)-1}. \\
\end{split}
\]

%    & = \sum_{\ell= r(\kappa)}^{k} q_{\ell+1} \bar N_{\ell}
%    - \sum_{\ell= r(\kappa)}^{k} q_{\ell} \bar N_{\ell} \\
%    & = \sum_{\ell= r(\kappa)}^{k-1} q_{\ell+1} \bar N_{\ell}
%    - \sum_{\ell=1}^k q_\ell (\bar N_{\ell-1} - \bar N_\ell)
%    - \sum_{\ell= r(\kappa)}^{k} q_{\ell} \bar N_{\ell} \\

%This is negative and decreasing as $r(\kappa)$ increases, 
%and is (nearly) the total migration rate to all patches differing at level $r$ (including those farther away).
%Furthermore, since $q(g) = q(g^{-1})$, we know that $\tilde q(\kappa)$ is real.

Lastly, for an arbitrary function $\mu$ we need to evaluate 
%$\|\mu\|_\kappa^2 = (d_\rho/\# G) \Tr\left( \hat \mu (\rho) \hat \mu(\rho)^\conj \right)$, 
%where $\kappa$ is the character of the irreducible representation $\rho$.
%This will require that we evaluate sums of the form
\[ \frac{1}{\# G} \sum_{\kappa: r(\kappa)=r} \| \mu \|_\kappa^2. \]
We do that by using the following lemma that follows immediately from orthogonality of characters.

\begin{lemma} \label{lemma:sum_of_squares} 
Let $H$ and $K$ be two finite Abelian groups.
For $f : H \otimes K \to \mathbb{C}$,
\[
\sum_{\kappa \in \tilde H}
\left| \sum_{(h,k) \in H \otimes K} f(h,k) \kappa(h) \right |^2
=
\# H \sum_{h \in H} \left| \sum_{k \in K} f(h,k) \right|^2.
\]
\end{lemma}

%For $u \in \hat G$ and $v \in \hat H$ let $\hat f(u,v) = \sum_{g,h} f(g,h) u(g) v(h)$ be the Fourier transform of $f$.  
%Then, if $\trrep$ is the trivial representation of $H$,
%\[
%    \sum_{u\in \hat G} d_u \Tr\left( \hat f(u,\trrep) \hat f(u,\trrep)^\conj \right) = \# G \sum_{g \in G} \left( \sum_{h \in H} f(g,h) \right)^2 ,
%\]
%where the sum on the left is over irreducible representations of $G$.
%\end{lemma}
%\begin{proof}
%Let $f_0(g) = \sum_{h \in H} f(g,h)$.  Then since $\trrep(h) = 1$ for all $h \in H$,
%\[
%  \hat f_0(u) = \sum_g f_0(g) u(g) = \sum_{g,h} f(g,h) u(g) \trrep(h) = \hat f(u,\trrep) .
%\]
%The lemma follows from Plancherel's relation applied to $f_0$ as a function on $G$.
%\end{proof}

% Theorem~\ref{thm:character_rep} gives us that
% \begin{align}
%     \chi(\delta) &\approx (\bar \mu - \frac{1}{2} \bar s) 
%             - \frac{1}{\delta \prod_j n_j} \sum_{\kappa \neq \trchar} d_\kappa^2 \left( \frac{\|\hat \mu(\kappa)\|^2 - (1/4) \hat s(\kappa)^2}{ \hat q(\kappa)} \right) ,
% \end{align}
% where $\trchar$ is the trivial character of $G$.

Using lemma \ref{lemma:sum_of_squares} applied to the decomposition of $G$ as $Z_r \otimes \bar Z_{r-1}$,
we get 
\[
    \sum_{\kappa \in \tilde Z_r}  \| \mu \|_\kappa^2 = \frac{N_r}{\# G}  \sum_{g \in Z_r} \left( \sum_{z \in \bar Z_{r-1}} \mu(g \gtimes z) \right)^2.
\]
Further decomposing $Z_{r+1}$ as $Z_{r} \otimes G_{r}$
and $\bar Z_{r-1}$ as $\bar Z_{r} \otimes G_{r}$,
and using $N_{r+1} = n_{r} N_{r}$ gives
\begin{eqnarray*}
    \sum_{\kappa : r(\kappa)=r}  \|  \mu \|_\kappa^2 
            &= &\sum_{\kappa \in \tilde Z_{r+1}}  \| \mu \|_\kappa^2 - \sum_{\kappa \in \tilde Z_{r}}  \| \mu \|_\kappa^2  \\
            &= &\frac{N_{r+1}}{\# G} \sum_{g \in Z_{r+1}} \left( \sum_{z \in \bar Z_r} \mu(g \gtimes z) \right)^2 
                - \frac{N_{r}}{\# G} \sum_{g \in Z_{r}} \left( \sum_{z \in \bar Z_{r-1}} \mu(g \gtimes z) \right)^2  \\
            &= &\frac{n_r N_{r+1}}{\# G} \sum_{g \in Z_{r}} \left( 
                    \frac{1}{n_r} \sum_{h \in G_r} \left( \sum_{z \in \bar Z_r} \mu(g \gtimes h \gtimes z) \right)^2  \right. \\
            & &\qquad \qquad   \left.      - \left( \frac{1}{n_r} \sum_{h \in G_r} \sum_{z \in \bar Z_r} \mu(g \gtimes h \gtimes z) \right)^2 \right). 
            % &= \left( \bar N_0^2 \right) \frac{1}{N_r} \sum_{g \in Z_{r}} \left( 
            %         \frac{1}{n_r} \sum_{h \in G_r} \left( \frac{1}{\bar N_r} \sum_{z \in \bar Z_r} \mu(g \gtimes h \gtimes z) \right)^2 \right. \\
            % & \qquad \qquad  \left.      - \left( \frac{1}{n_r} \sum_{h \in G_r} \frac{1}{\bar N_r} \sum_{z \in \bar Z_r} \mu(g \gtimes h \gtimes z) \right)^2 \right) ,
\end{eqnarray*}
To turn the remaining sums into averages, we need to pull out a factor of $N_r \bar N_r^2$,
leaving us with $n_r N_{r+1} N_r \bar N_r^2 = \prod_{\ell=1}^k n_\ell^2 = \# G^2$.
Therefore, recalling that
\[
    v_\mu(r) = \frac{1}{N_r}  \sum_{g \in Z_r} \left(
        \frac{1}{n_r} \sum_{h \in G_r} \left( \frac{1}{\bar N_r} \sum_{z \in \bar Z_r} \mu(g \gtimes h \gtimes z) \right)^2
                    - \left( \frac{1}{n_r} \sum_{h \in G_r} \frac{1}{\bar N_r} \sum_{z \in \bar Z_r} \mu(g \gtimes h \gtimes z) \right)^2 \right),
\]
% (note the relation to ANOVA; see Chapter 8 in Diaconis)
we have 
\[
    \sum_{\kappa : r(\kappa)=r}  \| \mu \|_\kappa^2 = \# G\times v_\mu (r) .
\]

The theorem follows once we note that
\[
\#\{\kappa: r(\kappa)=r\} 
= \# (\tilde Z_{r+1} \setminus \tilde Z_r)
= N_{r+1}-N_r.
\]
\qed

\end{document}